# Notes on Bolyai's *Appendix Demonstrating the Absolute Science of Space*

**Steven Rose**

**University College London**

## Preface

In the following pages I have reprinted the text of Bolyai's *Appendix*, including after each section some explanatory notes. Some of the theorems have been discussed in detail, with a particular focus on those in which Bolyai derives the identities of the hyperbolic plane and on his quadrature of the circle. In these theorems I have tried to fill in some of the details which Bolyai chose to omit, either through lack of space or because he thought them too obvious to include.

The translation of the *Appendix* is taken from an edition supervised by Prof. Ferenc Kárteszi, who also provided clear diagrams and many valuable notes. This edition is published in *János Bolyai*, *Appendix*, *The Theory of Space*, North Holland (1987), ISBN 0444865284. The current publisher, Elsevier, have kindly allowed this translation to be reprinted here.

## Historical Introduction

János Bolyai (1802-1860) was born in Koloszvar (Cluj) and received his early education in Maros-Vásárhely (Târgu-Mures), which were then in the Grand Principality of Transylvania, part of the Austrian Empire, but are now cities in Romania. His father, Farkas (Wolfgang) Bolyai, a talented mathematician in his own right, taught mathematics and science at the Evangelical Reformed College in Maros-Vásárhely where János became a pupil. János was an infant prodigy. At the age of ten he played first violin in the local string quartet. By the age of 13 he had mastered the calculus and even substituted on occasions for his father in class. The students said they liked his lessons better than those of his father. At fifteen he was sent to the Royal College of Engineers in Vienna. After completing the seven year course in four years he was commissioned into the engineering corps of the Austrian army where he served for ten years. Among his many accomplishments he was proficient in several languages: Hungarian, German, Romanian, Italian, French, Latin, Chinese and Tibetan. (It would be interesting to know how he acquired a knowledge of Tibetan in that era). He was also an excellent swordman and gained something of a reputation as a duellist in the army. On one occasion he defeated thirteen cavalry officers one after another, pausing only to play his beloved violin after every other duel. It is not reported whether he injured or killed anyone.

While in the army he devoted much of his spare time to mathematics, in particular to the theory of parallels. On hearing of his son's interest in this field, Farkas wrote:

*'You must not attempt this approach to parallels. I know this way to its very end. I have traversed this bottomless night, which extinguished all light and joy in my life. I entreat you, leave the science of parallels alone.'*

Ignoring his father's warning, Bolyai started out, like so many of his predecessors and indeed like his father before him, with the aim of solving the great geometrical conundrum of the age, which was to prove the parallel postulate. Euclid (c. 300 BC) had prefaced his *Elements* with five postulates:

POSTULATES.

Let the following be postulated:

1. To draw a straight line from any point to any point.
2. To produce a finite straight line continuously in a straight line.
3. To describe a circle with any centre and distance.
4. That all right angles are equal to one another.

The first three postulates underpin the use of compass and straightedge, the only tools permitted in Euclid's geometry. The fourth asserts the homogeneity of the plane, in that right angles, defined as equal adjacent angles, have equal magnitude wherever they are drawn in the plane.

His fifth postulate, however, known as the parallel postulate, was rather more controversial:

5. That, if a straight line falling on two straight lines make the interior angles on the same side less than two right angles, the two straight lines, if produced indefinitely, meet on that side on which are the angles less than the two right angles.

Later mathematicians were generally unwilling to accept Euclid's view of this statement as an assumption, possibly because it is less self-evident than the first four postulates, and preferred to treat it instead as an unproved theorem. So at various times over the next two thousand years mathematicians attempted to prove the parallel postulate, either in its original form or in the form of a logically equivalent statement such as Playfair's axiom. Playfair's axiom states that given a straight line and a point not on the line, only one straight line can be drawn through the point parallel to the given line. All these attempts were flawed, usually because they included an assumption which turned out to be equivalent to the parallel postulate, which had unexpectedly appeared in a new guise, rendering the argument circular.

An important development occurred in the eighteenth century when Gerolamo Saccheri (1667-1733), a Jesuit priest and Professor of Philosophy, Theology and Mathematics at the University of Pavia, investigated a quadrilateral whose base angles are right angles and whose sides about the base are equal, now called a Saccheri quadrilateral. He easily proved that the summit angles are equal. But are they right angles or are they acute or obtuse? Saccheri showed that the parallel postulate is equivalent to the proposition that the summit angles are right angles, which he called the hypothesis of the right angle (HRA). He called the proposition that the summit angles are obtuse the hypothesis of the obtuse angle (HOA) and the proposition that they are acute the hypothesis of the acute angle (HAA). Saccheri also showed that the HRA, the HOA and the HAA are equivalent, respectively, to the proposition that the sum of the angles of a triangle is equal to, greater than or less than two right angles. This meant that any attempt to prove the parallel postulate by contradiction would involve refuting not one but two alternatives. Saccheri, a devoted Euclidean as the title of his book *Euclid Vindicated* makes clear, succeeded in refuting the HOA. He believed that he had also disproved the HAA, which he claimed to be 'absolutely false, because it is repugnant to the nature of the straight line'. But his argument was unconvincing in that it depended on a questionable assumption about straight lines when infinitely extended. In the course of his argument he actually came close to discovering non-Euclidean geometry, but his passionate commitment to Euclidean geometry probably prevented him from pursuing the logic of his own discovery.

The Swiss mathematician and philosopher, Johann Lambert (1726-1777) was wary of the pitfalls facing anyone attempting to prove the parallel postulate:

*'Proofs of the Euclidean postulate can be developed to such an extent that apparently a mere trifle remains. But a careful analysis shows that in this seeming trifle lies the crux of the matter; usually it contains either the proposition that is being proved or a postulate equivalent to it.'*

Lambert was aware of Saccheri's work at second hand and in his (unpublished) *Theory of Parallel Lines* (1776) he undertook a careful examination of the HRA, the HOA and the HAA, discussing the logical implications of each one in turn. He was able to refute the HOA, though by a different method than Saccheri's, but he recognised the difficulty in disposing of the HAA:

*'But for the most part I looked for such consequences of the third hypothesis* (the HAA) *in order to see if contradiction might not come to light. From all this I saw that it is no easy matter to refute this hypothesis…*

Instead he considered what the implications of the HAA might be if it proved to be true. The HAA implies that the fourth angle of a quadrilateral with three right angles is acute. Lambert showed that the sides of this quadrilateral, now called a Lambert quadrilateral, are determined by the magnitude of its acute angle. That is to say, given a Lambert quadrilateral DABG with right angles at D, A and B and where AB = AD, the acute angle at G can fit no other quadrilaterals except those whose sides AB, AD have the absolute lengths AB, AD of the given figure. In other words, a particular side length corresponds to a given angle. Thus the HAA implies the existence of an absolute standard of distance (and, by extension, of area and volume) based on angular measure. This of course would be very convenient since it avoids the need to invent artificial units of measurement. On the other hand Lambert recognised that the HAA would also have serious disadvantages. For in another theorem he had shown that the HAA implies that the defect of a triangle, the amount by which its angle sum falls short of two right angles, is proportional to its area. This in turn implies that similar triangles do not exist. Unfortunately this would undermine the whole of plane trigonometry, which is based on Euclid's theorem (*Elements*, 6, 4) that the corresponding sides of similar triangles are always in the same ratio, which in turn allows a definite trigonometric ratio to be assigned to each angle:

… '*it would result in countless inconveniences. Trigonometric tables would have to be infinitely extended; the similarity and proportionality of figures would entirely lapse; no figure could be presented except in its absolute size; astronomy would become an evil task and so on.'*

By the end of the eighteenth century the theory of parallels, which Jean D'Alembert (1717-1783) had called the 'scandal' of geometry, was in disarray, since the whole of Euclidean geometry rested on a supposed theorem which no one could prove. Every attempt either to deduce the parallel postulate directly from the other four postulates, or alternatively to show that negating the parallel postulate leads to a contradiction, had ended in failure. And to mathematicians of the time one of the most puzzling aspects of the problem was the asymmetry whereby the HOA had been disproved while the HAA continued to resist refutation. As Carl Gauss (1777-1855) put it, the hypothesis of the acute angle is 'the reef on which all the wrecks occur'.

This was the situation which faced young Bolyai when he first embarked on a study of the theory of parallels. No doubt meeting the same impasse as his predecessors in an attempt to refute the HAA, he decided to explore the possibility of creating a geometry based on the HAA, which is now called hyperbolic geometry. This involves retaining Euclid's first four postulates while replacing the parallel postulate with a different postulate. This new postulate states that through a given point, not on a straight line, more than one straight line can be drawn which does not intersect the given line. By 1823 Bolyai had made sufficient progress in this new geometry to write to his father,

*'I am resolved to publish work on parallels as soon as I can put it in order, complete it, and the opportunity arises. I have not yet made the discovery but the path which I have followed is almost certain to lead to my goal, provided this goal is possible. I do not yet have it but I have found things so magnificent that I was astounded. It would be an eternal pity if these things were lost as you, my dear father, are bound to admit when you see them. All I can say now is that I have created a new and different world out of nothing.'*

His father advised him to publish his discoveries as soon as possible since, he said, an idea, when the time is ripe, tends to appear everywhere 'like violets in spring '. He added, prophetically as it turned out,

'*since all scientific striving is only a great war and one does not know when it will be replaced by peace, one must win, if possible; for here, preeminence comes to him who comes first.'*

Publication, however, was apparently delayed by Bolyai's inability to convince his father of the existence of a strange linear constant in his new geometry. Finally in 1832 his treatise, *Appendix Demonstrating the Absolute Science of Space* was published as an appendix to his father's text book, the *Tentamen.* Farkas Bolyai sent a copy to his old university friend, Gauss, the greatest mathematician in Europe. Gauss replied in the following terms:

*'If I commenced by saying that I am unable to praise this work, you would be surprised for a moment. But I cannot say otherwise. To praise it would be to praise myself. Indeed the whole contents of the work the path taken by your son, the results to which he is led, coincide almost entirely with my meditations, which have occupied my mind partly for the last thirty or thirty-five years. So I remain quite stupefied. So far as my own work is concerned, of which up till now I have put very little on paper, my intention was not to let it be published in my lifetime. Indeed the majority of people have not clear ideas upon the questions of which we are speaking, and I have found very few people who could regard with any special interest what I communicated to them on this subject. To be able to take such an interest it is first of all necessary to have devoted careful thought to the real nature of what is wanted and upon this matter almost all are uncertain. On the other hand it was my idea to write down all this later so that it should not perish with me. It is therefore a pleasant surprise for me that I am spared this trouble, and I am very glad that it is just the son of my old friend, who takes the precedence of me in such a remarkable manner.'*

Bolyai was devastated. He even suspected his father of having revealed his ideas to Gauss who was now dishonestly claiming priority. Bolyai even submitted his treatise to the military authorities in Vienna, but his work met with incomprehension. This was hardly surprising given that the *Appendix* appeared to undermine Euclidean geometry, which for over 2000 years had been considered as incontrovertibly true. His cause was not helped by the terse style of the *Appendix*, condensed into just 24 pages, which even Gauss found hard to read. Gauss for his part, though he privately described Bolyai as 'a genius of the first rank', never sought to further the career of his old friend's son, as he easily might have done.

Gauss repeated his claim of priority in a letter to C. F. Gerling in 1832, saying of Bolyai's work that it contained all his own ideas and results, 'expounded with great elegance'. How much credence can be placed in his claim of priority? There is no doubt that Gauss was already familiar with non-Euclidean geometry, for which he had even coined the name. But it seems that his dislike of controversy and the fear of the notoriety that might ensue if he were associated with a challenge to the unique authority of Euclid made him reluctant to publish his ideas. Such discoveries as he chose to reveal were communicated on the basis of the strictest confidence to a trusted circle on whose discretion he could rely.

One such person was Friedrich Wachter (1792-1817), a student of Gauss' at Göttingen. Following a conversation with Gauss on 'anti-geometry', Wachter wrote a letter in 1816 in which he referred to a sphere of infinite radius on whose surface Euclidean geometry would hold true even if the parallel postulate in the plane were proved false. This was a reference to the Euclidean property of the horosphere, which was later to play an important role in the work of Bolyai and Lobachevsky. But it is not clear whether this insight originally came from Gauss or Wachter.

Another correspondent was Franz Taurinus (1794-1874), who had developed a 'logarithmic spherical' geometry on a sphere of imaginary radius. Taking the first fundamental formula of spherical geometry

$$\cos a/k = \cos b/k \cos c/k + \sin b/k \sin c/k \cos A$$

for a spherical triangle ABC on a sphere of radius , Taurinus transformed the real radius $k$ into an imaginary radius $ik$, obtaining

$$\cosh a/k = \cosh b/k \cosh c/k - \sinh b/k \sinh c/k \cos A,$$

which turned out to be the same as the formula for a triangle in the hyperbolic plane. Taurinus also obtained the circumference of a circle of radius $r$ in this geometry as

$$2\pi k \sinh r/k$$

and its area as

$$2\pi k^2 (\cosh r/k - 1),$$

the same as in hyperbolic geometry.

In 1824 Gauss wrote to Taurinus about non-Euclidean geometry in the following terms:

*I have pondered it for over thirty years and I do not believe that anyone can ever have given more thought to it than I, though I have never published anything on it. The assumption that the sum of the three angles (of a triangle) is smaller than 180° leads to a geometry which is quite different from our (Euclidean) geometry but which is in itself completely consistent. I have satisfactorily constructed this geometry for myself so that I can solve every problem, except for the determination of one constant, which cannot be ascertained a priori.*

Then in the letter of 1832 to Farkas Bolyai, in which he acknowledged receipt of the *Appendix,* Gauss included an elegant proof that in hyperbolic geometry the area of a triangle is a function of its defect. In this letter he mentioned his intention of writing down his discoveries in non-Euclidean geometry. In fact he had already begun to write some *Meditations* on the subject. This material, containing a few theorems on parallels, was found among his papers after his death. But the absence of documentary evidence makes it impossible to ascertain with certainty the full extent of his knowledge. It is hard to believe that he had anticipated and proved every single theorem to be found in the *Appendix*, including Bolyai's astonishing quadrature of the circle. Surely Gauss would have mentioned something like that to one of his correspondents. What is certain is that on receiving a copy of Bolyai's treatise he abandoned his own planned memoir.

However what neither Gauss nor Bolyai knew in 1832 was that a Russian mathematician, Nikolai Lobachevsky (1792-1856), had made similar discoveries. He had written a paper in Russian on non-Euclidean geometry, called *A Concise Outline of the Foundations of Geometry*, which was published in 1829 in the *Kazan Messenger* but rejected when submitted to the St Petersburg Academy of Sciences. In 1840 he published another paper on the subject, entitled *Geometrical Investigations on the Theory of Parallels,* this time written in German in the hope of attracting a wider audience. In 1855 he published a longer, technical exposition, *Pangeometry*, written in French.

Gauss obtained a copy of the *Geometrical Investigations* and in 1846 wrote to his friend H.C. Schumacher praising Lobachevsky's skill as a geometer while, typically, maintaining that his treatise contained nothing that was new to him. Nevertheless he proposed Lobachevsky for membership of the Göttingen Academy (which is more than he ever did for Bolyai). In 1848 Bolyai, too, obtained a copy of Lobachevsky's *Geometrical Investigations* and expressed astonishment that another mathematician had arrived at the same results, though by a different route.

Sadly, neither Bolyai nor Lobachevsky received due recognition for their achievements in their lifetime. Bolyai took early retirement from the army in 1833 on the grounds of ill health and returned to Transylvania where he lived on a small pension. He continued to engage in mathematics but published nothing more. Lobachevsky, having started as a student at the University of Kazan, rose to become Lecturer, Professor and finally Rector of the University. But he was dismissed from his post in 1846, ostensibly due to his deteriorating health. Unable to walk and half-blind in his later years, he died in poverty.

Interest in hyperbolic geometry only revived after Gauss' death when his correspondence on the subject, which contained references to the work of Lobachevsky and Bolyai, was published. Mathematicians became interested in tracing the original papers, which within a short time were translated into French and Italian. The American mathematician, George Halsted, provided the first English translation of Lobachevsky's *Geometrical Investigations* in 1892 and Bolyai's *Appendix* in 1896.

Following the publication in 1868 of Bernhard Riemann's seminal paper *On the Principles Underlying Geometry*, non-Euclidean geometry was extended by the discovery of elliptic geometry. Riemann (1822-1866) had made a distinction between lines which are unbounded and lines which are infinite. It was found that if Euclid's second postulate, that a straight line can be infinitely extended, is modified to state that the straight line is finite but unbounded, while the parallel postulate is replaced by the assumption that no two straight lines are parallel, then it is possible to create a geometry in which the sum of the angles of a triangle is greater than two right angles, analogous to the geometry which holds true on the surface of a sphere. The term elliptic was coined by Felix Klein to denote this type of non-Euclidean geometry. But the most important aspect of Riemann's paper was his introduction of the concept of what is now called a manifold. This provided a wider geometrical framework into which non-Euclidean geometry, previously regarded by some mathematicians as a mere curiosity, could takes its place alongside Euclidean geometry.

The discovery of non-Euclidean geometry also had implications for physics. For over 2000 years Euclid's geometry had been regarded as synonymous with the geometry of space and therefore the only conceivable geometry. The discovery of non-Euclidean geometry left geometry in the uncomfortable position of being a deductive science based on assumptions which could not be validated either empirically or mathematically.

Paradoxically the very quandary in which geometry found itself released it from subordination to the physical world. Mathematicians felt free thereafter to create new geometries, however abstract they at first appeared to be. But once these geometries had been shown to be logically consistent, they could be applied to the physical world and tested to see if they fitted the observed phenomena. It was found, for example, that hyperbolic geometry can model Einstein's special theory of relativity (see Varičak, 1912).

APPENDIX

# THE ABSOLUTELY TRUE SCIENCE OF SPACE

EXPOUNDED INDEPENDENTLY OF THE CORRECTNESS OR FALSENESS (A PRIORI UNDECIDABLE FOR EVER) OF EUCLID'S AXIOM XI: FOR THE CASE OF FALSENESS, WITH A GEOMETRIC QUADRATURE OF THE CIRCLE

BY

JÁNOS BOLYAI

CAPTAIN-ENGINEER OF THE IMPERIAL AND ROYAL AUSTRIAN ARMY

## EXPLANATION OF SIGNS

Les us denote by

| | |
|---|---|
| $AB$ | the collection of all points on the straight line passing through the points $A$, $B$; |
| $\overrightarrow{AB}$ | the half-line (ray) from $A$ through $B$; |
| $ABC$ | the collection of all points in the plane containing the non-collinear points $A, B, C$; |
| $\|AB\|C$ | the half-plane having boundary $AB$ and containing $C$; |
| $\sphericalangle ABC$ | the smaller of the angular domains bounded by the arms $\overrightarrow{BA}$, $\overrightarrow{BC}$; |
| $ABCD$ | (if $D$ is an interior point of $\sphericalangle ABC$ and $A, B, CD$ do not intersect each other) the plane domain bounded by $\overrightarrow{BA}$, $BC$, $\overrightarrow{CD}$ and contained in $\sphericalangle ABC$; |
| $(AB, CD)$ | (if $AB$, $CD$ are coplanar and do not intersect) the plane domain bounded by $AB$ and $CD$; |
| $R$ | the right angle; |
| $\equiv$ | congruence*; |
| $AB \simeq CD$ | the relation $\sphericalangle CAB = \sphericalangle ACD$; |
| $x \to a$ | that $x$ tends to the limit $a$; |
| $\bigcirc r$ | the circumference of the circle of radius $r$; |
| $\odot r$ | the area of the circle of radius $r$. |

**§1**

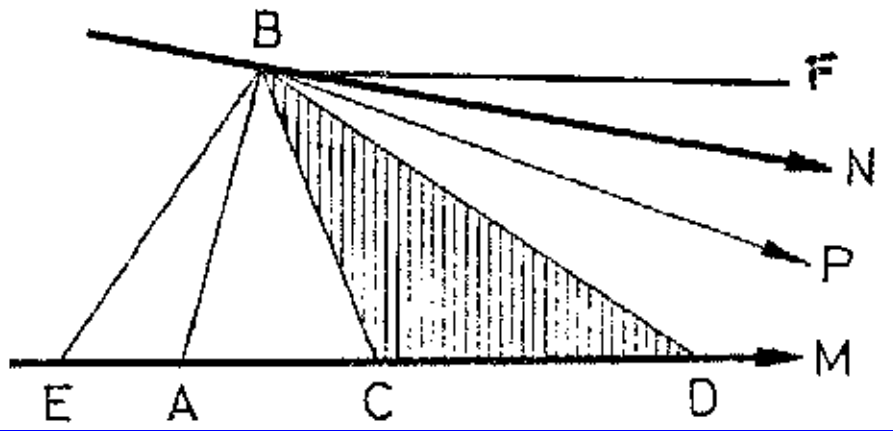

If $\overrightarrow{BN}$ does not, while any other $\overrightarrow{BP}$ in the angular domain $ABN$ does intersect $\overrightarrow{AM}$, we write

$$\overrightarrow{BN} \parallel \overrightarrow{AM}.$$

Obviously, each point $B$ outside the line $AM$ is the origin of one and only one $\overrightarrow{BN}$ with this property. Moreover,

$$\sphericalangle BAM + \sphericalangle ABN \leqq 2R.$$

For, if $BC$ rotates about $B$ until

$$\sphericalangle BAM + \sphericalangle ABC = 2R,$$

there will be a first position where $\overrightarrow{BC}$ does not intersect $\overrightarrow{AM}$, and in this position we have $\overrightarrow{BN} \parallel \overrightarrow{AM}$.

It is also clear that $\overrightarrow{BN} \parallel \overrightarrow{EM}$ for each point $E$ on the line $AM$ (assuming in all such cases that $M$ has been chosen so as to satisfy $\overline{AM} > \overline{AE}$).

If $C$ on $AM$ goes to infinity, and always $\overline{CD} = \overline{BC}$, then

$$\sphericalangle CDB = \sphericalangle CBD < \sphericalangle NBC.$$

But

$$\sphericalangle NBC \to 0.$$

Therefore

$$\sphericalangle ADB \to 0.$$

***Note on 1***

Bolyai defines **BN as *parallel* to AM if it does not intersect AM while every other straight line passing through B within ∠ABN (such as BP) does intersect AM.** This definition of parallelism is 'absolute' using Bolyai's term, or neutral in modern terminology, in that it makes no assumption about the validity or otherwise of the parallel postulate.

If the parallel postulate is accepted, then BN is the only line parallel to AM passing through B and the sum of angles ABN, BAM is two right angles.

If it is replaced by the hyperbolic postulate, there will be other lines passing through B outside of ∠ABN (such as BF) which also fail to intersect AM, generally called ultra-parallel lines nowadays. The sum of angles ABN, BAN will be less than two right angles. Furthermore the fact that ∠NBC $\to$ 0 implies that the parallel BN will meet AM at an infinitely distant point. In other words BN is asymptotic to AM.

§2

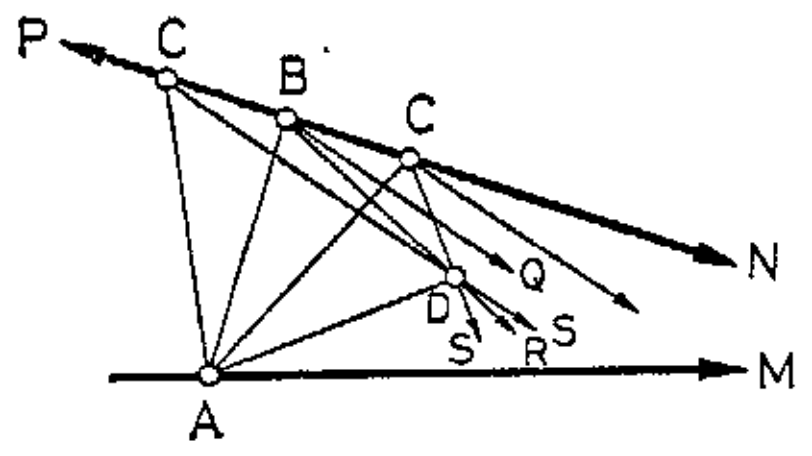

Figure 2

*If $\overrightarrow{BN} \| \overrightarrow{AM}$, then also $\overrightarrow{CN} \| \overrightarrow{AM}$.*

For let $D$ be an interior point of the plane domain $MACN$. If $C$ is on $\overrightarrow{BN}$, then in view of $\overrightarrow{BN} \| \overrightarrow{AM}$ the half-line $\overrightarrow{BD}$ intersects $\overrightarrow{AM}$ and, therefore, also $\overrightarrow{CD}$ intersects $\overrightarrow{AM}$. If, on the other hand, $C$ is on $\overrightarrow{BP}$, let $\overrightarrow{BQ} \| \overrightarrow{CD}$. By §1 the half-line $\overrightarrow{BQ}$ lies in the interior of $\sphericalangle ABN$ and intersects $\overrightarrow{AM}$. So $\overrightarrow{CD}$ intersects $\overrightarrow{AM}$.

Thus in both cases any $\overrightarrow{CD}$ in the interior of $\sphericalangle ACN$ does, while $\overrightarrow{CN}$ itself does not intersect $\overrightarrow{AM}$. Consequently, $\overrightarrow{CN} \| \overrightarrow{AM}$.

***Note on 2***

This theorem demonstrates that parallelism is independent of the choice of B on the half line BN parallel to AM. In other words, **a straight line parallel to another maintains the property of parallelism along all its points.**

§3

*If $\overrightarrow{BR} \| \overrightarrow{AM}$ and $\overrightarrow{CS} \| \overrightarrow{AM}$, but $C$ is not on the straight line $BR$, then $\overrightarrow{BR}$ and $\overrightarrow{CS}$ do not intersect each other.*

For, assuming that $D$ is a common point of $\overrightarrow{BR}$ and $\overrightarrow{CS}$, §2 yields $\overrightarrow{DR} \| \overrightarrow{AM}$ and $\overrightarrow{DS} \| \overrightarrow{AM}$. Thus by §1 the half-lines $\overrightarrow{DR}$ and $\overrightarrow{DS}$ coincide, so that $C$ would be a point of the line $BR$ contrary to the hypothesis.

***Note on 3***

In Euclidean geometry two straight lines parallel to a third are parallel to each other. This theorem only proves that that **two half lines parallel to a third do not intersect one another**, since it cannot be assumed that non-intersecting lines are the same as parallel lines according to the definition of parallelism given in Section 1.

The proof is by contradiction. Suppose that BR and CS intersect at a point D (see Figure 2). Then, given that by hypothesis BR and CS are parallel to AM, the result of Section 2 implies that DR and DS must also be parallel to AM. Now by Section 1 only one line through D can be parallel to AM. Therefore DR and DS coincide. This implies that CS and BR are collinear, in other words that C is a point on BR, contrary to the hypothesis.

## §4

*If $\sphericalangle MAN > \sphericalangle MAB$, then to any point B of $\overrightarrow{AB}$ there can be found a point C on $\overrightarrow{AM}$ such that*

$$\sphericalangle BCM = \sphericalangle NAM.$$

For §1 assures the existence of $D$ such that

$$\sphericalangle BDM > \sphericalangle NAM$$

and, consequently, of $P$ such that

$$\sphericalangle MDP = \sphericalangle MAN;$$

in this case $B$ lies in the plane domain $NADP$. If we shift $NAM$ along $\overrightarrow{AM}$ until $\overrightarrow{AN}$ takes the position $\overrightarrow{DP}$, then $\overrightarrow{AN}$ will at some time pass through $B$. Therefore necessarily

$$\sphericalangle BCM = \sphericalangle NAM$$

for some point $C$.

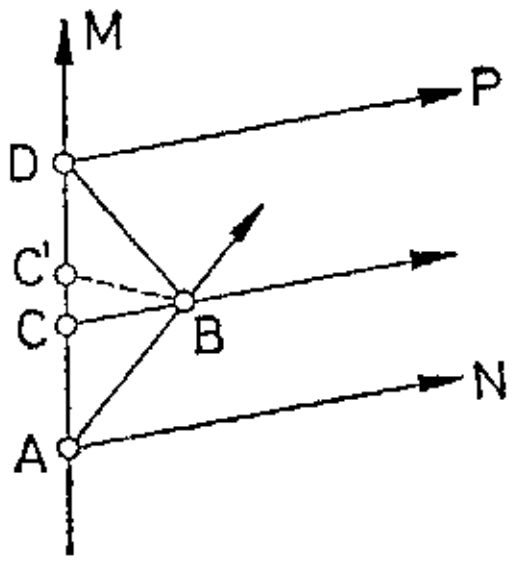

**Figure 3**

***Note on 4***

This theorem shows shows that **if ∠MAN is shifted along AM, then its other arm AN will sweep out the angular domain MAN in a continuous manner**. Bolyai implies (though he does not prove) that for any point B there is only one such C. For if there were two points, C and C′, such that

$$\angle BCM = \angle BC'M = \angle NAM,$$

then in view of the fact that

$$\angle BC'C + \angle BC'M = 2R,$$

it would follow that

$$\angle BCC' + \angle BC'C = 2R.$$

But this is impossible by Euclid's (neutral) theorem (*Elements*, 1, 17) which states that the sum of any two angles in a triangle is less than two right angles.

## §5

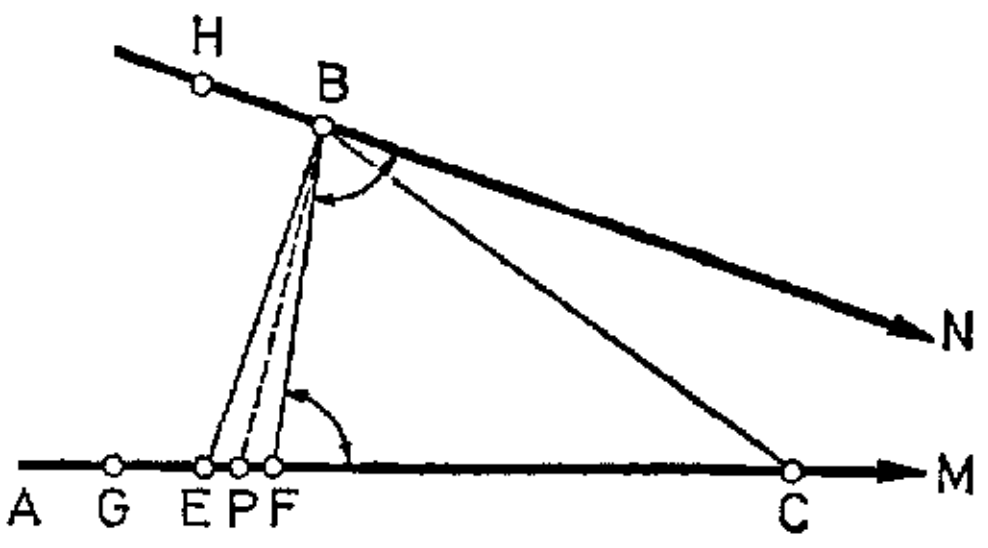

**Figure 4**

*If $\overrightarrow{BN} \parallel \overrightarrow{AM}$, then there is a point F on the line AM such that*

$$FM \simeq BN.$$

Really, by §1 there exists $\sphericalangle BCM$ which is greater than $\sphericalangle CBN$. Further, if $CE = CB$ then

$$EC \simeq BC.$$

Thus

$$\sphericalangle BEM \prec \sphericalangle EBN.$$

Let the point $P$ run through $EC$ and let $u = \sphericalangle BPM$, $v = \sphericalangle PBN$. It is clear that $u$ is initially smaller but finally greater than the corresponding $v$. However, $u$ increases from $\sphericalangle BEM$ to $\sphericalangle BCM$ in a *continuous* manner; in fact, by §4 there is no angle greater than $\sphericalangle BEM$ and smaller than $\sphericalangle BCM$ to which $u$ would not become equal at some time. Similarly, $v$ continuously decreases from $\sphericalangle EBN$ to $\sphericalangle CBN$. So there is a point $F$ on $EC$ having the property

$$\sphericalangle BFM = \sphericalangle FBN.$$

***Note on 5***

Here Bolyai demonstrates the existence of ***corresponding points*** on any pair of parallel lines such as BN and AM, two points B, F being defined as corresponding if the connecting transversal BF makes equal angles with the parallels i.e. if ∠BFM = ∠FBN. He begins by showing that there is a point E on AM such that ∠BEM < ∠EBN. He then allows a point P to run through EC. When P coincides with E, then ∠BPM < ∠PBN, but when it coincides with C, then ∠BPM > ∠PBN. Since ∠BPM increases continuously while ∠PBN decreases continuously, Section 4 implies that there must exist a point F on EC such that ∠BFM = ∠FBN. Thus F is the required point corresponding to B.

## §6

*If* $\overrightarrow{BN} \| \overrightarrow{AM}$, *whereas G and H are arbitrary points on the lines AM and BN, respectively, then* $\overrightarrow{HN} \| \overrightarrow{GM}$ *and* $\overrightarrow{GM} \| \overrightarrow{HN}$.

For **§1** yields $\overrightarrow{BN} \| \overrightarrow{GM}$ and by **§2** this implies $\overrightarrow{HN} \| \overrightarrow{GM}$.

Further, if we choose $F$ according to **§5** so that $FM \simeq BN$ then $MFBN = NBFM$. But $\overrightarrow{BN} \| \overrightarrow{FM}$; therefore $\overrightarrow{FM} \| \overrightarrow{BN}$ and by the preceding paragraph $\overrightarrow{GM} \| \overrightarrow{HN}$.

***Note on 6***

Section 6 demonstrates that **parallelism is a *symmetric* relation** i.e. if BN || AM, then AM || BN (see Figure 4). The proof depends on the existence of a point F on AM corresponding to B. By symmetry the existence of corresponding points implies that each line is parallel to the other.

## §7

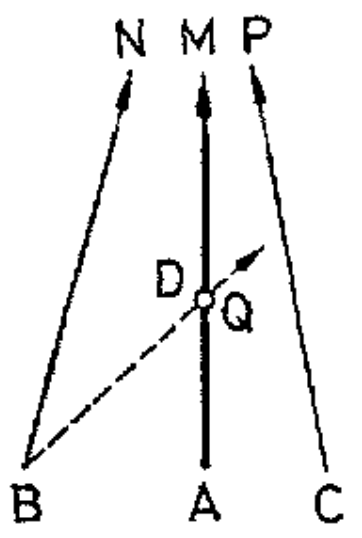

Figure 5

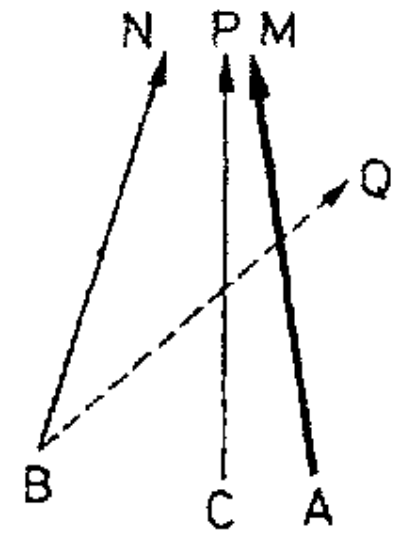

Figure 6

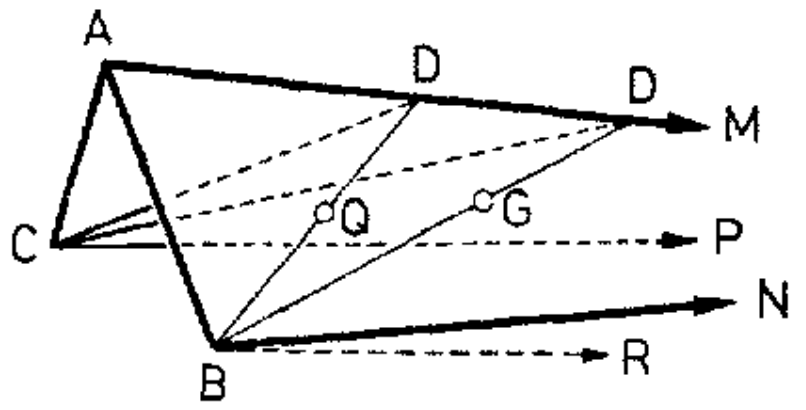

Figure 7

*If* $\overrightarrow{BN}\parallel\overrightarrow{AM}$ *as well as* $\overrightarrow{CP}\parallel\overrightarrow{AM}$ *and* $C$ *is not on the line* $BN$, *then* $\overrightarrow{BN}\parallel\overrightarrow{CP}$.

In fact, by §3 the half-lines $\overrightarrow{BN}$ and $\overrightarrow{CP}$ do not intersect each other. $\overrightarrow{AM}$, $\overrightarrow{BN}$ and $\overrightarrow{CP}$ may be either coplanar or not. In the first case $\overrightarrow{AM}$ may be either in the strip $(BN, CP)$ or outside it.

If $\overrightarrow{AM}$, $\overrightarrow{BN}$, $\overrightarrow{CP}$ are coplanar and $AM$ lies in the strip $(BN\ CP)$, then any $\overrightarrow{BQ}$ in $\sphericalangle NBC$ intersects the line $AM$ in some points $D$, since $\overrightarrow{BN}\parallel\overrightarrow{AM}$. According to §6

$$\overrightarrow{DM}\parallel\overrightarrow{CP},$$

so that $\overrightarrow{DQ}$ intersects $\overrightarrow{CP}$, and

$$\overrightarrow{BN}\parallel\overrightarrow{CP}.$$

If however, $BN$ and $CP$ are on the same side of $AM$, then one of them, say $CP$, lies between the other two, $BN$ and $AM$. Therefore any $\overrightarrow{BQ}$ in $\sphericalangle NBA$ intersects $\overrightarrow{AM}$ and, consequently, it intersects $\overrightarrow{CP}$. Thus

$$\overrightarrow{BN}\parallel\overrightarrow{CP}.$$

If the planes $MAB$, $MAC$ form an angle, then $\sphericalangle CBN$ has no common point with $\sphericalangle ABN$ outside $\overrightarrow{BN}$, and $\overrightarrow{AM}$ in $\sphericalangle ABN$ has no common point with $\overrightarrow{BN}$. Hence $\sphericalangle NBC$ has no common point with $\overrightarrow{AM}$. However, any half-plane $|BC|Q$ that contains $\overrightarrow{BQ}$, the latter lying in $\sphericalangle NBA$, intersects $\overrightarrow{AM}$. For $\overrightarrow{BN}\parallel\overrightarrow{AM}$ implies that $\overrightarrow{BQ}$ intersects $\overrightarrow{AM}$ in a point $D$. Let us turn the half-plane $|BC|D$ about $BC$ until it leaves $\overrightarrow{AM}$ for the first time. Then $|BC|D$ must coincide with $|BC|N$ and, for the same reason, with $|BC|P$. Therefore $\overrightarrow{BN}$ belongs to the half-plane $|BC|P$. Further, if $\overrightarrow{BR}\parallel\overrightarrow{CP}$ then using the relation $\overrightarrow{AM}\parallel\overrightarrow{CP}$ and following the line of argument given above it can be seen that $\overrightarrow{BR}$ belongs to the half-plane $|BA|M$. As $\overrightarrow{BR}\parallel\overrightarrow{CP}$, the half-line $\overrightarrow{BR}$ must belong also to $|BC|P$. Consequently, $\overrightarrow{BR}$ is the intersection of $\sphericalangle MAB$ and $\sphericalangle PCB$. This intersection is therefore identical with $\overrightarrow{BN}$. Hence

$$\overrightarrow{BN}\parallel\overrightarrow{CP}.$$

Thus if $\overrightarrow{CP}\parallel\overrightarrow{AM}$ and the point $B$ is outside the plane $CAM$, then the intersection of $\sphericalangle BAM$ and $\sphericalangle BCP$, that is the half-line $\overrightarrow{BN}$, is parallel to each of $\overrightarrow{AM}$ and $\overrightarrow{CP}$*.

***Note on 7***

Section 7 provides a proof that **parallelism is a *transitive* relation** i.e. that if two straight lines are parallel to the same line, they are parallel to each other. He examines two cases. The first case, in which the parallels are coplanar, is divided into two subcases: the two lines either lie on the same side of the third line or lie on either side. In the second case the lines form three different planes.

**§8**

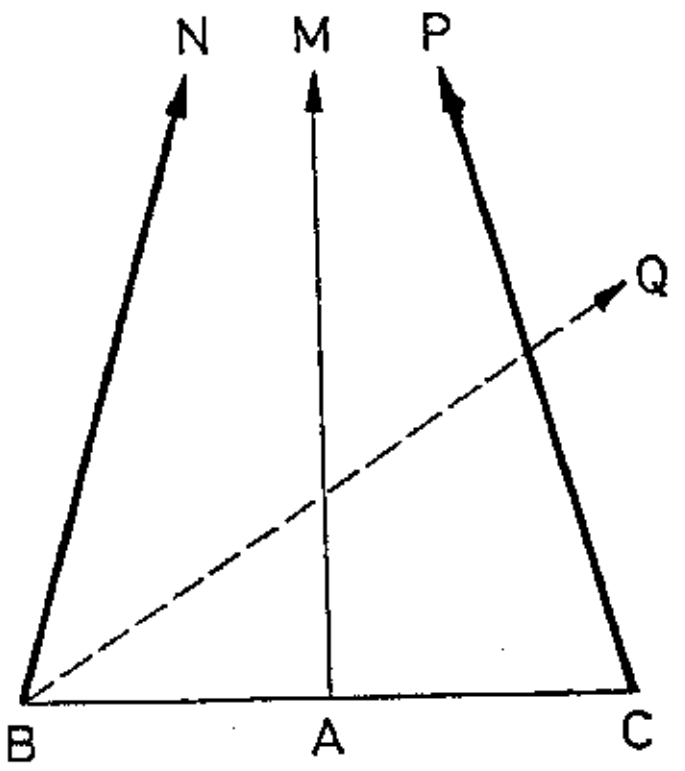

**Figure 8**

*If $\overrightarrow{BN}\|$ and $\simeq\overrightarrow{CP}$, or $\overrightarrow{BN}\|\simeq\overrightarrow{CP}$ for short, and $\overrightarrow{AM}$ in the domain NBCP perpendicularly bisects the distance $\overline{BC}$, then $\overrightarrow{BN}\|\overrightarrow{AM}$.*

For if $\overrightarrow{BN}$ intersects $\overrightarrow{AM}$ then, in view of $MABN=MACP$, the half-line $\overrightarrow{CP}$ intersects $\overrightarrow{AM}$ at the same point. The latter would be a common point of $\overrightarrow{BN}$ and $\overrightarrow{CP}$ although $\overrightarrow{BN}\|\overrightarrow{CP}$.

Any $\overrightarrow{BQ}$ in the angular domain $CBN$ intersects $\overrightarrow{CP}$. Therefore $\overrightarrow{BQ}$ intersects $\overrightarrow{AM}$ as well. Consequently,

$$\overrightarrow{BN}\parallel\overrightarrow{AM}.$$

***Note on 8***

In this important theorem, subsequently used in his quadrature of the circle, Bolyai proves that **the perpendicular bisector of the line which joins a pair of parallel lines at corresponding points is parallel to both lines.**

In Figure 8 BN | | CP, while B and C are corresponding points. AM is the perpendicular bisector of BC. Bolyai first proves by contradiction that BN and AM cannot intersect. For if they intersected, then by symmetry (since BA = AC) CP and AM would intersect at the same point, in which case BN and CP would also intersect, contrary to the hypothesis. Then he shows that any line BQ within ∠ABN, since it intersects CP, must necessarily intersect AM, implying that BN | | AM, as required.

## §9

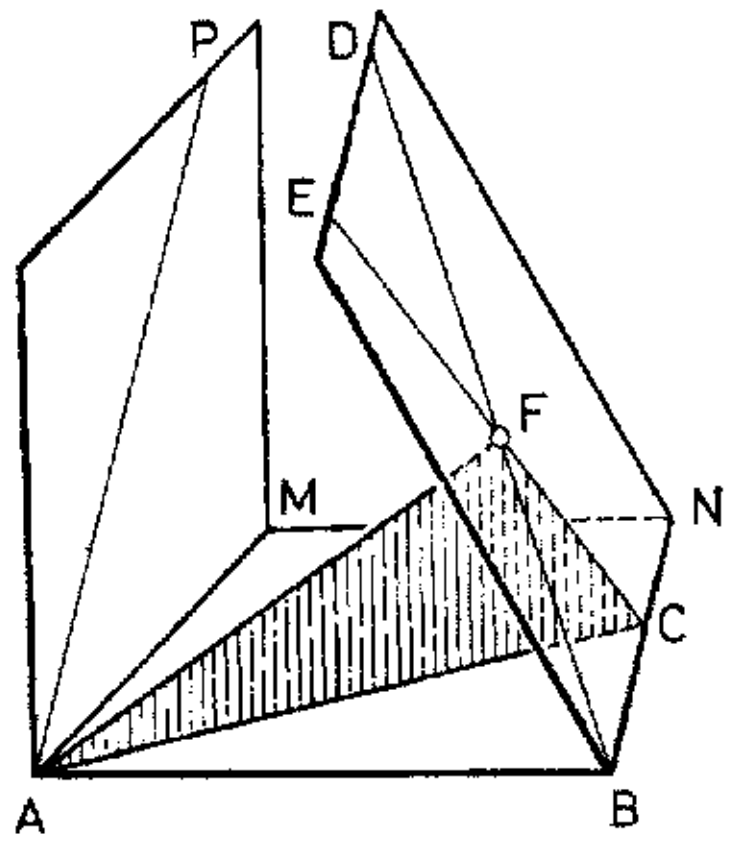

**Figure 9**

*If* $\overrightarrow{BN}\|\overrightarrow{AM}$ *and* $MAP \perp MAB$, *while the dihedral angle formed by* $|NB|D$ *and* $|NB|A$ *on that side of the plane* $MABN$ *which contains* $|MA|P$ *is smaller than* $R$, *then the half-planes* $|MA|P$ *and* $|NB|D$ *intersect each other.*

For let

$$BA \perp AM, \quad AC \perp BN$$

(whether $B$ and $C$ coincide or not) and

$$CE \perp NB$$

in $|NB|D$. Then the hypothesis implies that

$$\sphericalangle ACE < R,$$

and the segment $\overline{AF}$, perpendicular to $CE$, lies in $\sphericalangle ACE$.

Each of the half-planes $|AB|F$ and $|AM|P$ contains the point $A$; let their line of intersection be $AP$. Then

$$\sphericalangle BAP = \sphericalangle BAM = R,$$

since the planes $BAM$ and $MAP$ are perpendicular.

Finally, if $|AB|F$ is turned about $AB$ into $|AB|M$, then $\overrightarrow{AP}$ turns into $\overrightarrow{AM}$. As

$$AC \perp BN \quad \text{and} \quad \overline{AF} < \overline{AC},$$

it is clear that, after the rotation, $\overrightarrow{AF}$ will not reach $\overrightarrow{BN}$, and therefore $\overrightarrow{BF}$ will fall in the interior of $\sphericalangle ABN$. In this position, $\overrightarrow{BF}$ intersects $\overrightarrow{AP}$ since $\overrightarrow{BN}\|\overrightarrow{AM}$. Thus $\overrightarrow{BF}$ intersects $\overrightarrow{AP}$ in the original position too, and the point of intersection is a common point of $|MA|P$ and $|NB|D$. Consequently, the half-planes $|MA|P$ and $|NB|D$ intersect each other.

Hence it easily follows that *the half-planes* $|MA|P$ *and* $|NB|D$ *intersect whenever the sum of the angles they form with the plane domain* $MABN$ *is* $< 2R$.

***Note on 9***

In this striking theorem Bolyai proves a neutral three dimensional version of the parallel postulate which remains valid regardless of whether the parallel postulate in the plane is assumed to be true or not. The theorem shows that **if the sum of the dihedral angles which two planes make with a third plane is less than two right angles, the lines of intersection between the planes being parallel, then the first two planes will intersect.**

In Figure 9 BN is drawn parallel to AM. This implies that any line drawn inside $\angle$MAB such as AC will cut BN and similarly any line drawn inside $\angle$ABN will cut AM. $\angle$MAB is a right angle but no assumption is made about the magnitude of $\angle$ABN which will either be a right angle or an acute angle according to whether the parallel postulate is accepted or rejected. AC is drawn perpendicular to BN. (If the parallel postulate is assumed to be true, then $\angle$ABN will be a right angle and the points B and C will coincide. If $\angle$ABN is acute, then C will lie somewhere along BN, as shown in the figure).

At this point Bolyai introduces a third dimension. He draws planes PAM and BND, where PAM is perpendicular to plane ABNM but where the dihedral angle between plane BND and plane ABNM is acute. In other words the sum of the dihedral angles which the planes PAM and BND make with the base plane ABNM is less than two right angles. The aim is to show that the planes PAM and BND meet.

His proof requires three further constructions. Draw CE perpendicular to BN in the plane BND. Next draw AF perpendicular to CE meeting CE at F. (Since the dihedral angle between planes BND and ABNM is assumed to be acute, the point F will lie somewhere along CE, as shown). This creates triangle ACF (shaded) Finally join B and F. This creates a fourth plane, PABF.

Bolyai's argument then runs as follows. Suppose that the plane PABF is rotated about AB until it lies along ABNM. The fact that plane PAM is perpendicular to plane ABNM while $\angle$MAB is a right angle ensures that AP, after rotation, will lie along AM. What about AF? AC, being perpendicular to BN, represents the shortest distance between A and BN. But AF < AC, since AC is the hypotenuse of triangle ACF, so after rotation AF will not reach BN. This in turn implies that BF will lie inside $\angle$ABN. In this position BF must intersect AP since AM and BN are parallel by hypothesis. Therefore BF must intersect AP in the original position too, before rotation. But AP lies in plane MAP while BF lies in plane BND. Hence the two planes MAP and BND intersect, as required.

## §10

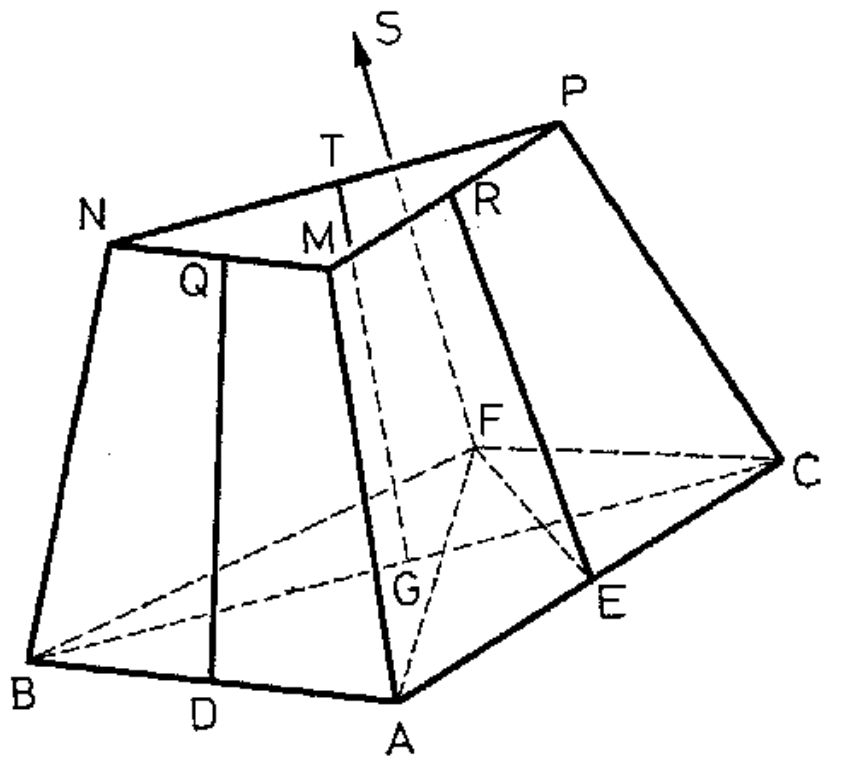

**Figure 10**

*If* $\overrightarrow{BN}$ *as well as* $\overrightarrow{CP}$ *are* $\parallel\simeq\overrightarrow{AM}$, *then* $\overrightarrow{BN}\parallel\simeq\overrightarrow{CP}$.

In fact, the planes $MAB$ and $MAC$ either form a dihedral angle or coincide.

*The first case.* Let the plane $QDF$ perpendicularly bisect $\overline{AB}$. Then $DQ\perp AB$, so that $\overrightarrow{DQ}\parallel\overrightarrow{AM}$ by §8. Similarly, if $ERS$ perpendicularly bisects $\overline{AC}$, then $\overrightarrow{ER}\parallel\overrightarrow{AM}$. Therefore $\overrightarrow{DQ}\parallel\overrightarrow{ER}$ by §7. Hence, in view of §9, it follows easily that the planes $QDF$ and $ERS$ intersect each other, and by §7 their line of intersection $FS$ satisfies

$$\overrightarrow{FS}\parallel\overrightarrow{DQ}.$$

As

$$\overrightarrow{BN}\parallel\overrightarrow{DQ},$$

also

$$\overrightarrow{FS}\parallel\overrightarrow{BN}.$$

Further, for any point $F$ on the line $FS$ we have

$$\overline{FB}=\overline{FA}=\overline{FC};$$

thus $FS$ lies in the plane $TGF$ which perpendicularly bisects $\overline{BC}$. But §7 and $\overrightarrow{FS}\parallel\overrightarrow{BN}$ imply that

$$\overrightarrow{GT}\parallel\overrightarrow{BN}.$$

A similar argument yields

$$\overrightarrow{GT}\parallel\overrightarrow{CP}.$$

Moreover, $\overrightarrow{GT}$ perpendicularly bisects $\overline{BC}$. So by §1

$$TGBN\equiv TGCP$$

and, consequently,

$$\overrightarrow{BN}\parallel\simeq\overrightarrow{CP}.$$

*The second case.* If $\overrightarrow{BN}$, $\overrightarrow{AM}$ and $\overrightarrow{CP}$ are in one plane, let the half-line $\overrightarrow{FS}$ lying outside this plane satisfy

$$\overrightarrow{FS}\parallel\simeq\overrightarrow{AM}.$$

According to the first case,

$$\overrightarrow{FS}\parallel\simeq\overrightarrow{BN}\quad\text{and}\quad\overrightarrow{FS}\parallel\simeq\overrightarrow{CP}.$$

Therefore

$$\overrightarrow{BN}\parallel\simeq\overrightarrow{CP}.$$

***Notes on 10***

This theorem shows that **the correspondence between points on parallel lines is a *transitive* relation** i.e. if BN || AM where N, M are corresponding points and if CP || AM where M, P are corresponding points, then BN || CP and B, C are also corresponding points.

§11

Denote by **F** the collection which consists of the point $A$ and all points $B$ such that $\overrightarrow{BN} \| \overrightarrow{AM}$ implies $BN \simeq AM$. The intersection of **F** with any plane containing the line $AM$ will be denoted by **L**.

On every line parallel to $\overrightarrow{AM}$ there is one and only one point of **F**. Clearly, $\overrightarrow{AM}$ divides **L** into two congruent parts. We say, $\overrightarrow{AM}$ is an axis of **L**. It is also obvious that in any plane which contains $\overrightarrow{AM}$ there is one and only one **L** having $\overrightarrow{AM}$ for axis. In the plane considered, this **L** will be said to correspond to $\overrightarrow{AM}$. Clearly, if **L** rotates about $\overrightarrow{AM}$, it describes an **F**. We say that $\overrightarrow{AM}$ is an axis of the **F** so obtained and that, in turn, this **F** corresponds to $\overrightarrow{AM}$.

***Note on 11***

Here Bolyai introduces a key figure, which he calls an **F surface,** though horosphere is the term generally used nowadays. (Lobachevsky called it a boundary surface, while Gauss used the term parasphere). The horosphere is a figure unique to hyperbolic geometry. Bolyai defines it as the set of mutually corresponding points in a bundle of parallel lines. It can be viewed as a sphere of infinite radius.

A plane containing any one of the parallels in the bundle intersects the horosphere in what Bolyai calls an **L line**, though the term horocycle is generally used nowadays. (Lobachevsky called it a limiting curve and Gauss a paracycle). A horocycle can be viewed as a circle of infinite radius.

If a horocycle is rotated about any of the parallel axes which it intersects, it generates a horosphere:

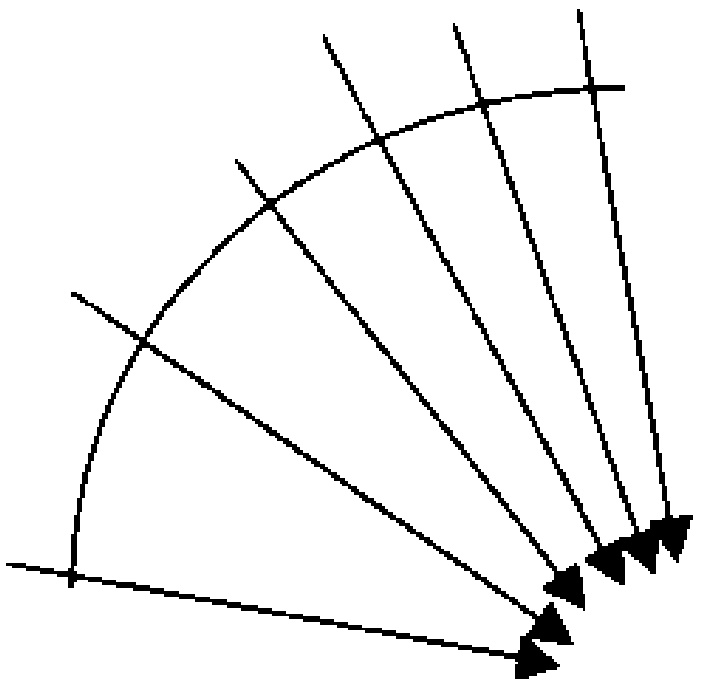

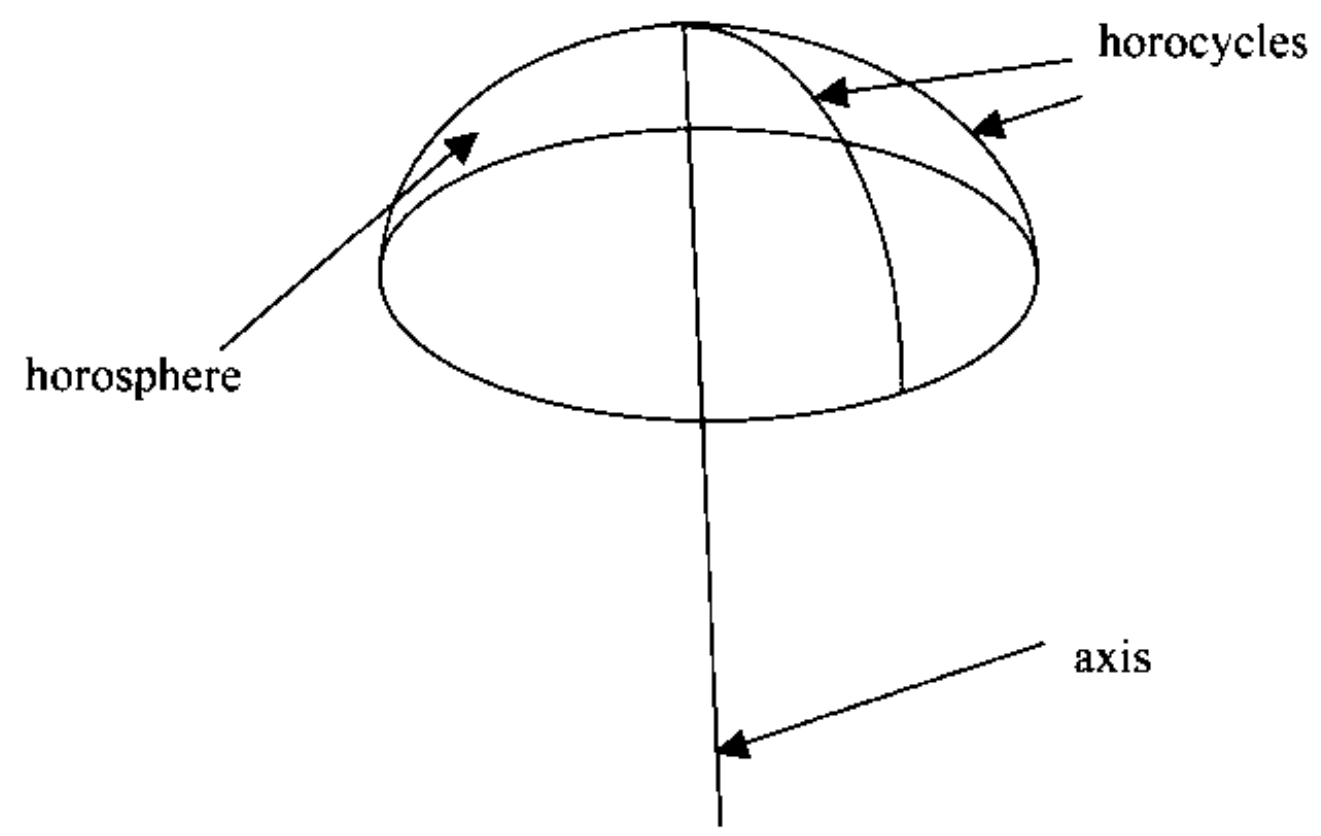

## §12

*If B is a point of* **L**, *where* **L** *corresponds to* $\overrightarrow{AM}$, *and*

$$\overrightarrow{BN} \parallel \simeq \overrightarrow{AM}$$

*in accordance with* **§11**, *then* **L** *corresponding to* $\overrightarrow{AM}$ *coincides with* **L** *corresponding to* $\overrightarrow{BN}$.

For better distinction, **L** corresponding to $\overrightarrow{BN}$ will be denoted by **l**. Let $C$ be a point of **l** and, in accordance with **§11**,

$$\overrightarrow{CP} \parallel \simeq \overrightarrow{BN}.$$

As $\overrightarrow{BN} \parallel \simeq \overrightarrow{AM}$, **§10** yields

$$\overrightarrow{CP} \parallel \simeq \overrightarrow{AM},$$

so that $C$ lies also on **L**. If, on the other hand, $C$ is a point of **L** and $\overrightarrow{CP} \parallel \simeq \overrightarrow{AM}$, then by **§10**

$$\overrightarrow{CP} \parallel \simeq \overrightarrow{BN}$$

and therefore by **§11** $C$ lies also on **l**. Consequently, **L** and **l** are identical and $\overrightarrow{BN}$ is an axis also for **L**. Thus, all axes of **L** are related by $\parallel \simeq$.

The same property of **F** can be proved in a similar way.

### *Note on 12*

This theorem confirms that **only one horocycle passes through the corresponding points of any pair of parallel lines.** Similarly only one horosphere passes through the mutually corresponding points on a bundle of parallel lines.

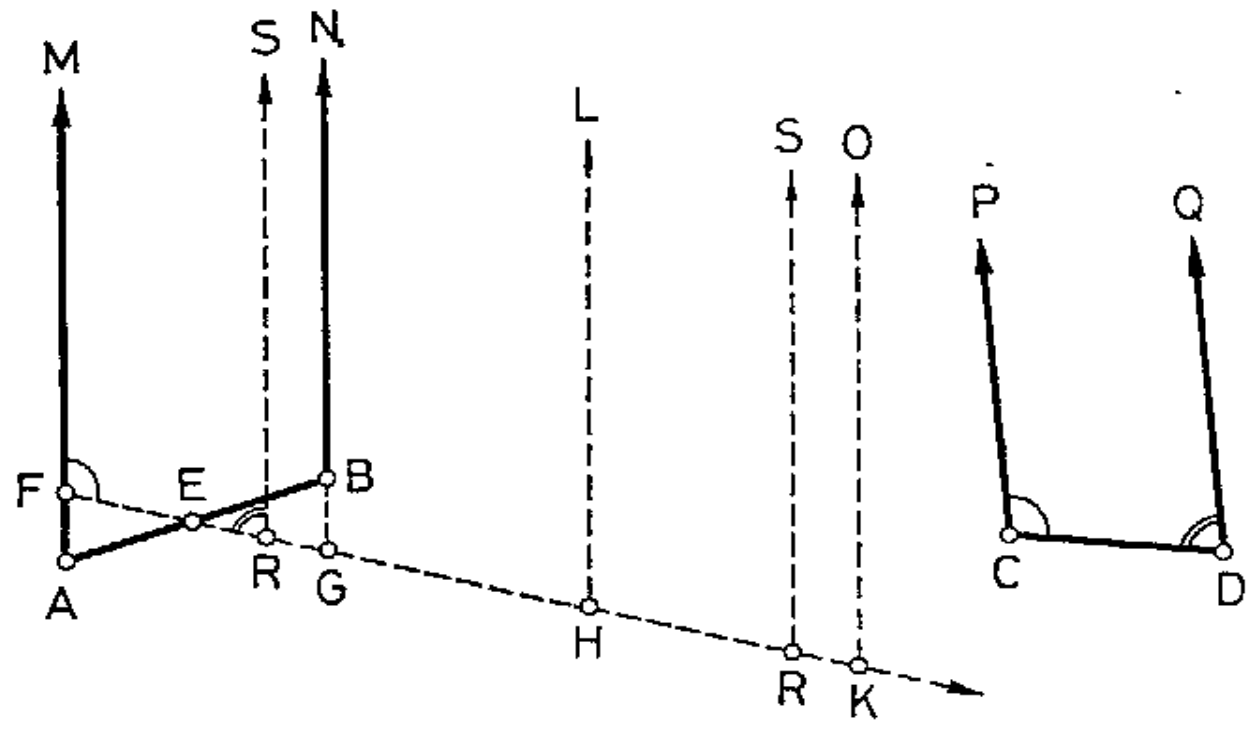

**Figure 11**

§13

*If* $\overrightarrow{BN}\parallel\overrightarrow{AM}$, $\overrightarrow{CP}\parallel\overrightarrow{DQ}$ *and* $\sphericalangle BAM+\sphericalangle ABN=2R$, *then* $\sphericalangle DCP+\sphericalangle CDQ=2R$.

For let $\overline{EA}=\overline{EB}$ and $\sphericalangle EFM=\sphericalangle DCP$, which can be achieved by §4. As

$$\sphericalangle BAM+\sphericalangle ABN=2R=$$
$$=\sphericalangle ABN+\sphericalangle ABG,$$

we have

$$\sphericalangle EBG=\sphericalangle EAF.$$

Therefore if also $\overline{BG}=\overline{AF}$ then $\Delta EBG\equiv\Delta EAF$ and $\sphericalangle BEG=\sphericalangle AEF$, so that $G$ is on the half-line $\overrightarrow{FE}$. Further

$$\sphericalangle GFM+\sphericalangle FGN=2R,$$

since $\sphericalangle EGB=\sphericalangle EFA$. Also, in view of §6, $\overrightarrow{GN}\parallel\overrightarrow{FM}$. So if

$$MFRS\equiv PCDQ$$

then by §7 $\overrightarrow{RS}\parallel\overrightarrow{GN}$. Unless $\overline{CD}=FG$, in which case the assertion is obviously true, $R$ lies on either the segment $\overline{FG}$ or its extension.

(I) In the first case, $\sphericalangle FRS$ is not greater than

$$2R-\sphericalangle RFM=\sphericalangle FGN,$$

since $\overrightarrow{RS}\parallel\overrightarrow{FM}$. On the other hand, $\overrightarrow{RS}\parallel\overrightarrow{GN}$ implies that $\sphericalangle FRS$ is not smaller than $\sphericalangle FGN$. Consequently

$$\sphericalangle FRS=\sphericalangle FGN,$$

and

$$\sphericalangle RFM+\sphericalangle FRS=\sphericalangle GFM+\sphericalangle FGN=2R.$$

Thus

$$\sphericalangle DCP+\sphericalangle CDQ=2R.$$

(II) If $R$ lies on the extension of $\overline{FG}$, then $\sphericalangle NGR=\sphericalangle MFR$. Let $MFGN\equiv$ $\equiv NGHL\equiv LHKO$ and so on until $\overline{FK}=\overline{FR}$ or $\overline{FK}>\overline{FR}$ for the first time. §7 yields $\overrightarrow{KO}\parallel\overrightarrow{HL}\parallel\overrightarrow{FM}$. If $K$ coincides with $R$, then by §1 $\overrightarrow{KO}$ coincides with $\overrightarrow{RS}$ and therefore

$$\sphericalangle RFM+\sphericalangle FRS=\sphericalangle KFM+\sphericalangle FKO=\sphericalangle KFM+\sphericalangle FGN=2R.$$

On the other hand, if $R$ lies within $\overline{HK}$, then according to (I) we have

$$\sphericalangle RHL+\sphericalangle HRS=2R=\sphericalangle RFM+\sphericalangle FRS=\sphericalangle DCP+\sphericalangle CDQ.$$

§14

*If* $\overrightarrow{BN} \| \overrightarrow{AM}, \overrightarrow{CP} \| \overrightarrow{DQ}$ *and* $\sphericalangle BAM + \sphericalangle ABN < 2R$, *then also* $\sphericalangle DCP + \sphericalangle CDQ < 2R$.

For, if $\sphericalangle DCP + \sphericalangle CDQ$ were not less than $2R$, then by **§1** it would be equal to $2R$. By **§13** this would imply

$$\sphericalangle BAM + \sphericalangle ABN = 2R$$

contrary to the assumption.

***Note on 13 and 14***

These two theorems show that **if in a single case the sum of the interior angles formed by a transversal between parallel lines is equal to or less than two right angles, then the sum of the interior angles between parallels will be equal to or less than two right angles in every case**. These theorems were first proved by Saccheri, who also showed that if in a single case the HOA is true, then in it is true in every case.

§15

*In possession of the results of* **§§ 13** *and* **14**, *denote by* $\Sigma$ *the system of geometry based on the hypothesis that* EUCLID'S *Axiom XI is true, and denote by* **S** *the system based on the opposite hypothesis.*

*All theorems we state without expressly specifying the system* $\Sigma$ *or* **S** *in which the theorem is valid are meant to be absolute, that is, valid independently of whether* $\Sigma$ *or* **S** *is true in reality.*

***Note on 15***

Here Bolyai emphasises that **theorems whose proof is independent of the parallel postulate are part of absolute geometry and are therefore true regardless of whether or not space is Euclidean**.

## §16

*If $\overrightarrow{AM}$ is an axis for some* L *then, in system* $\Sigma$, *this* L *is the straight line perpendicular to* $\overrightarrow{AM}$.

For, if $\overrightarrow{BN}$ is an axis of L starting from the point $B$ of L, then in $\Sigma$

$$\sphericalangle BAM + \sphericalangle ABN = 2\sphericalangle BAM = 2R,$$

so that

$$\sphericalangle BAM = R.$$

If $C$ is any point of the line $AB$ and if $\overrightarrow{CP} \| \overrightarrow{AM}$, then by **§13**

$$CP \simeq AM;$$

thus, according to **§11**, $C$ is a point of L.

*In system* S, *however, no three points $A$, $B$, $C$ of* L *or* F *are on a straight line.*

Really, one of the axes $\overrightarrow{AM}$, $\overrightarrow{BN}$, $\overrightarrow{CP}$, say $\overrightarrow{AM}$ lies between the other two, in which case by **§14** both $\sphericalangle BAM$ and $\sphericalangle CAM$ are smaller than $R$.

***Note on 16***

In this theorem Bolyai makes the point that in Euclidean geometry the equivalent of a horocycle is a straight line perpendicular to a parallel axis. **In hyperbolic geometry, however, no three points of a horocycle lie in a straight line.** This follows from the fact that, by Section 14, the angle formed by a transversal at a corresponding point is less than two right angles:

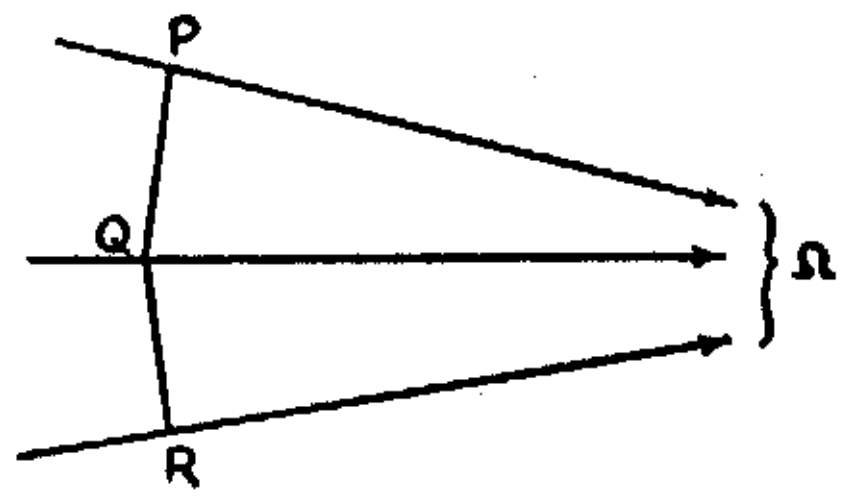

Therefore given three parallel lines with corresponding points P, Q and R, the sum of the adjacent angles RQΩ, PQΩ is less than two right angles. It follows that P, Q and R, through which a horocycle must pass, do not lie in a straight line.

## §17

*In system* S, *too,* L *is a line and* F *is a surface.*

For, in view of **§11**, any plane which is perpendicular to the axis $\overrightarrow{AM}$ and passes through some point of **F** intersects **F** in a circle whose plane, according to **§14**, is not perpendicular to any other axis $\overrightarrow{BN}$.

If we let **F** rotate about the axis $\overrightarrow{BN}$, then by **§12** all points of **F** remain in **F** and the intersection of **F** with any plane not perpendicular to $\overrightarrow{BN}$ describes a surface. By **§12**, however, taking two points $A$, $B$ of **F** we can place **F** congruently on itself so that $A$ falls on $B$. Thus **F** is a uniform surface.

Hence by **§§11—12** it follows that **L** is a uniform curve.*

***Note on 17***

This theorem establishes a number of results on the horosphere.

Firstly, **any plane which is perpendicular to an axis of the horosphere and passes through some point on the surface will intersect the horosphere in a circle.**

Secondly, by the result given in Section 14**, this plane cannot be perpendicular to any other axis of the horosphere** (otherwise the sum of the interior angles formed between two parallel axes would be  equal two right angles, which is not the case in hyperbolic geometry).

Thirdly**, if the horosphere, defined in Section 11 as the set of mutually corresponding points on a bundle of parallel lines, is rotated about any one of its axes, all its points remain on the horosphere**.

Lastly, since only one horocycle passes through any two points A, B on the horosphere (Section 12), the horosphere can be rotated about a suitable axis so as to place A on the position formerly occupied by B. A suitable axis in this case would lie in the plane which perpendicularly bisects the straight line joining A and B. This result implies that **the horosphere is a uniform surface**.

## §18

*In system* **S**, *any plane that passes through the point A of* **F** *and stands obliquely to the axis* $\overrightarrow{AM}$ *intersects* **F** *in a circle.*

For let $A, B, C$ be three points of the intersection, and let $\overrightarrow{BN}$ and $\overrightarrow{CP}$ be axes. The planes $AMBN$, $AMCP$ form an angle; otherwise the plane determined by $A, B, C$ would contain, according to §16, the axis $\overrightarrow{AM}$, which contradicts the assumption. Therefore by §10, the planes perpendicularly bisecting the distances $\overline{AB}$ and $\overline{AC}$ intersect each other in an axis $\overrightarrow{FS}$ of **F**. Thus

$$\overline{FB} = \overline{FA} = \overline{FC}.$$

Further let

$$AH \perp FS.$$

Let the plane $FAH$ rotate about $FS$. Then $A$ describes the circle of radius $\overline{HA}$ passing through $B$ and $C$, and this circle is in both **F** and the plane $ABC$. Moreover, by §16, **F** and $ABC$ have no common points besides those of the circle $\circ\overline{HA}$.

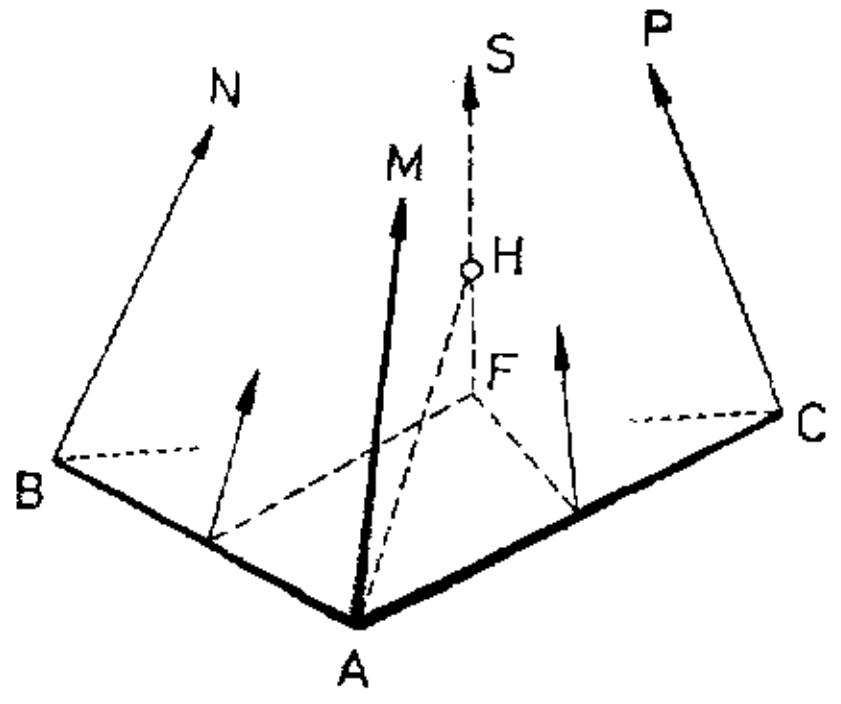

**Figure 12**

It is also clear that if the arc $\overset{\frown}{FA}$ of an **L** rotates as a radius about the point $F$ within **F** then its endpoint describes the circle $\circ\overline{HA}$.

### *Note on 18*

Bolyai proves in this theorem that **a plane passing through a point A on the horosphere intersects the horosphere in a circle providing that the plane does not contain the axis AM**. An important corollary of the theorem is that the circle containing he points A, B, C, which is generated by rotating the straight line HA about the axis FS and the circle generated by rotating the horocyclic arc FA (not shown in the figure) about FS are one and the same.

### §19

*In system* S, *a line BT which is perpendicular to the axis* $\overrightarrow{BN}$ *of some* L *and lies in the plane of* L *is a tangent to* L.

For, by **§14**, $B$ is the only point of $\overrightarrow{BT}$ that lies in L. If, however, $\overrightarrow{BQ}$ lies in $\sphericalangle TNB$, then the plane containing $BQ$ and perpendicular to the plane $TBN$ intersects, by **§18**, the **F** corresponding to the axis $\overrightarrow{BN}$ in a circle whose centre lies manifestly on $\overrightarrow{BQ}$. If $\overline{BQ}$ is a diameter of this circle, then $\overrightarrow{BQ}$ obviously intersects the L corresponding to the axis $\overrightarrow{BN}$ in the point $Q$.

***Note on 19***

Bolyai shows that **a line perpendicular to the axis of a horocycle at the point where it intersects the horosphere, assuming the line lies in the plane of the horocycle, is tangent to the horocycle**. His proof relies on Section 14, which implies that in hyperbolic geometry the sum of interior angles formed when a transversal meets two parallel lines at corresponding points is less than two right angles. Since these angles are by definition equal, each one is less than a right angle. Consequently a line BT which is perpendicular to the axis BN of a horocycle cannot meet any line parallel to BN at a corresponding point and therefore cannot meet the horocycle passing through B a second time. It follows that BT is tangent to the horocycle at B, as required.

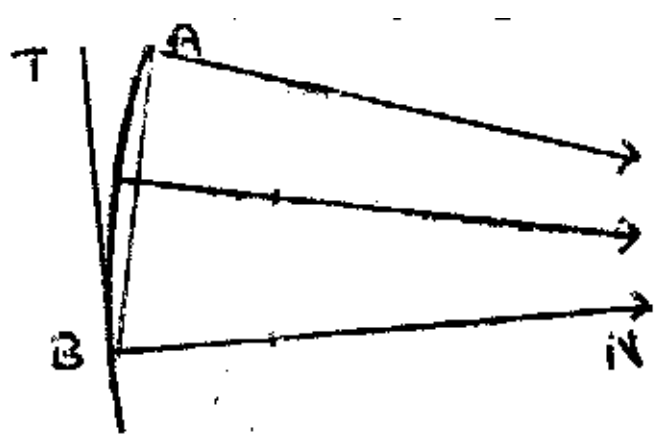

This theorem also implies that the axis BN is perpendicular to the horocycle BA and by extension is perpendicular to the horosphere on whose surface BA lies. To see this, suppose that the chord BA rotates anticlockwise about B. As ∠ABN increases, point A will approach point B and BA will approach the tangent BT. When A coincides with B, ∠ABN will a right angle.

Thus in Sections 11-19 Bolyai has shown that the horosphere and a plane can be related in four ways. Either the plane does not intersect the horosphere, or is tangent to it, or intersects it in a circle, or (in the case of the principal plane containing an axis) in a horocycle.

### §20

*According to* **§§ 11** *and* **18**, *any two points of* **F** *determine an* L*-line of* **F**. *As*, however, by **§§16** and **19** L is perpendicular to each of its axes, *the angle between any two* L*-lines in* **F** *is equal to the dihedral angle between the planes that contain the arms and are perpendicular to* **F**.

***Note on 20***

This theorem states that **the angle between any two horocycles on the horosphere is equal to the dihedral angle between the planes containing the horocycles.** The dihedral angle between two intersecting planes is defined as the angle between the perpendiculars drawn along the planes from any point on their line of intersection. Now by Section 11 each plane formed by a pair of parallel axes intersects the horosphere in a horocycle and by Section 19 this plane is perpendicular to the horosphere. When two planes intersect at a point on the horosphere, the tangents to the horocycles at this point are perpendicular to the line of intersection of the planes. Therefore the angle between the horocycles is equal to the dihedral angle between the planes.

**§21**

*Two* **L**-*lines* $\overrightarrow{AP}$, $\overrightarrow{BQ}$ *that lie in the same* **F** *and whose interior angles with a third* **L**-*line AB have a sum less than* 2*R intersect each other*. (In the surface **F**, *AP* denotes the **L**-line through the points *A* and *P*, while $\overrightarrow{AP}$ denotes that half of it which starts from *A* and contains *P*.)

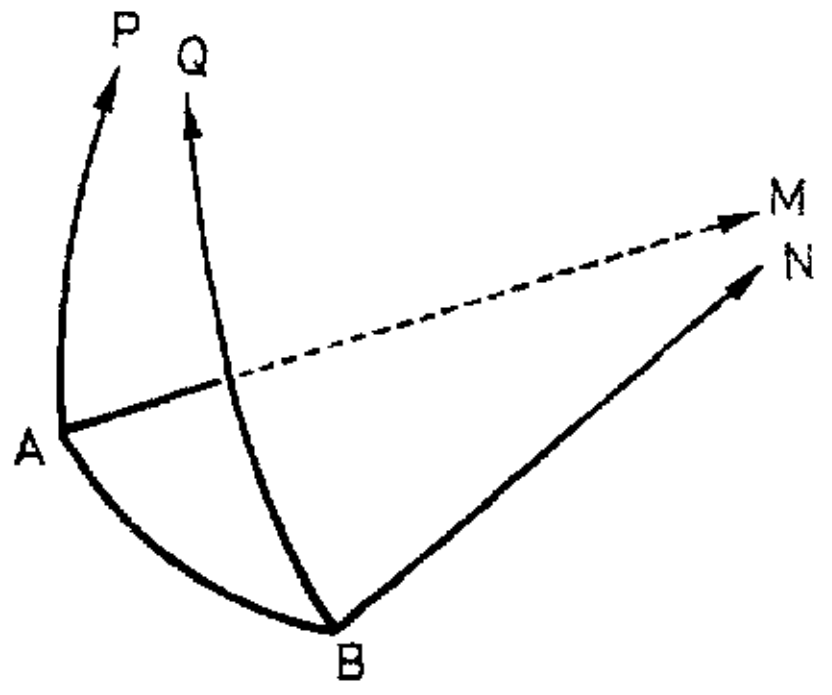

**Figure 13**

For, if $\overrightarrow{AM}$ and $\overrightarrow{BN}$ are axes of **F**, then $|AM|P$ and $|BN|Q$ intersect each other by §9, and **F** intersects their line of intersection by §§ 7 and **11**. Consequently, also $\overrightarrow{AP}$ and $\overrightarrow{BQ}$ intersect each other.

Hence it is clear that *in the surface* **F**, *if we let* **L**-*lines play the role of straight lines*, Axiom XI as well as the whole geometry and trigonometry of the plane are absolutely valid. Thus, based on the foregoing, trigonometrical functions can be defined just as in system $\Sigma$. Also, the circumference of the circle whose radius along an **L**-line of **F** is $r$ equals $2\pi r$ and its area in **F** equals $\pi r^2$ (here $\pi$ denotes half the circumference of the circle of unit radius, that is, the well-known number 3.1415926...).

***Note on 21***

This theorem reveals the astonishing fact that **the geometry of the horosphere is Euclidean**. That is to say, if horocycles play the role of straight lines, the figures formed on the surface of the horosphere satisfy Euclid's five postulates viz: only one horocycle can be drawn between two points, a horocycle can be infinitely extended, a circle of any radius can be swept out by a horocycle and all right angles between horocycles are equal on the horosphere.

Bolyai does not actually mention the first four postulates here, perhaps feeling that he had already provided sufficient proof that they can be applied to the horosphere. The first two follow from his definition of the horocycle in Section 11. The third follows from Section 16. The fourth follows from the fact that angles between horocycles are considered right angles if and only if the plane angles between their tangents are right angles. Since all right angles in the plane are equal by Euclid's fourth postulate, which holds in neutral geometry, all right angles are equal on the horosphere.

What about the parallel postulate? Bolyai provides an explicit proof that this holds true on the horosphere. Suppose that AM is parallel to BN and that AP and BQ are horocycles on the surface of the horosphere. Now by Section 20 the horocyclic angle PAB is equal to the dihedral angle between plane PAM and plane MABN, while the horocyclic angle QBA is equal to the dihedral angle between plane QBN and plane MABN. Suppose that the sum of the interior angles PAB, QBA is less than two right angles. It follows that the sum of the dihedral angles between the planes is less than two right angles. But it is proved in Section 9 that if the sum of the dihedral angles which two planes make with a third plane is less than two right angles, the lines of intersection between the planes being parallel, then the two planes will intersect. Thus planes PAM, QBN meet and by the same token the horocycles AP, AQ lying in those planes must also meet. This proves the parallel postulate in the form given by Euclid, namely that, given two lines, if the sum of the interior angles formed on one side by a transversal is less than two right angles, then the two lines will meet on that side. Only in this case the lines are horocycles rather than straight lines.

It follows that the whole of Euclidean geometry is valid on the horosphere if horocycles are substituted for straight lines. So, for example, the sum of the angles of a horocyclic triangle is two right angles and the area of a circle swept out on the horosphere by a horocyclic radius of length $r$ is given by the Euclidean formula $\pi r^2$. Moreover the trigonometric functions together with all the theorems of plane trigonometry can applied to horocyclic triangles.

**§22**

*If $AB$ is an* **L**-*line corresponding to the axis $\overrightarrow{AM}$ while $C$ lies on $\overrightarrow{AM}$, and $\sphericalangle CAB$ (whose arms are the axis $\overrightarrow{AM}$ and the* **L**-*formed $\overrightarrow{AB}$) moves to infinity first along $\overrightarrow{AB}$ and afterwards along $\overrightarrow{BA}$, then the path $CD$ of $C$ is the* **L**-*line corresponding to the axis $\overrightarrow{CM}$.*

Really, denote by **L**′ the L-line corresponding to $\overrightarrow{CM}$. Let $D$ be any point of the path $CD$. Let

$$\overrightarrow{DN} \parallel \overrightarrow{CM},$$

and $B$ a point of **L** on the line $DN$. Then

$$BN \simeq AM \quad \text{and} \quad \overline{AC} = \overline{BD},$$

whence

$$DN \simeq CM,$$

so that $D$ lies on **L**′. If, however, $D$ lies on **L**′,

$$\overrightarrow{DN} \parallel \overrightarrow{CM}$$

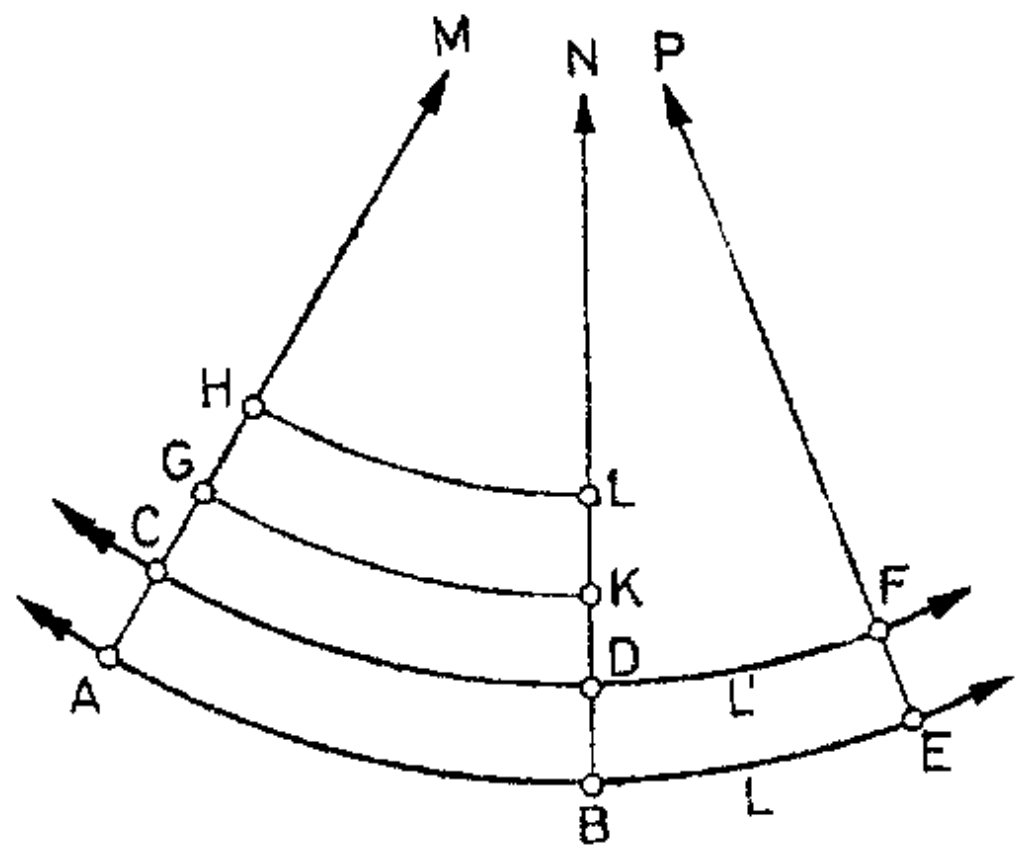

**Figure 14**

and $B$ is a point of both **L** and $DN$, then

$$AM \simeq BN \quad \text{and} \quad CM \simeq DN,$$

which obviously imply that

$$\overline{BD} = \overline{AC},$$

$D$ is on the path of $C$, and **L**′ coincides with $CD$. For such **L**′, we write

$$\mathbf{L}' \parallel \mathbf{L}.$$

***Note on 22***

This theorem establishes the existence of **concentric horocycles**.

## §23

If $CDF$ and $ABE$ are L-formed lines such that by §22

$$CDF \parallel ABE,$$

further if

$$\overset{\frown}{AB} = \overset{\frown}{BE},$$

and $\overrightarrow{AM}$, $\overrightarrow{BN}$, $\overrightarrow{EP}$ are axes, then obviously

$$\overset{\frown}{CD} = \overset{\frown}{DF}.$$

Also, if $A$, $B$, $E$ are three arbitrary points of the line $AB$ and

$$\overset{\frown}{AB} = n \cdot \overset{\frown}{CD},$$

then

$$\overset{\frown}{AE} = n \cdot \overset{\frown}{CF}$$

and therefore

$$\overset{\frown}{AB} : \overset{\frown}{CD} = \overset{\frown}{AE} : \overset{\frown}{CF};$$

the latter relation holds even if $\overset{\frown}{AB}$, $\overset{\frown}{AE}$, $\overset{\frown}{CD}$ are incommensurable. Thus *the ratio* $\overset{\frown}{AB} : \overset{\frown}{CD}$ *is independent of the length of* $\overset{\frown}{AB}$ *and is completely determined by the distance* $\overline{AC}$. *Whenever we denote the length* $\overline{AC}$ *by a small letter* (say $x$), *the ratio* $\overset{\frown}{AB} : \overset{\frown}{CD}$ *will be denoted by the capital letter of the same name* ($X$).

***Note on 23***

Here Bolyai establishes that **the ratio *X* of the lengths of two concentric horocyclic arcs bounded between the same parallels is purely a function of the radial distance *x* between them**.

## §24

*For any x and y, with the notation introduced in* §23,

$$Y = X^{\frac{y}{x}}.$$

For, either is one of the values $x$, $y$ a multiple of the other (say, $y$ a multiple of $x$ or not.

If

$$y = nx,$$

then let

$$x = \overline{AC} = \overline{CG} = \overline{GH}$$

and so on, until

$$\overline{AH} = y.$$

Further let

$$\widehat{CD} \parallel \widehat{GK} \parallel \widehat{HL}.$$

§23 yields

$$X = \widehat{AB}:\widehat{CD} = \widehat{CD}:\widehat{GK} = \widehat{GK}:\widehat{HL}.$$

Hence

$$\frac{\widehat{AB}}{\widehat{HL}} = \left(\frac{\widehat{AB}}{\widehat{CD}}\right)^n,$$

that is,

$$Y = X^n = X^{\frac{y}{x}}.$$

If $x$ and $y$ are multiples of one and the same $z$, say

$$x = mz, \quad y = nz,$$

then by the foregoing

$$X = Z^m, \quad Y = Z^n$$

and, consequently,

$$Y = X^{\frac{n}{m}} = X^{\frac{y}{x}}.$$

The result can be easily extended to the case where $x$ and $y$ are incommensurable.

If $q = y - x$, then obviously $Q = Y:X$.

It is also clear that in system $\Sigma$ for any $x$ we have

$$X = 1.$$

In system S, however, we always have

$$X > 1,$$

and for any arcs $\widehat{AB}$, $\widehat{ABE}$ there is a $\widehat{CDF}$ such that

$$\widehat{CDF} \parallel \widehat{ABE} \quad \text{and} \quad \widehat{CDF} = \widehat{AB};$$

hence

$$(AM, BN) \equiv (AM, EP),$$

although the latter can be any multiple of the former. However strange this result should be, it obviously still does not prove the absurdity of system S.

***Note on 24***

In this section Bolyai derives the equation expressing the ratio of the arc lengths of concentric horocycles lying between the same parallels in terms of their radial distance from a given horocycle. If $X$ represents the ratio between a given arc and another arc separated by a radial distance $x$, while $Y$ represents the ratio between the given arc and a third arc separated by a radial distance $y$ (the two arcs lying in the direction of parallelism with respect to the given arc), Bolyai's theorem states that

$$Y = X^{y/x}.$$

His proof runs as follows. Suppose, with reference to Figure 14, that AB is the given arc, AC = CG = GH = $x$, while AH = $y$ where $y = nx$. If

$$\text{arc AB} : \text{arc CD} = \text{arc CD} : \text{arc GK} = \text{arc GK} : \text{arc HL} = X$$

while

$$\text{arc AB} : \text{arc HL} = Y,$$

then

$$\text{arc AB} / \text{arc HL} = (\text{arc AB} / \text{arc CD}) \,.\, (\text{arc CD} / \text{arc GK}) \,.\, (\text{arc GK} / \text{arc HL})$$

$$= (\text{arc AB}/ \text{arc CD})^n.$$

In other words

$$Y = X^n$$

$$= X^{y/x}.$$

Furthermore if AH = z, AG = $v$ and GH = $y$, so that $y = z - v$, then

$$\text{arc GK} / \text{arc HL} = (\text{arc AB} / \text{arc HL}) / (\text{arc AB} / \text{arc HL}),$$

so

$$Y = Z : V,$$

a result which Bolyai will need in Section 29.

Bolyai makes the point that these identities hold true even if $x/y$ is an irrational number, a result which can be proved by taking limits.

This theorem implies the existence of a linear constant in hyperbolic geometry, though Bolyai only introduces this constant explicitly in Section 30. Suppose that arc AB : arc CD = $a$ for some $a > 1$ and let AC = 1 (the length of AC being the unit of distance corresponding to $a$). Suppose that HK is a concentric arc where AH = $x$. Then the theorem implies that arc AB : arc HK = $a^x$ . Denote arc AB by $s$ and arc HK by $s_x$. Then

$$s / s_x = a^x.$$

If now $k$ is chosen so that $a = e^{1/k}$, then

$$s / s_x = e^{x/k}.$$

This is the origin of the constant $k$ which appears in all the formulas of hyperbolic geometry. Putting $x = k$ in the previous equation gives

$$s / s_k = e.$$

Thus the linear constant $k$ (which is a physical distance of unknown magnitude, not just a number) can be defined as the radial distance separating two concentric horocyclic arcs whose lengths are in the ratio $e$ : 1.

It follows that Bolyai's $X$ and $Y$ can be written, respectively, as $e^{x/k}$ and $e^{y/k}$.

In the last part of this section Bolyai asserts that in Euclidean geometry $X = 1$ (since the distance between parallels is constant), whereas in hyperbolic geometry $X > 1$ (since parallels approach one another asymptotically). This fact implies one of the (many) paradoxes of hyperbolic geometry, which is illustrated in the figure below:

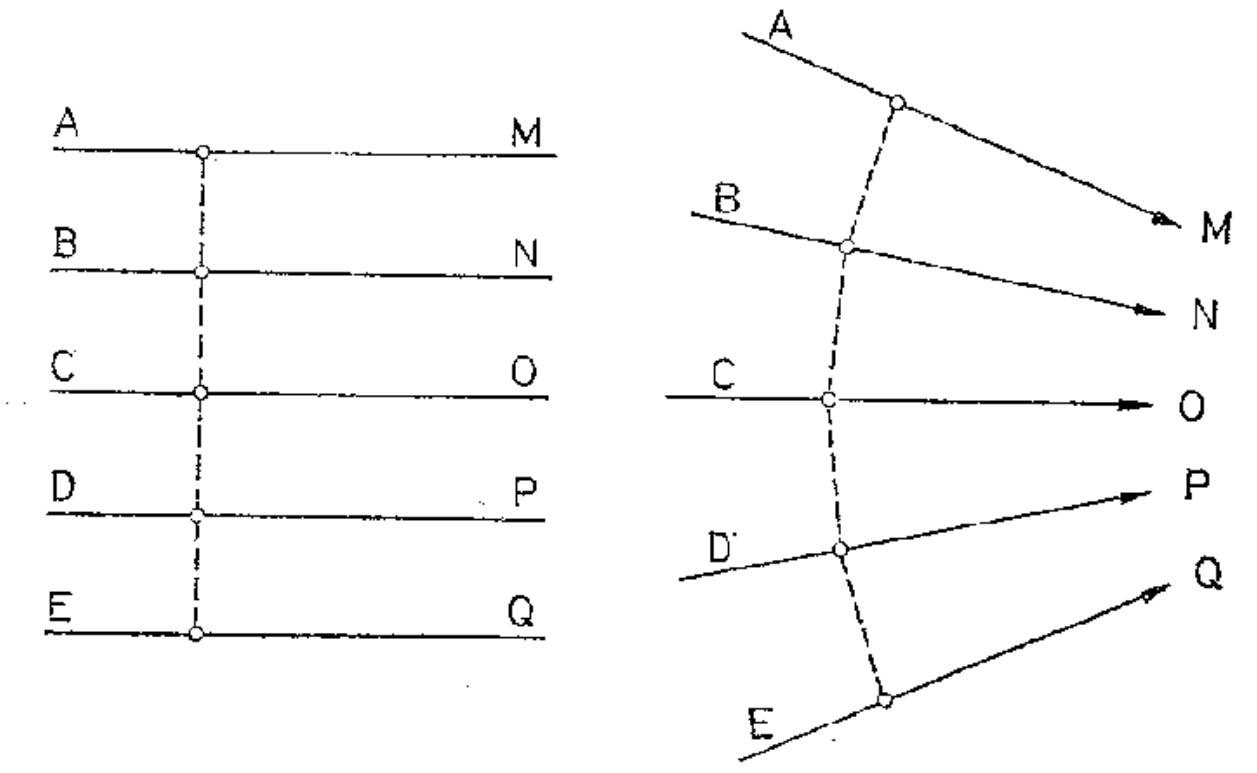

Bolyai makes the point that in Euclidean geometry a strip such as (AM, BN) is clearly not congruent to the strip (AM, EQ), which is many times greater, whereas in hyperbolic geometry all strips are congruent. For example, if the strip (AM, BN) is extended to the left, an arc AB can be found equal to arc AE so that (AM, BN) can be superimposed upon (AM, EQ). The paradox is the fact that in hyperbolic geometry the part is congruent with the whole, though, as Bolyai states, this does not prove that hyperbolic geometry is absurd.

## §25

*In any rectilinear triangle, the circumferences of the circles of radii equal to the sides are to each other as are the sines of the angles opposite to them.*

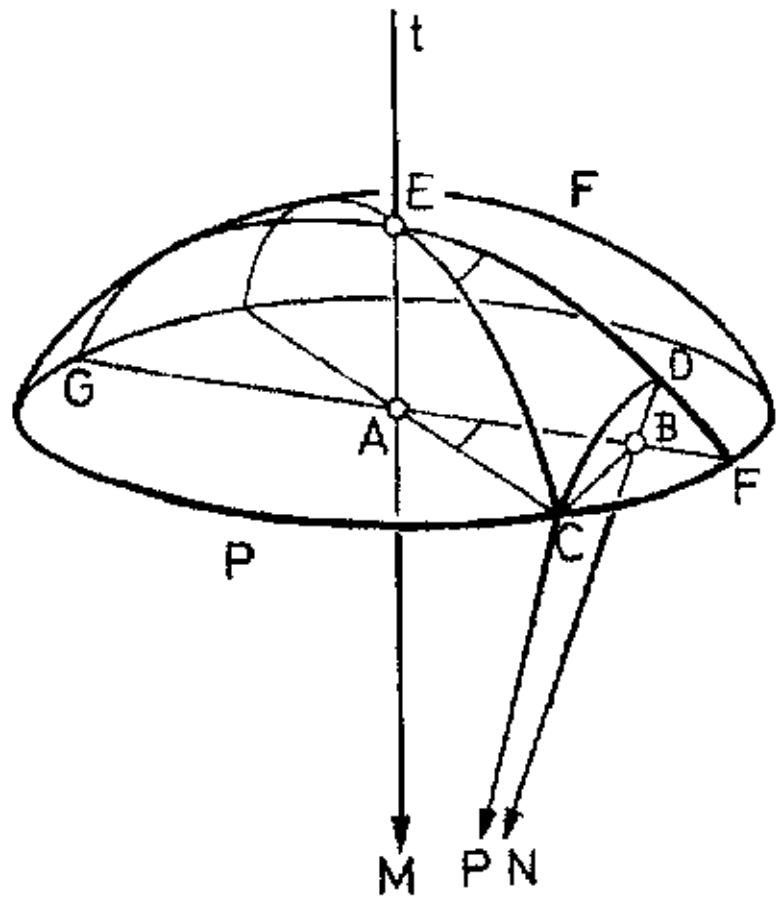

**Figure 15**

For let

$$\sphericalangle ABC = R, \quad \overrightarrow{AM} \perp BAC$$

and

$$\overrightarrow{BN} \parallel \overrightarrow{AM}, \quad \overrightarrow{CP} \parallel \overrightarrow{AM}.$$

Then the plane $CAB$ is perpendicular to the plane $AMBN$ and, since $CB \perp BA$, also $CB$ is perpendicular to $AMBN$, so that the planes $CPBN$, $AMBN$ are perpendicular to each other. Let the **F** corresponding to the axis $\overrightarrow{CP}$ intersect the lines $BN$, $AM$ in the points $D$, $E$ and the strips $(CP, BN)$, $(CP, AM)$, in the L-formed arcs $\widehat{CD}$, $\widehat{CE}$, $\widehat{DE}$, $(BN, AM)$ respectively. Then, by §20, $\sphericalangle CDE$ will be equal to the dihedral angle between $NDC$ and $NDE$, that is $R$, and a similar argument yields

$$\sphericalangle CED = \sphericalangle CAB.$$

According to §21, in the triangle $CED$ formed of L-lines we have*

$$\widehat{EC} : \widehat{DC} = 1 : \sin DEC = 1 : \sin CAB.$$

Again by §21, for the circles drawn on **F**

$$\widehat{EC} : \widehat{DC} = \circ\widehat{EC} : \circ\widehat{DC},$$

and by §18

$$\circ\widehat{EC} : \circ\widehat{DC} = \circ\overline{AC} : \circ\overline{BC}.$$

Thus we also have

$$\circ\overline{AC} : \circ\overline{BC} = 1 : \sin CAB,$$

from which the validity of the assertion for any triangle follows.

***Note on 25***

In this theorem Bolyai exploits the Euclidean properties of the horosphere to derive a neutral version of the sine law. The Euclidean sine law states that the sides of a triangle are in the same ratio as the sines of the opposite angles. In contrast, **the neutral sine law (valid both in Euclidean and hyperbolic geometry) states that in any triangle the circumference of the circles whose radii are equal to the sides are in the same ratio as the sines of the opposite angles.**

He proves in respect of the right angled triangle ABC that

$$\boldsymbol{{}^{\circ}AB : {}^{\circ}BC = 1 : \sin \angle CAB},$$

where ${}^{\circ}$AB, ${}^{\circ}$ BC denote, respectively, the circumference of the circles with radius AB, BC.

Bolyai's proof depends on three important facts. Firstly, the circle generated by the horocyclic arc EC as it rotates about the axis EM and the circle generated by the segment AC as it rotates about the same axis are one and the same. (AC and EC are related by the fact that AC is equal to half the length of the chord joining the endpoints of the horocyclic arc whose length is twice the length of arc EC, an important relationship which Bolyai uses in Section 28). Secondly, ∠ECF = ∠CAB since both angles are equal to the dihedral angle between planes MECP, MEDN. Thirdly, ∠EDC = ∠ABC = *R* by construction, so that *sin R* = 1.

Then

$$
\begin{aligned}
1 : \sin \angle CAB &= 1 : \sin \angle CED \\
&= \text{arc } EC : \text{arc } DC \quad \text{(by the Euclidean sine rule valid on the horosphere)} \\
&= {}^{\circ}\text{arc } EC : {}^{\circ}\text{arc } DC \quad \text{(since circles are in proportion to their radii)} \\
&= {}^{\circ} AB : {}^{\circ} BC ,
\end{aligned}
$$

since ${}^{\circ}$arc EC and ${}^{\circ}$ AB are the same circle, while a similar argument applied to a slightly different configuration would imply that ${}^{\circ}$arc DC and ${}^{\circ}$ BC are also the same circle. Bolyai makes the point that this theorem, having been proved in respect of a right angled triangle, can be extended to any triangle.

Bolyai goes on to make repeated use of the neutral sine rule to derive the trigonometrical identities of hyperbolic geometry. In his hands the rule allows the trigonometrical functions of the Euclidean plane to be extended to the hyperbolic plane.

## §26

*In any spherical triangle, the sines of the sides are to one another as are the sines of the angles opposite to them.*

For let

$$\sphericalangle ABC = R,$$

and let the plane $CED$ be perpendicular to the radius $OA$ of the sphere. Then the plane $CED$ is perpendicular to the plane $AOB$ and, as also $BOC$ is perpendicular to $BOA$,

$$CD \perp OB.$$

However, according to §25, in the triangles $CEO$ and $CDO$ we have

$$\circ\overline{EC} : \circ\overline{OC} : \circ\overline{DC} = \sin \widehat{AC} : 1 : \sin \widehat{BC} = \sin COE : 1 : \sin COD,$$

and again by §25

$$\circ\overline{EC} : \circ\overline{DC} = \sin CDE : \sin CED.$$

Consequently,

$$\sin \widehat{AC} : \sin \widehat{BC} = \sin CDE : \sin CED.$$

But

$$\sphericalangle CDE = R = \sphericalangle CBA$$

and

$$\sphericalangle CED = \sphericalangle CAB,$$

so that

$$\sin \widehat{AC} : \sin \widehat{BC} = 1 : \sin A.$$

*Thus the spherical trigonometry which can hence be deduced has obtained a foundation independent of Axiom XI.*

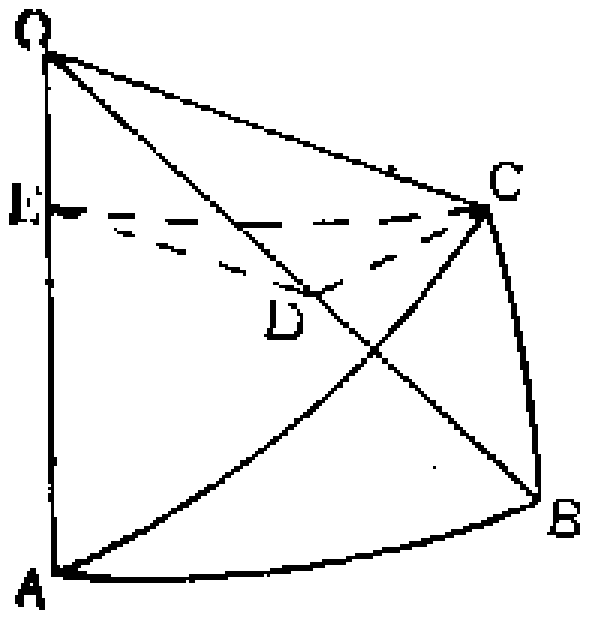

**Figure 16**

***Note on 26***

In this section Bolyai derives the spherical law of sines for a right angled spherical triangle using only neutral geometry. This proves that **spherical trigonometry, in contrast to plane trigonometry, is independent of the parallel postulate.** Given a spherical triangle ABC with a right angle at B, whose sides *a, b, c* are arcs of great circles on a sphere centre O, he shows that

***sin ∠*A *= sin a / sin b*.**

In Figure 16 the plane triangle CED is drawn perpendicular to the radius of the sphere OA. Then the planes CED and BOC are both perpendicular to the plane AOB. It follows that triangle CEO has a right angle at E and triangle CED has a right angle at D.

The proof runs as follows. In triangles CED, CEO the neutral sine rule implies that

$$^{\circ}EC : {}^{\circ}OC : {}^{\circ}DC = \mathit{sin}\ \angle COE : 1 : \mathit{sin}\ \angle COD$$

$$= \mathit{sin}\ \text{arc AC} : 1 : \mathit{sin}\ \text{arc BC},$$

where the sine of an arc denotes the sine of the angle subtended by the arc at the centre of the sphere.

But in triangle CED alone the neutral sine rule implies that

$$^{\circ}EC : {}^{\circ}DC = \mathit{sin}\ \angle CDE : \mathit{sin}\ \angle CED$$

Consequently

$$\mathit{sin}\ \text{arc AC} : \mathit{sin}\ \text{arc BC} = \mathit{sin}\ \angle CDE : \mathit{sin}\ \angle CED.$$

But

$$\angle CDE = R = \angle CBA.$$

Furthermore, since $\angle$CED and $\angle$CAB are both equal to the dihedral angle between the planes OAC and OAB,

$$\angle CED = \angle CAB = \angle A.$$

Hence

$$\mathit{sin}\ \text{arc AC} : \mathit{sin}\ \text{arc BC} = 1 : \mathit{sin}\ \angle A,$$

or

$$\mathit{sin}\ \angle A = \mathit{sin}\ a\ /\ \mathit{sin}\ b,$$

which is the required sine law for a right angled spherical triangle.

## §27

*If AC and BD are perpendicular to AB and $\sphericalangle CAB$ moves along the line AB, then denoting the path of the point C by $\overset{\frown}{CD}$ we have*

$$\overset{\frown}{CD} : \overline{AB} = \sin u : \sin v.$$

In fact, let

$$DE \perp CA.$$

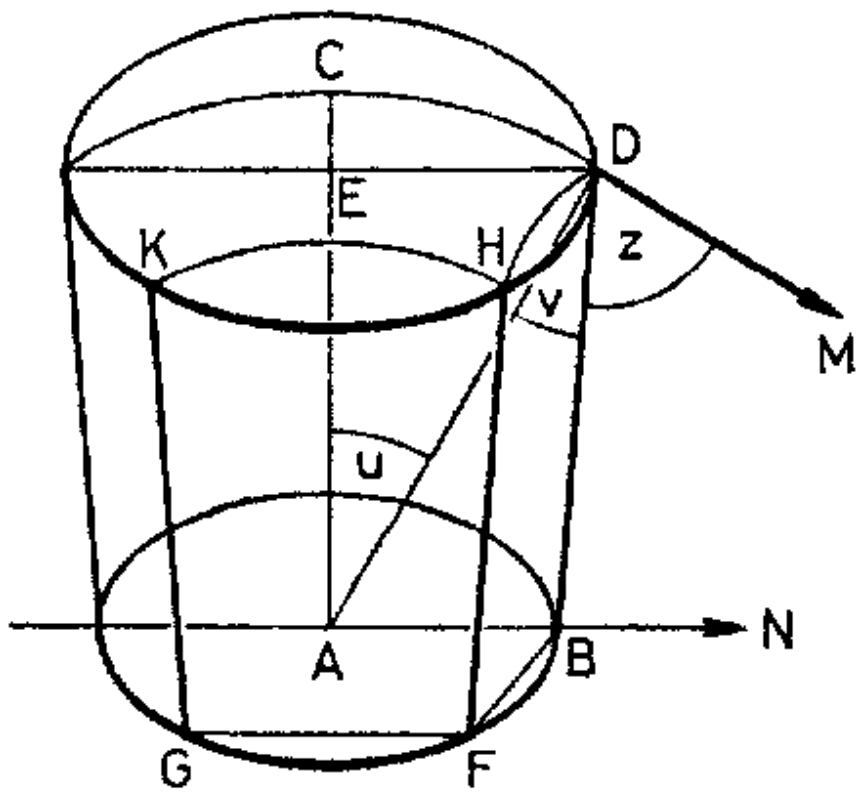

**Figure 17**

For the triangles *ADE, ADB* §25 yields

$$\circ\overline{ED} : \circ\overline{AD} : \circ\overline{AB} = \sin u : 1 : \sin v.$$

If the domain *BACD* rotates about *AC*, then *B* and *D* describe $\circ\overline{AB}$ and $\circ\overline{ED}$, respectively; the surface described by the arc $\overset{\frown}{CD}$ mentioned above will be denoted by $\odot\overset{\frown}{CD}$. Further, let *BFG*... be a polygon inscribed in $\odot\overline{AB}$. The planes containing the sides $BF$, $FG$, ..., respectively, and perpendicular to the plane of $\odot \overline{AB}$ intersect $\odot \overset{\frown}{CD}$ in a curvilinear polygon with the same number of sides. Similarly to §23, it can be proved that

$$\overset{\frown}{CD}:\overline{AB} = \overset{\frown}{DH}:\overline{BF} = \overset{\frown}{HK}:\overline{FG} = \ldots,$$

which yields

$$(\overset{\frown}{DH}+\overset{\frown}{HK}+\ldots):(\overline{BF}+\overline{FG}+\ldots) = \overset{\frown}{CD}:\overline{AB}.$$

If each of the sides $\overline{BF}$, $\overline{FG}$, ... tends to the limit 0, then obviously

$$\overline{BF}+\overline{FG}+\ldots \to \circ\overline{AB}$$

and

$$\overset{\frown}{DH}+\overset{\frown}{HK}+\ldots \to \circ\overline{ED}.$$

Thus, also,

$$\circ\overline{ED}:\circ\overline{AB} = \overset{\frown}{CD}:\overline{AB}.$$

But, as we have just seen,

$$\circ\overline{ED}:\circ\overline{AB} = \sin u:\sin v.$$

Consequently,

$$\overset{\frown}{CD}:\overline{AB} = \sin u:\sin v.$$

If we remove $AC$ from $BD$ to infinity, then the ratio $\overset{\frown}{CD}:\overline{AB}$, that is $\sin u:\sin v$, remains constant. By §1, however,

$$u \to R,$$

and if $\overrightarrow{DM} \| \overrightarrow{BN}$ then

$$v \to z.$$

Hence

$$\overset{\frown}{CD}:\overline{AB} = 1:\sin z.$$

For the path $\overset{\frown}{CD}$ in question we write

$$\overset{\frown}{CD} \| \overline{AB}.$$

***Note on 27***

In this theorem Bolyai obtains the rectification of the equidistant curve. In Euclidean geometry the locus of points which all have the same perpendicular distance from a straight line (and are on the same side) is itself a straight line, but in hyperbolic geometry the equidistant to a straight line is a curve.

In Figure 17 the given straight line is AB. Suppose that AC is drawn perpendicular to AB. If the right angle ∠CAB moves along AB until A coincides with B, the point C will trace out an equidistant curve CD. Bolyai's theorem states that

**arc CD : AB = 1 : *sin z*,**

where *z* is the angle of parallelism corresponding to the segment BD.

The neutral sine rule applied to triangles AED, ADB implies that

$${}^{\circ}ED : {}^{\circ}AD : {}^{\circ}AB = \sin u : 1 : \sin v.$$

If the domain BACD is rotated about AC, then B and D will trace out $^{\circ}$AB and $^{\circ}$ED respectively, while arc CD will sweep out an equidistant surface. Let BFG … be a polygon inscribed in $^{\circ}$AB. The planes containing the sides BF, FG, … and perpendicular to the plane of $^{\circ}$AB will intersect the equidistant surface swept out by arc CD in a curvilinear polygon with the same number of sides. It can be shown that

$$\text{arc CD : AB = arc DH : BF = arc HK : FG = … ,}$$

$$\text{= (arc DH + arc HK + … ) : (BF + FG + …).}$$

If each of the sides BF, FG tends to zero, then

$$\text{BF + FG + …} \rightarrow {}^{\circ}\text{AB}$$

and

$$\text{arc DH + arc HK + …} \rightarrow {}^{\circ}\text{ED},$$

so

$$^{\circ}\text{ED : }^{\circ}\text{AB = arc CD : AB.}$$

But

$$^{\circ}\text{ED : }^{\circ}\text{AB} = \sin u : \sin v.$$

Consequently

$$\text{arc CD : AB} = \sin u : \sin v.$$

However Bolyai goes on to express this ratio in terms of the angle of parallelism corresponding to BD.

If AC recedes from BD to an infinitely distant point, the ratio arc CD : AB, that is $\sin u : \sin v$, remains constant. But $\angle$DAE $\rightarrow$ $\angle$BAE, so $u \rightarrow R$. At the same time AD becomes the asymptotic left hand parallel to AB so that $v$ becomes the angle of parallelism corresponding to the segment BD. Now if DM || BN, $\angle$BDM = $z$ is the right hand angle of parallelism corresponding to BD. Given that the left hand and right hand angles of parallelism are equal, it follows that $v \rightarrow z$. Hence the length of the equidistant arc CD is given by

$$\text{arc CD : AB} = 1 : \sin z.$$

In Section 31 Bolyai makes use of an alternative form of this formula,

**arc CD = AB $\cosh$ (BD/$k$)**.

This form can be justified by invoking the identity obtained by Bolyai in Section 29 on the angle of parallelism. The identity implies that the angle of parallelism $z$ is related to its corresponding segment BD by the equation $\cot \tfrac{1}{2} z = e^{\mathrm{BD}/k}$. Then

$$\sin z = 2 \sin \tfrac{1}{2} z \cos \tfrac{1}{2} z \,/\, (\sin^2 \tfrac{1}{2} z + \cos^2 \tfrac{1}{2} z)$$

$$= 2 \cot \tfrac{1}{2} z \,/\, (\cot^2 \tfrac{1}{2} z + 1)$$ (dividing numerator and denominator by $\sin^2 \tfrac{1}{2} z$)

$$= 2\, e^{\mathrm{BD}/k} \,/\, (e^{2\,\mathrm{BD}/k} + 1)$$

$$= 2/(e^{\mathrm{BD}/k} + e^{-\mathrm{BD}/k})$$

$$= 1/\cosh(\mathrm{BD}/k),$$

whence

arc CD = AB $\cosh$ (BD/$k$).

In modern terms, if $s$ denotes the length of an arc of an equidistant curve, $a$ its projection on the base line and $b$ the common distance of all points of the curve from the base line, then

$$\boldsymbol{s = a \cosh b/k.}$$

## §28

*If $\overrightarrow{BN} \| \simeq \overrightarrow{AM}$ and C lies on $\overrightarrow{AM}$, whereas $\overline{AC} = x$, then with the notation of* **§23** *we have*

$$X = \sin u : \sin v.$$

For, if $CD$ and $AE$ are perpendicular to $BN$, and $BF$ is perpendicular to $AM$, then it can be seen as in **§27** that

$$\circ \overline{BF} : \circ \overline{CD} = \sin u : \sin v.$$

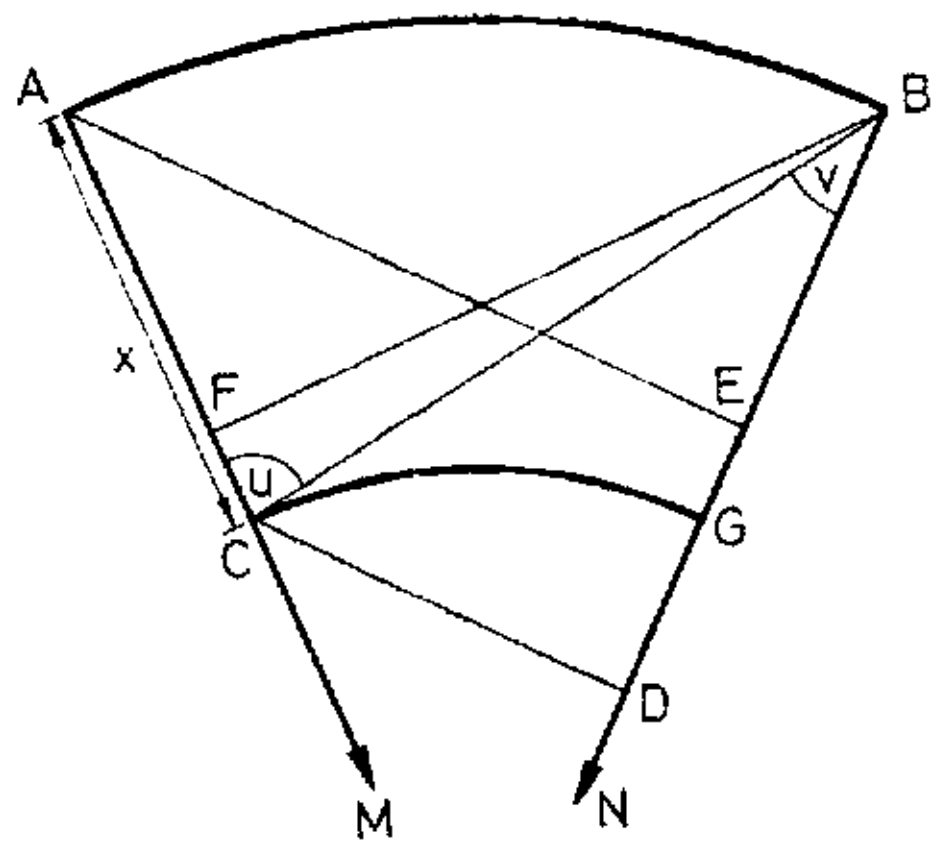

**Figure 18**

But, obviously, $\overline{BF} \simeq \overline{AE}$. So

$$\circ \overline{EA} : \circ \overline{DC} = \sin u : \sin v.$$

On the other hand, in the **F**-surfaces corresponding to $\overrightarrow{AM}$ and $\overrightarrow{CM}$ (that intersect the strip $(AM, BN)$ in the arcs $\widehat{AB}$ and $\widehat{CG}$) **§21** yields

$$\circ \overline{EA} : \circ \overline{DC} = \widehat{AB} : \widehat{CG} = X.$$

As a result,

$$X = \sin u : \sin v.$$

***Note on 28***

Section 28 provides a lemma which Bolyai subsequently uses to derive the identity relating the angle of parallelism to its corresponding segment.

In Figure 18 AM || BN, while AB and CD are concentric horocyclic arcs separated by a radial distance AC = *x*. BF is drawn perpendicular to AM, while AE and CD are drawn perpendicular to BN. ∠ACB = *u* and ∠CBG = *v*.

The lemma states that

**arc AB : arc CG = *X***

***= sin u : sin v.***

Bolyai’s proof runs as follows. The neutral sine rule applied to triangles BFC, BCD implies that

$$\sin u : {}^{\circ}\mathrm{BF} = 1 : {}^{\circ}\mathrm{BC}$$

while

$$\sin v : {}^{\circ}\mathrm{CD} = 1 : {}^{\circ}\mathrm{BC},$$

hence

$$\sin u : \sin v : {}^{\circ}\mathrm{BF} : {}^{\circ}\mathrm{CD}.$$

But by symmetry AE = BF, so

$$\sin u : \sin v = {}^{\circ}\mathrm{AE} : {}^{\circ}\mathrm{CD}.$$

Now in the light of Section 25 AE can be regarded as the semichord corresponding to the horocyclic arc AB, so that ${}^{\circ}$AE = ${}^{\circ}$arc AB, while a similar argument implies that ${}^{\circ}$CD = ${}^{\circ}$arc CG. Therefore

$$X = \text{arc AB} : \text{arc CG}$$

$$= \sin u : \sin v.$$

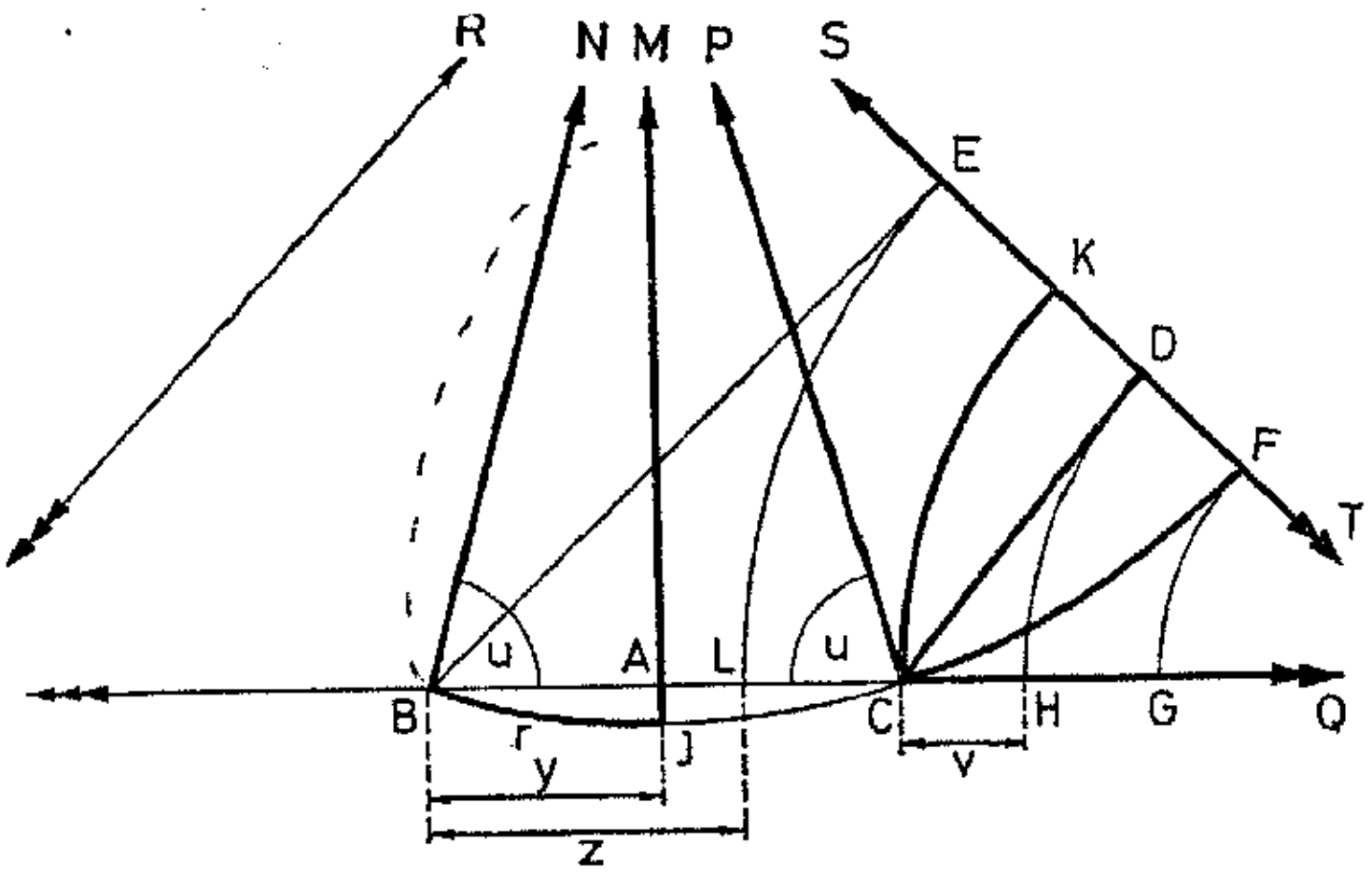

**Figure 19**

## §29

*If* $\sphericalangle BAM = R$, $\overline{AB} = y$ *and* $\overrightarrow{BN} \parallel \overrightarrow{AM}$, *then in system* S

$$Y = \operatorname{ctg} \frac{1}{2} u.$$

For, if

$$\overline{AB} = \overline{AC} \quad \text{and} \quad \overrightarrow{CP} \parallel \overrightarrow{AM}$$

so that

$$BN \simeq CP,$$

further if

$$\sphericalangle PCD = \sphericalangle QCD,$$

then by §**19** there is a $\overrightarrow{DS}$ perpendicular to $\overrightarrow{CD}$ and parallel to $\overrightarrow{CP}$. Hence by §**1**

$$\overrightarrow{DT} \parallel \overrightarrow{CQ}.$$

Moreover, let $BE$ be perpendicular to $\overrightarrow{DS}$. §7 yields $\overrightarrow{DS} \parallel \overrightarrow{BN}$; thus by §6 $\overrightarrow{BN} \parallel \overrightarrow{ES}$.

As, however, $\overrightarrow{DT} \parallel \overrightarrow{CQ}$, we have $\overrightarrow{BQ} \parallel \overrightarrow{ET}$. Consequently, in view of §1,

$$\sphericalangle EBN = \sphericalangle EBQ.$$

Let the arc $\widehat{BCF}$ be contained in the **L**-line corresponding to the axis $\overrightarrow{BN}$, and let $\widehat{FG}$, $\widehat{DH}$, $\widehat{CK}$ and $\widehat{EL}$ be some arcs of the **L**-lines corresponding to the axes $\overrightarrow{FT}$, $\overrightarrow{DT}$, $\overrightarrow{CQ}$ and $\overrightarrow{ET}$, respectively. From **§22** it is clear that

$$\overline{HG} = \overline{DF} = \overline{DK} = \overline{HC}.$$

Hence

$$\overline{CG} = 2 \cdot \overline{CH} = 2v.$$

The relation

$$\overline{BG} = 2 \cdot \overline{BL} = 2z$$

can be verified similarly. But

$$\overline{BC} = \overline{BG} - \overline{CG}.$$

Therefore

$$y = z - v.$$

Consequently, according to **§24**,

$$Y = Z : V.$$

Finally, in view of **§28**,

$$Z = 1 : \sin \frac{1}{2} u$$

and

$$V = 1 : \sin \left( R - \frac{1}{2} u \right).$$

Thus

$$Y = \operatorname{ctg} \frac{1}{2} u.$$

***Note on 29***

In this theorem Bolyai derives the fundamental identity of hyperbolic trigonometry, expressing the functional relationship between the angle of parallelism and its corresponding segment as

$$Y = \cot \tfrac{1}{2} u,$$

where *u* is the angle of parallelism corresponding to a segment of length *y* and *Y* is the ratio of two concentric horocyclic arcs separated by a radial distance *y* i.e. $Y = e^{y/k}$ (see Section 24).

This identity is usually written as

$$\tan \tfrac{1}{2} \Pi (y) = e^{-y/k},$$

where $\Pi (y)$ (in Lobachevsky's notation) denotes the angle of parallelism corresponding to a segment of length *y.*

In Figure 19 BN is parallel to AM while AM is perpendicular to BC, so ∠NBA = $u$ is the angle of parallelism corresponding to the segment AB = $y$. AC = AB while CP is parallel to AM, so ∠PCA = $u$ as well.

DS is drawn perpendicular to CD and parallel to CP. If ∠PCD = ∠QCD, being the angle of parallelism corresponding to CD, it follows that DT is parallel to BQ.

BE is drawn perpendicular to DS, so ∠EBN = ∠EBQ = ½ $u$, being the angle of parallelism corresponding to BE.

A horocycle centred at an infinitely distant point in the direction of N is drawn through B, C and F, while hocycles centred at an infinitely distant point in the direction of Q are drawn through G and F, through H and D, through C and K, through L and E and finally through B (dotted). By symmetry FD = DK. (The symmetry comes from the fact that the points F and K are defined on ST by horocycles centred in opposite directions passing through the same point C. In fact it can be shown that FD = DK = *log cosh* CD). This in turn implies that GH = HC, so CG = 2HC = 2$v$. A similar argument implies that BG = 2BL = 2$z$.

But BC = BG − CG, so 2$y$ = 2$z$ − 2$v$ or $y = z - v$. Consequently, by Section 24, $Y = Z : V$.

Now considering the horocycle LE and the horocycle through B, the lemma of Section 28 implies that

$$Z = \sin \angle \mathrm{BES} : \sin \mathrm{EBL}$$
$$= 1 : \sin \tfrac{1}{2} u.$$

But considering the horocycles HD and CK, the lemma implies that

$$V = \sin \angle \mathrm{CDK} : \sin \angle \mathrm{DCH}$$
$$= 1 : \sin (R - \tfrac{1}{2} u) \quad (\text{since } \angle \mathrm{DCH} = \tfrac{1}{2} \angle \mathrm{PCG} = \tfrac{1}{2} (2R - u))$$
$$= 1 : \cos \tfrac{1}{2} u.$$

Therefore

$$Y = Z : V$$
$$= (1/\sin \tfrac{1}{2} u) / (1/\cos \tfrac{1}{2} u)$$
$$= \cot \tfrac{1}{2} u$$

as required.

Bolyai also employs an alternative form of this identity, which he expresses as

$$\boldsymbol{\cot u = \tfrac{1}{2}\,(Y - Y^{-1})}$$

i.e. $\cot u = \sinh y$.

This form can be justified as follows. Starting from

$$\cot \tfrac{1}{2} u = Y = e^{y/k}$$

it follows that

$$u = 2 \cot^{-1} (e^{y/k}).$$

If $\cot^{-1} (e^{y/k})$ is written as $q$, then

$$\cot u = \cot 2q$$

$$= (\cot^2 q - 1) / 2 \cot q$$ (using the cotangent double angle formula),

whence

$$\cot u = \{\cot^2 [\cot^{-1} (e^{y/k})] - 1\} / 2 \cot [\cot^{-1} (e^{y/k})]$$

$$= (e^{2y/k} - 1) / 2\, e^{y/k}$$

$$= \tfrac{1}{2} (e^{y/k} - e^{-y/k})$$ (multiplying numerator and denominator by $e^{-y/k}$)

$$= \sinh y/k.$$

## §30

Making use of **§25** it is easy to see that, for solving the problem of plane trigonometry in system S, the expression for the circumference of the circle in terms of the radius is needed. This expression, in turn, can be obtained by the rectification of the L-line.

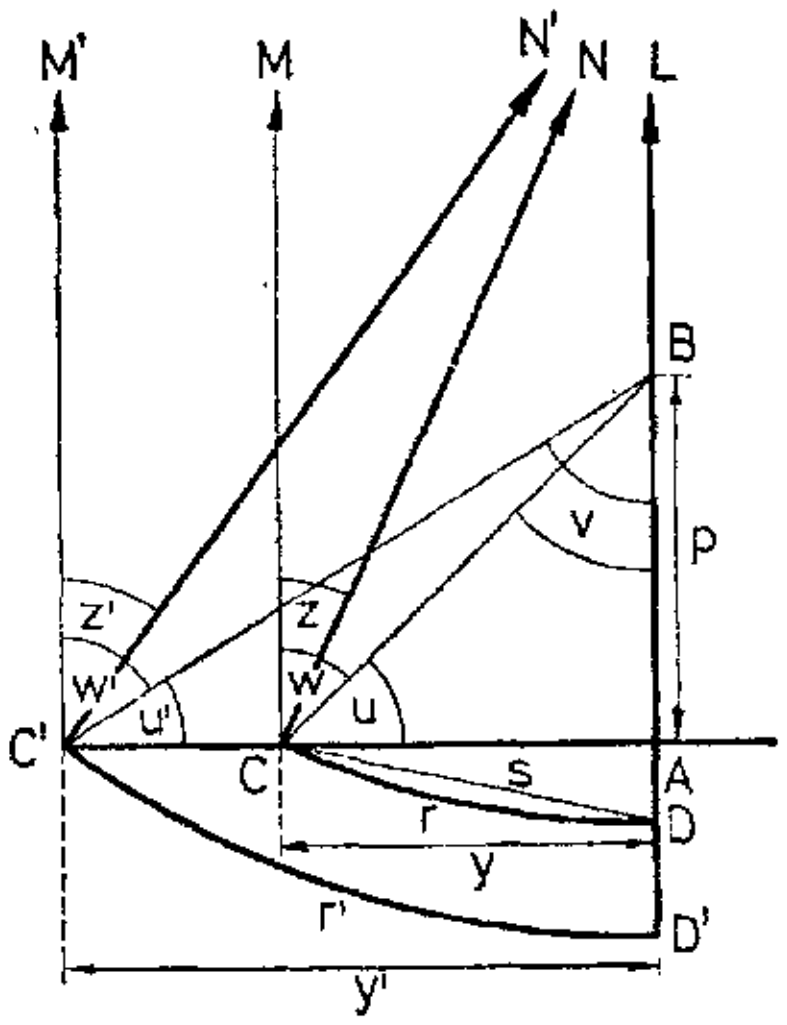

**Figure 20**

Let $\overrightarrow{AB}$, $\overrightarrow{CM}$, $\overrightarrow{C'M'}$ be perpendicular to $\overrightarrow{AC}$, and let $B$ be any point of $\overrightarrow{AB}$. By **§25**

$$\sin u : \sin v = \circ p : \circ y$$

and

$$\sin u' : \sin v' = \circ p : \circ y'.$$

Consequently,

$$\frac{\sin u}{\sin v} \cdot \circ y = \frac{\sin u'}{\sin v'} \cdot \circ y'.$$

But according to **§27**

$$\sin v : \sin v' = \cos u : \cos u',$$

so that

$$\frac{\sin u}{\cos u} \cdot \circ y = \frac{\sin u'}{\cos u'} \cdot \circ y',$$

or

$$\circ y : \circ y' = \operatorname{tg} u' : \operatorname{tg} u = \operatorname{tg} w : \operatorname{tg} w'.$$

Moreover, let

$$\overrightarrow{CN} \parallel \overrightarrow{AB}, \quad \overrightarrow{C'N'} \parallel \overrightarrow{AB},$$

and let $\overset{\frown}{CD}$, $\overset{\frown}{C'D'}$ be L-lines perpendicularly intersecting the straight line $AB$. By **§21** we also have

$$\circ y : \circ y' = r : r';$$

thus

$$r:r' = \text{tg } w:\text{tg } w'.$$

Now let the distance $p$ measured from $A$ increase indefinitely. Then

$$w \to z \quad \text{and} \quad w' \to z'.$$

Hence also

$$r:r' = \text{tg } z:\text{tg } z'.$$

The constant $r:\text{tg } z$ which is independent of $r$ will be denoted by $k$. If $y \to 0$, then

$$\frac{r}{y} = \frac{k \cdot \text{tg } z}{y} \to 1$$

and therefore

$$\frac{y}{\text{tg } z} \to k.$$

§29 yields

$$\text{tg } z = \frac{1}{2}(Y - Y^{-1}).$$

Thus

$$\frac{2y}{Y - Y^{-1}} \to k$$

or, by §24,

$$\frac{2yK^{\frac{Y}{K}}}{K^{\frac{2y}{k}} - 1} \to k.$$

But it is well known that, for $y \to 0$ the limit of this expression is

$$\frac{k}{\log K}.$$

Consequently,

$$\frac{k}{\log K} = k$$

and

$$K = e = 2.7182818 \ldots.$$

Here, too, this number seems to have outstanding significance. Thus if, from now onwards, $k$ denotes the distance whose corresponding $K$ is just equal to $e$, then

$$r = k \cdot \text{tg } z.$$

On the other hand, in §21 we pointed out that

$$\circ y = 2\pi r.$$

So by §24

$$\circ y = 2\pi k \cdot \text{tg } z = \pi k(Y - Y^{-1}) =$$

$$= \pi k\left(e^{\frac{y}{k}} - e^{-\frac{y}{k}}\right) = \frac{\pi y}{\log Y}(Y - Y^{-1}).$$

***Note on 30***

Section 30 provides the rectification of the horocycle. That is to say, Bolyai obtains the formula for the length of a horocyclic arc *r* in terms of its corresponding semichord *y*. His theorem states that

$$\boldsymbol{r = k\ \tfrac{1}{2}\,(Y - Y^{-1})}$$

$$\text{i.e. } r = k \sinh y/k,$$

where *k* is the linear constant which here makes its first explicit appearance in the *Appendix* (though in the Latin original the constant, rather confusingly, is denoted by *i*).

In Figure 20 AL, CM and C'M' are perpendicular to AC. B is any point on AL, where AB = *p*. CN is drawn parallel to AL, where ∠MCN = *z* is the complement of the angle of parallelism corresponding to the segment AC. C'N' is also drawn parallel to AL, where ∠M'C'N' = *z'* is the complement of the angle of parallelism corresponding to the segment AC'. A horocycle of length *r* is drawn through C and D, while a horocycle of length *r'* is drawn through C' and D'. The semichord AC, corresponding to the horocyclic arc CD, is of length *y*. Similarly the semichord AC', corresponding to the horocyclic arc C'D', is of length *y'.* (It is implied, but not shown in the figure, that chords of length 2*y* and 2y', respectively, would join the endpoints of horocycles of length 2*r* and 2*r'*).

The proof begins with the neutral sine theorem applied to triangles BCA and BC'A, which implies that

$$\sin u : \sin v = {}^{\circ}p : {}^{\circ}y$$

and

$$\sin u' : \sin v' = {}^{\circ}p : {}^{\circ}y'.$$

Therefore

$$(\sin u / \sin v)\ .\ {}^{\circ}y = (\sin u' / \sin v')\ .\ {}^{\circ}y'\ .$$

But, as shown Section 27,

$$\sin w : \sin v = \sin w' : \sin v' = \text{constant},$$

which implies, given that *w* is the complement of *u* and *w'* is the complement of *u'*, that

$$\cos u : \sin v = \cos u' : \sin v'$$

$$\text{i.e. } \sin v : \sin v' = \cos u : \cos u'.$$

It follows that

$$(\sin u / \cos u) \, . \, {}^{\circ}y = (\sin u' / \cos u') \, . \, {}^{\circ}y'$$

or

$$\tan u \, . \, {}^{\circ}y = \tan u' \, . \, {}^{\circ}y',$$

so

$${}^{\circ}y : {}^{\circ}y' = \tan u' : \tan u$$

$$= \tan w : \tan w' \quad \text{(since } u \text{ and } u' \text{ are the complements of } w \text{ and w').}$$

But, as shown in Section 25,

$${}^{\circ}y = {}^{\circ}r.$$

Therefore

$$r : r' = \tan w : \tan w' \, .$$

Now as B recedes from A along AL, $p \to \infty$ and $w \to z$ while $w' \to z'$, so

$$r : r' = \tan z : \tan z'$$

$$\text{i.e. } r / \tan z = r' / \tan z'.$$

Since $r$ is independent of $z$, the quotient $r / \tan z$ is a constant, denoted by $k$. Thus

$$r = k \tan z.$$

Now given that $\angle$MCA is a right angle, $\tan z = \cot \angle$NCA. But $\angle$NCA is the angle of parallelism corresponding to AC = $y$. Using the alternative form of the fundamental identity of parallelism given in Section 29, it follows that

$$\tan z = ½ \, (Y - Y^{-1})$$

$$= \sinh y/k$$

and so

$$r = k \sinh y/k,$$

which rectifies the horocycle as required.

Bolyai also investigates the nature of the linear constant $k$. Starting from the identity

$$r = k \tan z,$$

he divides both sides by $y$ to obtain

$$r/y = k \tan z / y.$$

But it is clear from Figure 20 that as $y \to 0$, $r/y \to 1$, so

$$k \tan z / y \to 1,$$

whence

$$y / \tan z \to k.$$

But given that $\tan z = ½ (Y - Y^{-1})$, it follows that

$$2y / Y - Y^{-1}) \to k.$$

Now, by Section 24, $Y$ can be expressed as $K^{y/k}$. So multiplying the numerator and denominator of the last equation by $K^{y/k}$, he obtains

$$2y K^{y/k} / (K^{2y/k} - 1) \to k.$$

However the limit of the left hand side as $y \to 0$ is

$$k / \log K,$$

whence

$$k / \log K = k,$$

implying that

$$K = e.$$

But $K$ represents the ratio of two concentric horocyclic arcs separated by a radial distance $k$. Thus **the linear constant *k* can be defined as the radial distance between two horocyclic arcs whose lengths lie in the ratio of *e* : 1** (see the Note to Section 24).

Bolyai also uses his rectification of a horocyclic arc to obtain the formula for the circumference of a circle in the hyperbolic plane. The circumference of a circle swept out on the horosphere by the rotation of a horocyclic arc $r$ is given by the Euclidean formula, $2\pi r$. But Section 25 shows that this circle is the same as the circle generated in the hyperbolic plane by the rotation of the corresponding semichord $y$. In view of the fact that $r = k \sinh y/k$, the circumference of the circle is given by

$$\mathbf{2\pi r = 2\pi k \sinh y/k.}$$

## §31

In system **S**, the knowledge of three equations is sufficient for solving every *rectilinear right triangle*, and this renders easy to solve any triangle. Namely, we only need to know the equations that specify the relation

(I) between $a$, $c$, $\alpha$,

(II) between $a$, $\alpha$, $\beta$,

(III) between $a$, $b$, $c$,

where $a$, $b$ denote the legs, $c$ the hypotenuse, and $\alpha$, $\beta$ the angles opposite to the legs. From them, of course, the three remaining equations can be deduced by elimination.

(I) According to §§ **25** and **30**

$$1:\sin\alpha = (C-C^{-1}):(A-A^{-1}) = \left(e^{\frac{c}{k}}-e^{-\frac{c}{k}}\right):\left(e^{\frac{a}{k}}-e^{-\frac{a}{k}}\right).$$

This is the equation for $a$, $c$, $\alpha$.

(II) If $\overrightarrow{BM}\|\overrightarrow{CN}$, then §**27** yields

$$\cos\alpha:\sin\beta = 1:\sin u,$$

while §**29** yields

$$1:\sin u = \frac{1}{2}(A+A^{-1}).$$

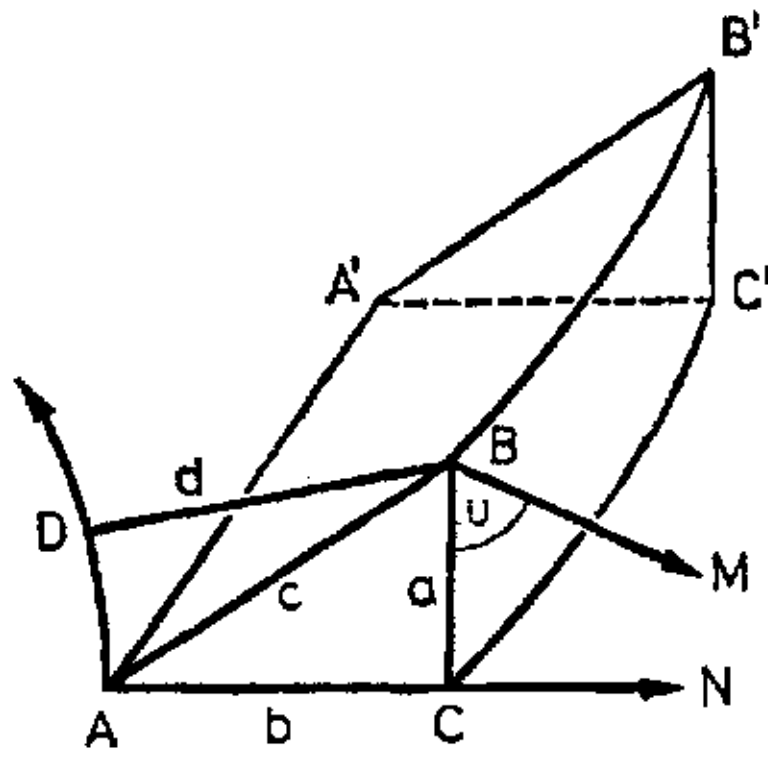

**Figure 21**

Consequently,

$$\cos\alpha:\sin\beta = \frac{1}{2}(A+A^{-1}) = \frac{1}{2}\left(e^{\frac{a}{k}}+e^{-\frac{a}{k}}\right).$$

This is the equation for $a$, $\alpha$, $\beta$.

(III) If $AA'$ is perpendicular to the plane $BAC$, if (with the notation introduced in §**27**) $\widehat{BB'}\|\overline{AA'}$ and $\widehat{CC'}\|\overline{AA'}$, finally if the plane $B'A'C'$ is perpendicular to the line $AA'$, then in the same way as in §**27** it follows that

$$\frac{\widehat{BB'}}{\widehat{CC'}} = \frac{1}{\sin u} = \frac{1}{2}(A+A^{-1}),$$

$$\frac{\widehat{CC'}}{\overline{AA'}} = \frac{1}{2}(B+B^{-1}),$$

and

$$\frac{\widehat{BB'}}{\widehat{AA'}} = \frac{1}{2}(C+C^{-1}).$$

Hence

$$\frac{1}{2}(C+C^{-1}) = \frac{1}{2}(A+A^{-1})\cdot\frac{1}{2}(B+B^{-1})$$

and, consequently,

$$(e^{\frac{c}{k}}+e^{-\frac{c}{k}}) = \frac{1}{2}(e^{\frac{a}{k}}+e^{-\frac{a}{k}})(e^{\frac{b}{k}}+e^{-\frac{b}{k}}).$$

This is the equation for $a$, $b$, $c$.

If

$$\sphericalangle CAD = R$$

and

$$BD \perp AD,$$

then

$$\circ c:\circ a = 1:\sin\alpha$$

and

$$\circ c:\circ d = 1:\cos\alpha,$$

where $d=\widehat{BD}$. Therefore, denoting $\circ x\cdot\circ x$ by $\circ x^2$, obviously

$$\circ a^2+\circ d^2 = \circ c^2.$$

But, in view of §27 and (II),

$$\circ d = \circ b\cdot\frac{1}{2}(A+A^{-1}),$$

which yields

$$(e^{\frac{c}{k}}-e^{-\frac{c}{k}})^2 = \frac{1}{4}(e^{\frac{a}{k}}+e^{-\frac{a}{k}})^2(e^{\frac{b}{k}}-e^{-\frac{b}{k}})^2+(e^{\frac{a}{k}}-e^{-\frac{a}{k}})^2.$$

This is another equation for $a$, $b$, $c$; its right-hand side can be easily brought to symmetric form.

Finally, from

$$\frac{\cos\alpha}{\sin\beta} = \frac{1}{2}(A+A^{-1})$$

and

$$\frac{\cos\beta}{\sin\alpha} = \frac{1}{2}(B+B^{-1})$$

by the aid of (III) we obtain

$$\operatorname{ctg}\alpha\cdot\operatorname{ctg}\beta = \frac{1}{2}(e^{\frac{c}{k}}+e^{-\frac{c}{k}}).$$

This is the equation for $\alpha$, $\beta$, $c$.

***Note on 31***

In this section Bolyai asserts that three fundamental equations suffice to solve any right angled triangle and therefore any triangle in hyperbolic geometry. In respect of a right angled triangle ABC with sides *a*, *b* and hypotenuse *c* and opposite angles $\alpha$, $\beta$ and *R* these identities are:

(1) $1 : \sin\alpha = (C - C^{-1}) : (A - A^{-1})$

$$= (e^{c/k} - e^{-c/k}) : (e^{a/k} - e^{-a/k}).$$

This result follows from the neutral sine rule, which implies that in triangle ABC

$$\sin R : \sin\alpha = {}^{\circ}c : {}^{\circ}a$$

and so, by the result derived in Section 30 relating the circumference of a circle to its radius,

$$\mathbf{1 : \sin\alpha = \sinh c/k : \sinh a/k.}$$

(2) If BM || CN in Figure 21 shown below, then

$$\boldsymbol{\cos\alpha : \sin\beta = \cosh a/k}:$$

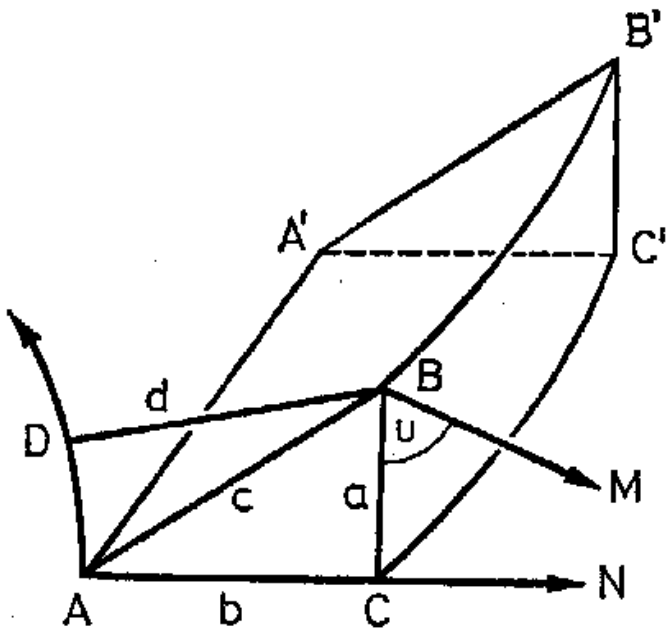

The argument presented in the course of the rectification of the equidistant (Section 27) implies that

$$\sin\angle\text{DAB} : \angle\sin\text{ABC} = 1 : \sin u.$$

But $\sin\angle\text{DAB} = \cos\angle\text{BAC}$ and by Section 27 again, $\sin u = 1/\cosh \text{BC}/k$, therefore

$$\cos\alpha : \sin\beta = \cosh a/k.$$

(3) $\boldsymbol{\cosh c/k = \cosh a/k \,.\, \cosh b/k}$

In the figure shown above, planes ABC and A'B'C' are perpendicular to AA', while BB' and CC' are both equidistant curves to AA'.

The result given in Section 27 implies that

arc BB'/arc CC' = 1/*sin u* = *cosh a/k*

arc CC'/AA' = *cosh b/k*

arc BB'/AA' = *cosh c/k*

It follows that

arc BB'/AA' = (arc BB'/CC') . (arc CC'/AA'),

whence

$$\cosh c/k = \cosh a/k \, . \, \cosh b/k.$$

This is the hyperbolic equivalent of the (Euclidean) Pythagorean theorem

$$c^2 = a^2 + b^2.$$

It can be seen by expanding *cosh x/k* in a power series

$$\cosh x/k = 1 + (x/k)^2/2! + (x/k)^4/4! + (x/k)^6/6! + \ldots ,$$

that the hyperbolic formula transposes into the Pythagorean theorem if the sides *a, b, c* are taken as infinitely small in comparison to the linear constant *k*, demonstrating once again that the theorems of Euclidean geometry hold true over infinitesimal domains in the hyperbolic plane.

Bolyai derives two further equations for a right angled triangle ABC. The first states that

$$(e^{c/k} - e^{-c/k})^2 = \tfrac{1}{4}\,(e^{a/k} - e^{-a/k})^2\,(e^{b/k} + e^{-b/k})^2 + (e^{a/k} - e^{-a/k})^2.$$

This is justified as follows. If, in Figure 21 shown above, $\angle$BCA = $\angle$CAD = $\angle$ADB = *R*, then ADBC is a quadrilateral with three right angles and an acute angle at C, a figure known as a Lambert quadrilateral. The neutral sine rule applied to triangles ABD, ABC implies that

$$^{\circ}c : {}^{\circ}a = 1 : \sin\alpha$$

and

$$^{\circ}c : {}^{\circ}d = 1 : \cos\alpha,$$

whence, by cross multiplication and using the identity $\sin^2 x + \cos^2 x = 1$,

$$^{\circ}a^2 + {}^{\circ}d^2 = {}^{\circ}c^2,$$

where $^{\circ}c^2$ denotes $^{\circ}c \, . {}^{\circ}c$

But if BD = *d*, the result obtained in Section 27 implies that

$$^{\circ}d = {}^{\circ}b \cosh a/k,$$

and by (1) above

$$^{\circ}c : {}^{\circ}a = \sinh c/k : \sinh a/k.$$

Therefore

$$\sinh^2 c/k = \cosh^2 a/k \sinh^2 b/k + \sinh^2 a/k.$$

The second formula states that

$$\cot \alpha \cot \beta = ½ (e^{c/k} + e^{-c/k}).$$

This identity is based on the fact that by (2) above

$$\cos \alpha : \sin \beta = \cosh a/k,$$

while by the same token

$$\cos \beta : \sin \alpha = \cosh b/k.$$

Hence, using (3) above,

$$\cot \alpha \cot \beta = \cosh c/k.$$

## §32

It remains to indicate the way of solving problems in system S and finally, this being accomplished for some examples most often encountered, to describe frankly what our theory offers.

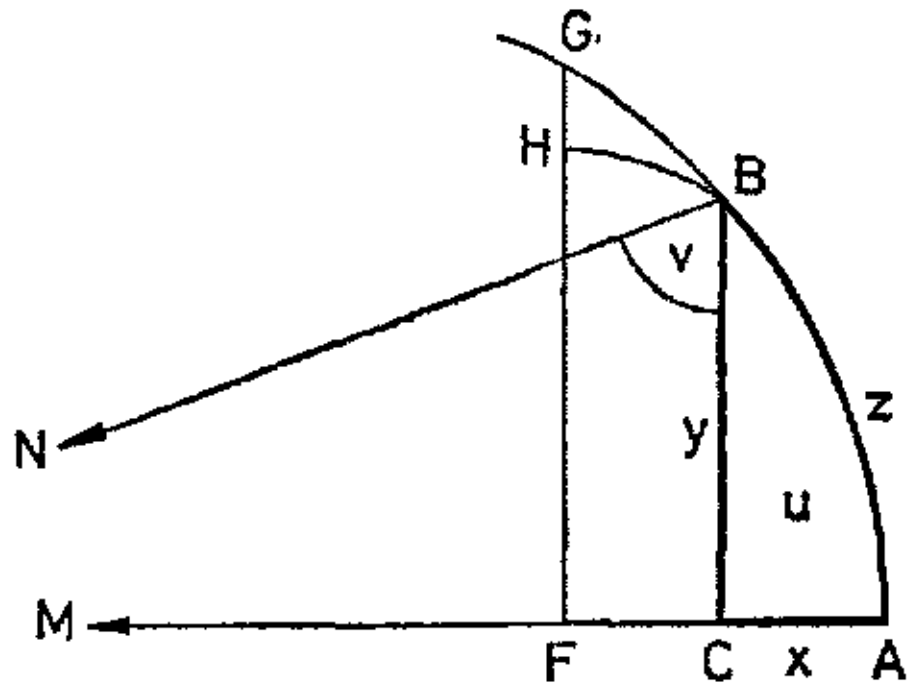

Figure 22

(I) Let $\overset{\frown}{AB}$ be a plane curve whose equation in rectangular coordinates reads $y=f(x)$. Denote by $dz$ any increment of $z$, and by $dx$, $dy$, $du$ the increments, corresponding to $dz$, of $x$, $y$ and the area $u$, respectively. Further, with the notation introduced in §27, let $\overset{\frown}{BH} \| \overline{CF}$. Making use of §31, express $\overset{\frown}{BH}/dx$ as a function of $y$. Find the limit of $dy/dx$ as $dx$ tends to the limit 0 (which will be tacitly assumed whenever taking limits of this kind). From these, the limit of $dy/\overset{\frown}{BH}$ and, consequently, the value tg $HBG$ become known. Since $\sphericalangle HBC$ is neither greater nor smaller than the right angle, it is a right angle. As a result, we obtain the tangent of the curve $\overset{\frown}{BG}$ at the point $B$ in terms of $y$.

***Note on 32, I***

The fact that **the theorems of Euclidean geometry are valid over an infinitesimal domain in the hyperbolic plane** enables Bolyai to apply the calculus to problems of area and volume involving circles, equidistants and horocycles.

Figure 22 shows an arbitrary curve ABG whose equation is given by $y = f(x)$. The origin is at A, AC = $x$, BC = $y$ and the length of the curve AB is denoted by $z$. CF = $dx$, HG = $dy$ and BG = $dz$ represent small increments of $x$, $y$ and $z$ respectively. Arc BH is the infinitesimal curve equidistant to CF.

(II) It can be proved that

$$\frac{dz^2}{dy^2+\widehat{BH}^2} \to 1.$$

Hence we can calculate the limit of $dz/dx$ and find $z$ as a function of $x$ by integration.

In system S, it is possible to deduce the equation of any curve given concretely; for instance, that of the L-line.

If $\overrightarrow{AM}$ is an axis of L, then any half-line $\overrightarrow{CB}$ starting from the point $C$ of $\overrightarrow{AM}$ intersects L, since by §19 every straight line through $A$, excepting $AM$, intersects L. If, however, $\overrightarrow{BN}$ is also an axis, then by §28

$$X = 1 : \sin CBN$$

and by §29

$$Y = \operatorname{ctg}\frac{1}{2}CBN.$$

Hence

$$Y = X+\sqrt{X^2-1}$$

or

$$e^{\frac{y}{k}} = e^{\frac{x}{k}}+\sqrt{e^{\frac{2x}{k}}-1}.$$

This is the equation required. It yields

$$\frac{dy}{dx} \to X(X^2-1)^{-\frac{1}{2}}.$$

But

$$\frac{\widehat{BH}}{dx} \to 1 : \sin CBN = X.$$

Thus

$$\frac{dy}{\widehat{BH}} \to (X^2-1)^{-\frac{1}{2}},$$

$$1+\left(\frac{dy}{\widehat{BH}}\right)^2 \to X^2(X^2-1)^{-1},$$

$$\left(\frac{dz}{\widehat{BH}}\right)^2 \to X^2(X^2-1)^{-1},$$

$$\frac{dz}{\widehat{BH}} \to X(X^2-1)^{-\frac{1}{2}}$$

and

$$\frac{dz}{dx} \to X^2(X^2-1)^{-\frac{1}{2}}.$$

By integration we obtain

$$z = k(X^2-1)^{\frac{1}{2}} = k \operatorname{ctg} CBN$$

in accordance with §30.

***Note on 32, II***

Bolyai begins by asserting that as $x \to 0$,

$$dz^2 / (dy^2 + \text{arc BH}^2) \to 1.$$

This can be justified by treating BHG as an infinitesimal right angled triangle with hypotenuse *dz*, height *dy* and base BH, to which the Euclidean Pythagorean formula applies. Therefore

$$dz^2 = dy^2 + \text{arc BH}^2.$$

However, as shown in the Note to Section 27, the equidistant arc BH can be expressed as CF *cosh* BC/*k* = *dx cosh y*/*k*. Thus the element of arc in rectangular coordinates can be expressed in modern notation as

$$dz^2 = dy^2 + \cosh^2 y/k\, dx^2.$$

Choosing to treat the curve $y = f(x)$ as a horocycle, Bolyai proceeds to obtain the length of arc AB = *z* in terms of ∠CBN = *v*. His theorem states that

$$\mathbf{z = k \cot \angle CBN}.$$

To prove this he adapts the configuration shown in Figure 18 in Section 28, where AM || BN, BC = *y*, ∠CBN = *v* and ∠ACB = *u*, a right angle:

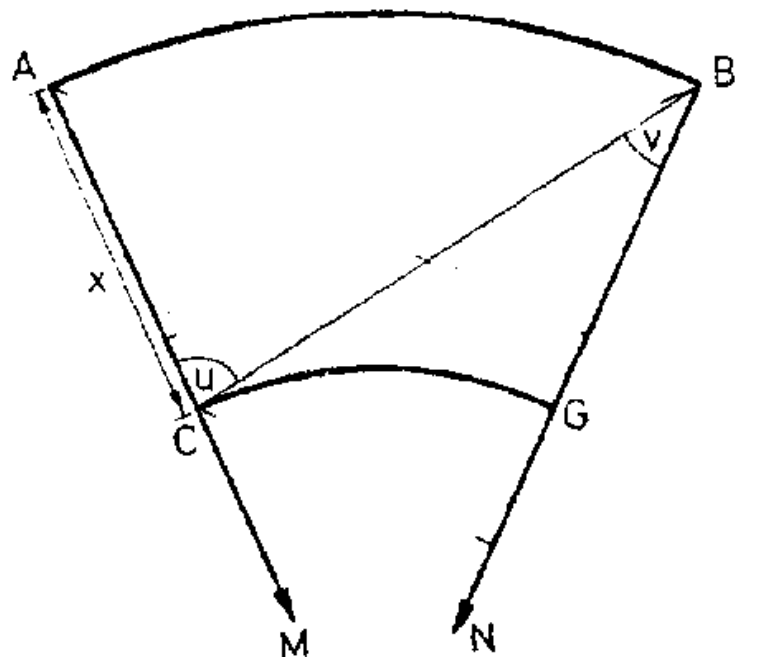

The lemma of Section 28 implies that

$$X = 1/\sin v.$$

Furthermore, given that BC = *y* is the segment corresponding to the angle of parallelism ∠CBN, the fundamental identity obtained in Section 29 gives

$$Y = \cot \tfrac{1}{2} v.$$

Now

$$\cot \tfrac{1}{2} v = (1 + \cos v) / \sin v$$

$$= 1/\sin v + \cos v/\sin v$$

$$= 1/\sin v + \sqrt{(1 - \sin^2 v)}/\sin v$$

$$= 1/\sin v + \sqrt{(1/\sin^2 v - 1)}.$$

Therefore

$$Y = X + \sqrt{(X^2 - 1)}$$

$$\text{i.e.} \quad e^{y/k} = e^{x/k} + (e^{2x/k} - 1)^{1/2},$$

whence

$$dy/dx \rightarrow X\,(X^2 - 1)^{-½}.$$

But by Section 27

$$\text{arc BH} : \text{CF} = \text{arc BH} : dx$$

$$= 1 : \sin v$$

$$= X.$$

Thus

$$dy/\text{arc BH} = dy/dx\ .\ dx/\text{arc BH}$$

$$= X\,(X^2 - 1)^{-½}.\ 1/X$$

$$= (X^2 - 1)^{-½},$$

whence

$$(dy/\text{arc BH})^2 = (X^2 - 1)^{-1}$$

and

$$1 + (dy/\text{arc BH})^2 = 1 + 1/(X^2 - 1)$$

$$= X^2/(X^2 - 1).$$

Now by the earlier result on the element of arc,

$$dz^2 = dy^2 + \text{arc BH}^2,$$

so

$$dz^2/\text{arc BH}^2 = dy^2/\text{arc BH}^2 + 1$$

$$= X^2/(X^2 - 1),$$

whence

$$dz/\text{arc BH} = X/(X^2 - 1)^{½}.$$

Now

$$dz/dx = dz/\text{arc BH}\ .\ \text{arc BH}/dx$$

$$= X/(X^2 - 1)^{½}.\ X$$

$$= X^2/(X^2 - 1)^{½}.$$

Integration then gives

$$z = k\,(X^2 - 1)^{1/2}$$

$$= k\,(1/\sin^2 v - 1)^{1/2}$$

$$= k \cot \angle\mathrm{CBN},$$

as required.

(III) For the area of the domain *HFCBH* we obviously have

$$\frac{du}{dx} \to \frac{HFCBH}{dx}.$$

This depends only on $y$ and must first be expressed in terms of $y$; then $u$ is obtained by integration.

Denoting (see Fig. 17) $\overline{AB}$ by $p$, $\overline{AC}$ by $q$, $\widehat{CD}$ by $r$ and the area of the domain *CABDC* by $s$, it can be shown as in (II) that

$$\frac{ds}{dq} \to r = \frac{1}{2}p(e^{\frac{q}{k}} + e^{-\frac{q}{k}}),$$

and hence by integration

$$s = \frac{1}{2}pk(e^{\frac{q}{k}} - e^{-\frac{q}{k}}).$$

This can be derived also without integration.

If we deduce the equation of the circle using §31, (III), or that of the straight line using §31, (II), or that of a conic section using the arguments above, then the areas enclosed by these lines can also be calculated.

Clearly, the area of the surface $t$ parallel to the plane figure $p$ at distance* $q$ is to $p$ as are the second powers of the homologous line segments or, more explicitly, as

$$\frac{1}{4}(e^{\frac{q}{k}} + e^{-\frac{q}{k}})^2 : 1.$$

It is easy to see that calculating the volume in a similar way requires two integrations (here, in fact, even the differential can only be obtained by integration), and that first of all the volume of the solid enclosed by $p$, $t$ and all lines perpendicular to the plane of $p$ and connecting the boundaries of $p$ and $t$ must be determined. We find (whether by means of integration or without it) that this volume equals

$$\frac{1}{8}pk(e^{\frac{2q}{k}} - e^{-\frac{2q}{k}}) + \frac{1}{2}pq.$$

In S, also the surface area of a solid can be obtained, as well as the curvature, evolute, evolvent of an arbitrary curve, etc. As to the curvature in system S, it may either be that of an L-line, or characterised by the radius of a circle or by the distance between a curve which is parallel to a straight line and that line. Really, in view of the foregoing it is easy to show that there are no uniform plane curves other than L-lines, circles, and curves parallel to a straight line.

***Note on 32, III***

In this sub-section Bolyai introduces the calculation of areas and volumes. Referring to Figure 22, he observes that the rate of change of area $du/dx$ beneath a curve $y = f(x)$ is equal to $y$, which is of course a statement of the first part of the fundamental theorem of calculus (FTC).

He then applies this to the calculation of area beneath the equidistant curve shown below, giving the area $s$ of CABD as

$$\boldsymbol{s = \tfrac{1}{2}\, pk\, (e^{q/k} + e^{-q/k})}.$$

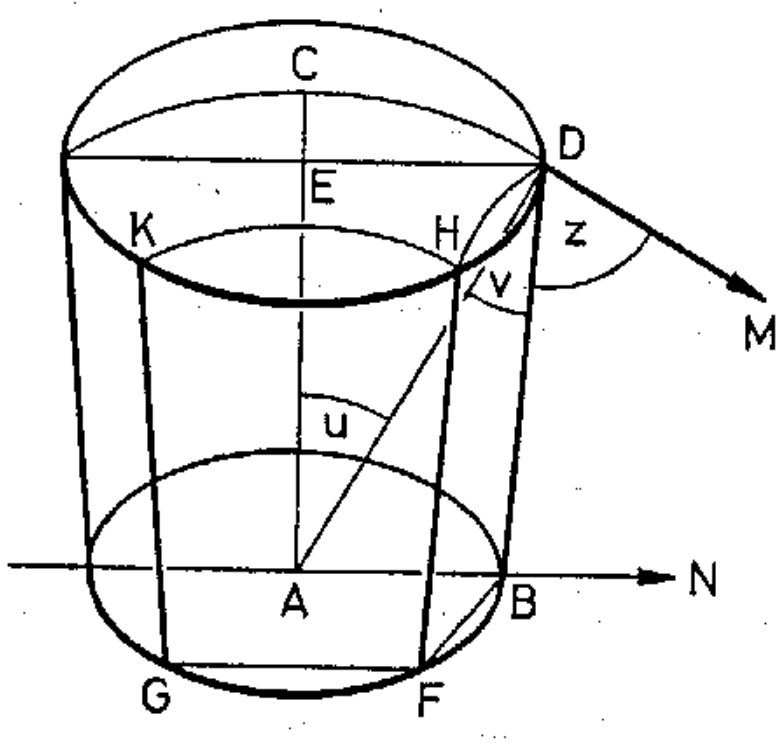

If AB = $p$, AC = $q$, arc CD = $r$ and area CABD = $s$ with the origin at A, then by the FTC

$$ds/dq \to r,$$

where, by the result obtained in Section 27, $r = p \cosh q/k$. Therefore the area $s$ is given by

$$\boldsymbol{s} = p \int \cosh q/k \, dq$$

$$= pk \sinh q/k.$$

The equidistant surface $t$ compared to its projection in the plane $p$ at a distance $q$ is in the ratio of the square of the line segments in the former compared to the latter i.e. $\cosh^2 q/k : 1$. On this basis Bolyai gives the volume of the solid beneath an equidistant surface $t$ with base $p$ whose sides of length $q$ are perpendicular to the base as

$$\boldsymbol{1/8\, pk\, (e^{2q/k} - e^{-2q/k}) + \tfrac{1}{2}\, pq}.$$

This can be justified as follows:

$$\text{Volume} = p \int_0^q \cosh^2 z/k \; dz$$

$$= p \int_0^q \frac{\cosh\left(\frac{2z}{k} + 1\right)}{2} dz$$

$$= \tfrac{1}{4}\; pk \sinh 2q/k + \tfrac{1}{2}\, pq.$$

(IV) As in (III), for the area and circumference of the circle we have

$$\frac{d(\odot x)}{dx} \to \circ x.$$

Owing to **§30**, integration yields

$$\odot x = \pi k^2 (e^{\frac{x}{k}} - 2 + e^{-\frac{x}{k}}).$$

***Note on 32, IV***

In this sub-section Bolyai gives the area $A$ of a circle of radius $x$ as

$$\boldsymbol{A = \pi k^2 (e^{x/k} - 2 + e^{-x/k}).}$$

He observes that the derivative of the area of a circle is equal to the circumference, implying that the area is equal to the integral of the circumference. Now the circumference of a circle of radius $x$ is given at the end of Section 30 as $2\pi k \sinh x/k$.

Therefore the area $A$ of a circle of radius $x$ is

$$A = 2\pi k \int_0^x \sinh \frac{z}{k}\, dz$$

$$= 2\pi k^2 (\cosh x/k - 1)$$

$$= \pi k^2 (e^{x/k} - 2 + e^{-x/k}).$$

(V) For the area* $CABDC = u$ enclosed by the **L**-line $\overset{\frown}{AB} = r$, by the **L**-line $\overset{\frown}{CD} = y$ parallel to it, and by the segments $\overline{AC} = \overline{BD} = x$ we have

$$\frac{du}{dx} \to y.$$

From **§24**

$$y = re^{-\frac{x}{k}},$$

and hence, by integration,

$$u = rk(1 - e^{-\frac{x}{k}}).$$

If $x \to \infty$, then $e^{-\frac{x}{k}} \to 0$, so that in **S**

$$u \to rk.$$

In what follows, the size of the domain $MABN$ will always mean this limit.

It can be established in a similar way that if $p$ is a figure in **F** then the volume enclosed by $p$ and all axes starting from the boundary of $p$ equals $\frac{1}{2}pk$.

***Note on 32, V***

Bolyai proceeds to obtain the area bounded by two horocyclic arcs:

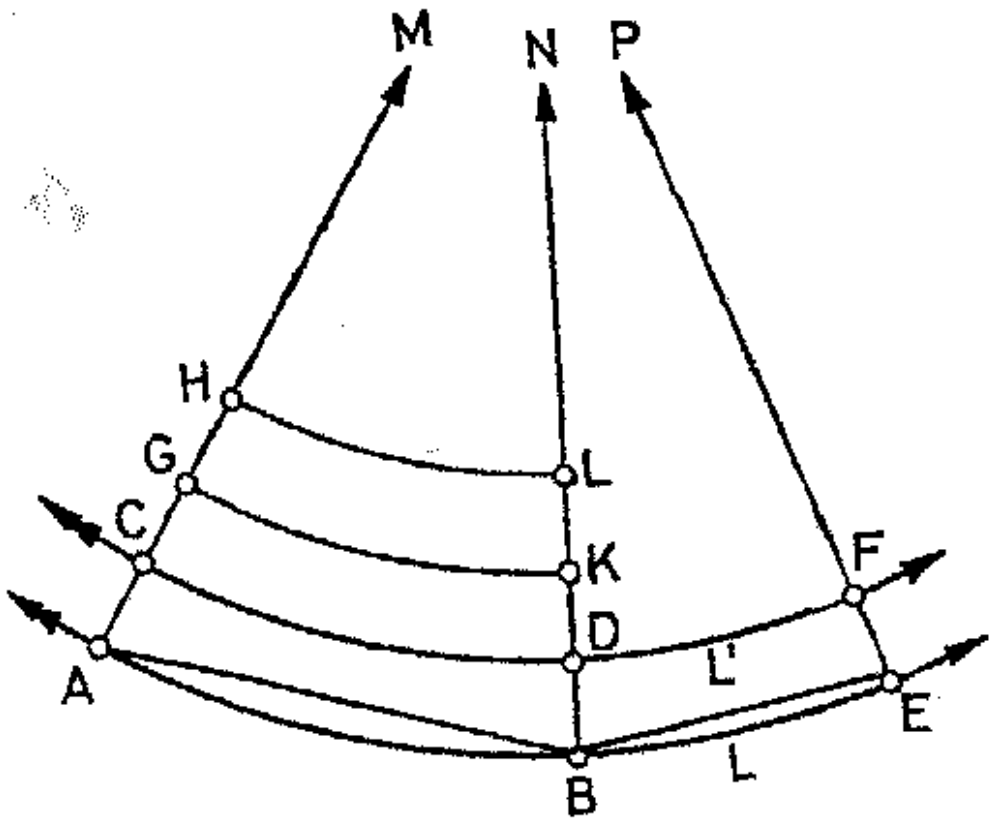

If arc AB = $r$, arc CD = y while AC = BD = x, then the (FTC) implies in respect of the area $u$ of CABD that

$$du/dx = y,$$

where by Section 24

$$y = r\,e^{-x/k}.$$

Then

$$u = r\int_0^x e^{-z/k}\,dz$$

$$= rk\,(1 - e^{-x/k}).$$

Bolyai observes that if $x \to \infty$, then $e^{-x/k} \to 0$ and so the area bounded by arc AB and the parallels AM, BN is

$$\boldsymbol{u \to rk}.$$

He then asserts that the volume $V$ enclosed by a figure of area $p$ on the horosphere and its parallel axes extending to an infinitely distant point is given by

$$\boldsymbol{V = \tfrac{1}{2}\,pk}:$$

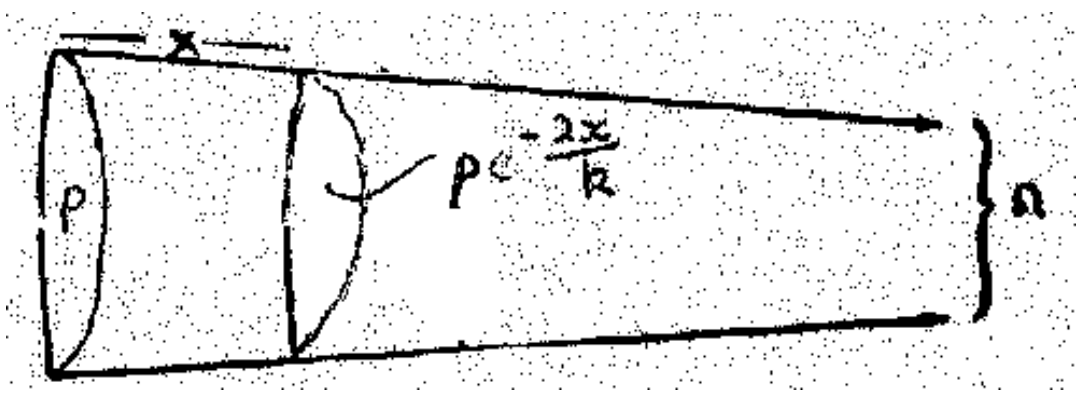

This result can be justified as follows. If the figure on the left hand horosphere has area *p,* the area of the corresponding figure on the right hand horosphere at a distance *x* will be in the ratio of the square of the corresponding line elements, so its area is given by $pe^{-2x/k}$. By the FTC the rate of change of volume with respect to *x* is

$$dV/dx = pe^{-2x/k},$$

so the volume *V* of the solid between the concentric horospheres is given by

$$V = p \int_0^x e^{-2z/k}\, dz$$

$$= \tfrac{1}{2}\, pk\,(1 - e^{-2x/k}).$$

As $x \to \infty$, $V \to \tfrac{1}{2}pk$, as required.

(VI)** If the central angle of the spherical cap $z$ is $2u$, the circumference of the great circle is $p$, and the circular arc corresponding to central angle $u$ is $\widehat{FC} = x$, then by **§25**

$$1 : \sin u = p : \circ\overline{BC}.$$

Thus

$$\circ\overline{BC} = p \sin u.$$

On the other hand,

$$x = \frac{pu}{2\pi} \quad \text{and} \quad dx = \frac{p\,du}{2\pi}.$$

Further

$$\frac{dz}{dx} \to \circ\overline{BC},$$

so that

$$\frac{dz}{du} \to \frac{p^2}{2\pi} \sin u,$$

and by integration

$$z = \frac{1 - \cos u}{2\pi} p^2.$$

Consider the **F**-surface that contains the circle $p$ passing through the centre $F$ of the spherical cap. The line*** $t$ intersects **F** perpendicularly at the point $E$. The planes $FEM$ and $CEM$ containing the radii $\overline{AF}$ and $\overline{AC}$ are perpendicular to the surface **F** and intersect it in the curves *FEG* and *CE*. Furthermore, consider the L-arc $\overset{\frown}{CD}$ starting from *C* and perpendicular to the curve *FEG*, and the L-arc $\overset{\frown}{CF}$.

According to **§20** we have

$$\sphericalangle CEF = u$$

and by **§21**

$$\frac{\overset{\frown}{FD}}{p} = \frac{1-\cos u}{2\pi}.$$

Therefore

$$z = \overset{\frown}{FD} \cdot p.$$

But from **§21** it follows that

$$p = \pi \cdot \overset{\frown}{FG}.$$

Thus

$$z = \pi \cdot \overset{\frown}{FD} \cdot \overset{\frown}{FG}.$$

On the other hand, again by **§21**,

$$\overset{\frown}{FD} \cdot \overset{\frown}{FG} = \overset{\frown}{FC}^2.$$

Consequently, for the circular domain of the **F**-surface we have

$$z = \pi \cdot \overset{\frown}{FC}^2 = \odot \overset{\frown}{FC}.$$

Now let* $\overset{\frown}{BJ} = \overset{\frown}{CJ} = r$. By **§30**

$$2r = k(Y - Y^{-1}).$$

Consequently, according to **§21**, for the circular domain of **F**

$$\odot 2r = \pi k^2 (Y - Y^{-1})^2.$$

On the other hand, (IV) yields

$$\odot 2y = \pi k^2 (Y^2 - 2 + Y^{-2}).$$

Thus for the circular domain of **F** we have

$$\odot 2r = \odot 2y.$$

Therefore *the area of the spherical cap z is equal to the area of the circle described with the chord* $\overline{FC}$ *as radius.*

Hence the surface area of the whole sphere is

$$\odot \overline{FG} = \overset{\frown}{FG} \cdot p = \frac{p^2}{\pi}.$$

Therefore *the surface areas of two spheres are to each other as are the second powers of the circumferences of their great circles.*

***Note on 32, VI***

Bolyai begins by obtaining the area of a spherical cap in terms of the angle it subtends at the centre of the sphere :

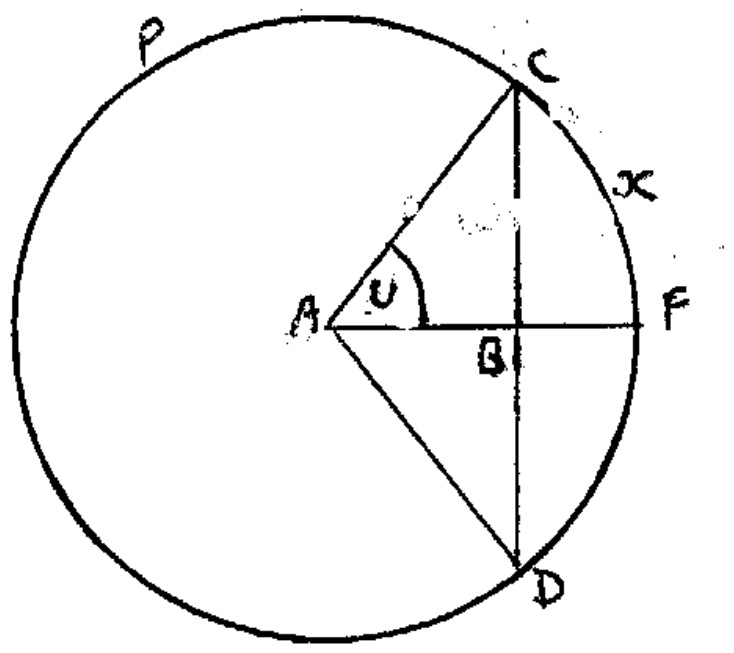

In the figure FCD represents a cap of area *z* subtended by $\angle$CAD = 2*u* at the centre of the sphere. Arc FC = *x*, AC = *r* and BC = *y*, while the circumference of the sphere = *p*. Bolyai's theorem states that the area of the cap is given by

$$z = p^2(1 - \cos u)/2\pi.$$

His proof is as follows. In triangle ACB the neutral sine rule implies that

$$1/\sin u = {}^{\circ}AC / {}^{\circ}BC$$

$$= p / {}^{\circ}BC,$$

hence

$${}^{\circ}BC = p \sin u.$$

But

$$x/p = u/2\pi,$$

so

$$x = pu/2\pi,$$

whence

$$dx/du = p/2\pi.$$

Now the surface area *z* of the cap FCD is given by $z = \int 2\pi\, y\, dx$ (where *dx* would normally be written *ds* nowadays). Therefore

$$dz/dx = 2\pi\, y$$

$$= {}^{\circ}BC.$$

Then

$$dz/du = dz/dx \,.\, dx/du$$

$$= {}^{\circ}\mathrm{BC}\; p\,/2\pi$$

$$= p^2 \sin u\,/2\pi$$

and so by integration

$$\boldsymbol{z} = (1/2\pi)\, p^2 \int_u^0 \sin\theta \, d\theta$$

$$= p^2(1 - \cos u)/2\pi,$$

as required.

Bolyai then uses this result to obtain the area $z$ of a horospherical cap swept out by a horocyclic arc as it rotates about the axis of the horosphere:

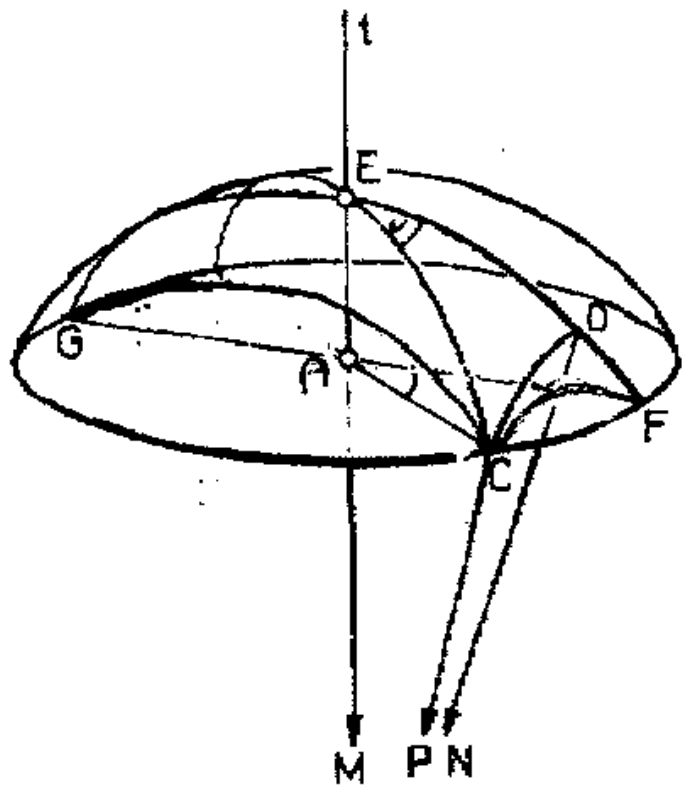

In the figure the plane of the circle CFG with circumference $p$ is perpendicular to the axis ME. By hypothesis plane PCDN is perpendicular to MEDN, so that ∠CDE = $R$. The dihedral angle between planes MECP and MEDN is equal to $u$, so ∠CEF = ∠CAF = $u$. Furthermore ∠FCG = $R$, being the angle in a horocyclic semicircle where arc FEG is the diameter. Arc EF = arc EC = $r$.

Bolyai's theorem states that **the area *z* of the horospherical cap swept out by the arc FC is equal to the area of the plane circle whose radius is the chord FC.** (The cap itself is not shown in the figure). His proof, with a few details filled in, runs as follows.

$$\text{arc FD} = \text{arc EF} - \text{arc ED}$$

$$= r - r\cos u$$

$$= r\,(1 - \cos u)$$

$$= 2\pi\, r\,(1 - \cos u)/2\pi$$

$$= p\,(1 - \cos u)/2\pi,$$

hence

$$\text{arc FD}/p = (1 - \cos u)/2\pi.$$

But the formula for a spherical cap given previously, being part of neutral geometry, can be applied to a horospherical cap, so

$$z = p^2(1 - \cos u)/2\pi$$

$$= \text{arc FD} \,.\, p$$

Given that the geometry of the horosphere is Euclidean, $p = \pi$ arc FG, so

$$z = \pi \text{ arc FD} \,.\, \text{arc FG}.$$

Now arc CD can be regarded as a mean proportional between arc FD and arc GD, so by Euclid's theorem (*Elements,* 6, 11) arc $CD^2$ = arc FD . arc GD. Furthermore the Pythagorean theorem implies that arc $FC^2$ = arc $CD^2$ + arc $FD^2$. Therefore

$$\text{arc FC}^2 = \text{arc FD} \,.\, \text{arc GD} + \text{arc FD}^2$$

$$= \text{arc FD (arc GD + arc FD)}$$

$$= \text{arc FD} \,.\, \text{arc FG}.$$

Thus the area of the horospherical cap with radius FC is given by the Euclidean formula

$$z = \pi \text{ arc FC}^2.$$

At this point Bolyai compares the area of a circle swept out on the horosphere by a horocycle of length $2r$ with the area of a circle in the hyperbolic plane of radius $2y$, where $2y$ is the length of the chord joining the endpoints of the horocycle.

By the result given in Section 30

$$2r = k\,(e^{y/k} - e^{-y/k}).$$

But Section 21 implies that the area $A$ of a circle on the horosphere is given by the Euclidean formula $\pi r^2$. So if the radius of this circle is $2r$, the area is

$$A = \pi\,(2r)^2$$

$$= \pi k^2\,(e^{2y/k} - 2 - e^{-2y/k}).$$

On the other hand Section 32, IV gives for the area of a circle of radius $2y$ in the hyperbolic plane

$$\pi\, k^2\,(e^{2y/k} - 2 - e^{-2y/k}),$$

which is the same. Bolyai concludes that **a cap on the horosphere of radius 2*r* is equal in area to the plane circle whose radius is equal to the chord of length 2*y* joining the endpoints of the horocycle***.* It follows that the area $z$ of the horospherical cap of radius arc FC is equal to the area of the plane circle of radius FC, where FC is the chord joining the endpoints of the arc FC, as required.

The surface area of the whole sphere is therefore given by the circle whose radius is FG. Hence

area° FG = area° arc FG

= arc FG . $p$

= $2r \,.\, 2\pi r$

= $4\pi r^2$.

(VII) We find in a similar way that in system S the volume of the sphere of radius $x$ is

$$\frac{1}{2}\pi k^3(X^2 - X^{-2}) - 2\pi k^2 x.$$

The area of the surface obtained by rotation of the arc* $\overset{\frown}{CD}$ about $\overline{AB}$ is

$$\frac{1}{2}\pi k p(Q^2 - Q^{-2}),$$

whereas the volume of the solid described by the figure $CABDC$ equals

$$\frac{1}{4}\pi k^2 p(Q - Q^{-1})^2.$$

*For the sake of brevity, we shall not explain how everything we have presented beginning from* (IV) *can be deduced also without integration.*

It can be proved that *if $k$ tends to infinity then the limit of any expression containing* $k$ (hence based on the hypothesis that $k$ exists) *coincides with the value of the same quantity valid in* $\Sigma$ (which system involves the non-existence of any $k$) *unless an identity is obtained.* Beware of the impression, however, as if *the system itself could be altered*; it is completely *determined in itself and by itself.* Only the *hypothesis* can be altered, as far as we are not led to a contradiction. Thus, *assuming* that in such expressions $k$ denotes the unique value whose $K$ is equal to $e$ if system S is true in reality, while the expression is thought to be replaced by the limit mentioned above if system $\Sigma$ is actually valid, then it is clear that all expressions obtained from the hypothesis of the reality of S are, in this sense, absolutely valid, though it remains perfectly unknown whether $\Sigma$ is, or is not, fulfilled.

Thus, for instance, the expression obtained in §30 yields (either by or without differentiation) the value

$$\circ x = 2\pi x,$$

well known is system $\Sigma$. From §31, (I) in the usual way we deduce

$$1 : \sin\alpha = c : a.$$

Moreover, from (II),

$$\frac{\cos\alpha}{\sin\beta} = 1,$$

that is

$$\alpha + \beta = R.$$

The first equation of (III) becomes an identity, hence is valid in $\Sigma$, though does not determine anything in it; the second, however, implies

$$c^2 = a^2 + b^2.$$

*These are the well-known fundamental equations of plane trigonometry in system* $\Sigma$.

* See Fig. 17.

Further, according to §32, in system $\Sigma$ the area as well as the volume appearing in (III) are equal to

$$pq,$$

(IV) yields

$$\odot x = \pi x^2,$$

from (VII) it follows that the volume of the sphere of radius $x$ is

$$\frac{4}{3}\pi x^3,$$

etc. Obviously, the theorems announced at the end of (VI) are also unconditionally true.

***Note on 32, VII***

Bolyai begins this subsection with three integrations. The first gives the volume $V$ of a sphere as

$$\mathbf{V = \tfrac{1}{2}\,\pi k^3\,(X^{2} - X^{-2}) - 2\pi k^2 x.}$$

This is based on the fact that the volume of a sphere is equal to the integral of its surface area. The surface area of a sphere with horocyclic radius $r$ is $4\pi r^2$. But $r = k \sinh y/k$, where $y$ is the semichord corresponding to a horocyclic arc $r$. Therefore the volume $V$ of a sphere of radius $x$ in hyperbolic geometry is given by

$$V = 4\pi k^2 \int_0^x \sinh^2 y/k \, dy$$

$$= 2\pi k^2 \int_0^x \cosh\,(2y/k - 1)\, dy$$

$$= \pi k^3 \sinh 2x/k - 2\pi k^2 x$$

or, writing $e^{x/k}$ as $X$,

$$V = \tfrac{1}{2}\,\pi k^3\,(X^{2} - X^{-2}) - 2\pi k^2 x.$$

The second integration gives the area of the surface of revolution generated by the rotation of an equidistant curve about its projection on a plane. If, with reference to the figure shown below, AB = *p* and AC = *q*, implying that CD = *p cosh q/k* as proved in Section 27, then the area *A* of the surface generated by the rotation of arc CD about AB is given by

$$\boldsymbol{A = \tfrac{1}{2}\,\pi\, pk\,(e^{2q/k} - e^{-2q/k})}.$$

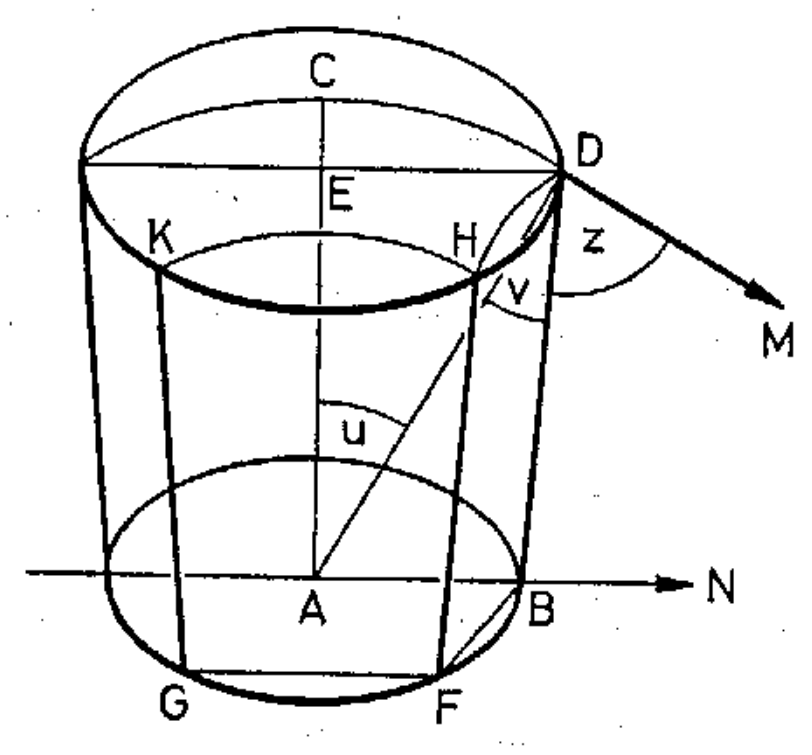

The justification for this is as follows. An element of arc *ds* of CD is

$$ds = \cosh q/k \; dp.$$

The circumference *C* of the vertical circle of radius *q*, according to Section 30, is

$$C = 2\pi k \;\sinh q/k.$$

Then the element of surface area *dA* generated by rotating this element of arc CD about AB is given by

$$dA = 2\pi k \cosh q/k \,.\, \sinh q/k \; dp$$

$$= \tfrac{1}{2}\,\pi k\,(e^{2q/k} - e^{-2q/k})\; dp$$

and so, by integration,

$$A = \tfrac{1}{2}\,\pi pk\,(e^{2q/k} - e^{-2q/k})\,.$$

The third integration is frankly rather puzzling. Bolyai's original Latin text gives for the volume *V* of the cylindrical solid generated by rotating CABD about AB

$$\boldsymbol{V = \tfrac{1}{2}\;\pi k^2 p\,(Q^2 + Q^{-2})},$$

whereas the volume given in Kárteszi's edition is

$$\boldsymbol{V = \tfrac{1}{4}\,\pi k^2 p\,(Q - Q^{-1})^2},$$

which is different.

It would appear that the modern version is correct. The volume of this solid can be regarded as the sum of the surface areas of an infinite number of concentric cylindrical shells of diminishing radius. Since the surface area of the solid is purely a function of the radius of revolution *q* times a constant *p*, the volume can be obtained as the integral of this function with respect to *q*.

The surface area, as shown above, is ½ $\pi kp \sinh 2q/k$. Therefore the volume is

$$V = \tfrac{1}{2} \pi kp \int_0^q \sinh 2z/k \, dz$$

$$= \tfrac{1}{4} \pi k^2 p (\cosh 2q/k - 1)$$

$$= \tfrac{1}{4} \pi k^2 p \sinh^2 q/k$$

$$= \tfrac{1}{4} \pi k^2 p (Q - Q^{-1})^2,$$

as stated in the modern version.

Bolyai points out in this subsection that the identities of hyperbolic geometry *S* coincide with those of Euclidean geometry Σ if the linear constant *k* is taken as infinitely great (see the note on Section 31).

## §33

It remains to make clear the implications of our theory, as promised in **§32.**

(I) It rests undecided whether system *Σ* or some system **S** is true in actual fact.

(II) Everything we have deduced from the hypothesis that Axiom XI is false is absolutely valid (always in the sense of §32) and therefore, in this sense, does not lean on any assumption. So, there is an a priori plane trigonometry in which only the true system is unknown, hence only the absolute magnitudes of the expressions remain undetermined, but on the basis of a single known case, obviously, the whole system could be fixed. On the other hand, spherical trigonometry was established in **§26** in an absolute way. Further, on the surface **F** we have a geometry completely analogous to the plane geometry of system *Σ*.

(III) If we knew that *Σ* is valid, there would remain no open question in this respect. On the other hand, if we knew that *Σ* is not valid, then (see §31) starting e.g. from the sides *x*, *y* and the angle between them, each concretely given, it would obviously be impossible in itself and by itself to resolve the triangle in an absolute manner, i.e. a priori determine the remaining angles and the ratio of the third side to both of the given ones; this could only be done by determining *X* and *Y*, for which purpose we would need a concrete value of *a* whose corresponding *A* is known. In the latter case, *k* would be a natural unit of length (just as *e* is the base of natural logarithms). Assuming the existence of *k*, we shall show how it can be constructed, at least for practical use, as accurately as possible.

(IV) In the sense of (II) and (III), all spatial problems can apparently be settled by a recent method of analysis which deserves high appreciation if applied within proper limits.

(V) Finally, there comes something not at all disagreeable to the gentle reader: assuming that S, and not *Σ*, is valid in reality, we construct a rectilinear figure of area equal to that of a circle.

***Note on 33***

In this section Bolyai makes several general points:

1. It is still to be determined whether the geometry of space is Euclidean or whether it accords with some other system.

2. Everything deduced from the assumption that the parallel postulate is false is valid regardless of whether or not space is Euclidean. A single known case would establish the validity of the whole system. Spherical geometry, however, is independent of the parallel postulate and is therefore part of absolute or neutral geometry. The geometry of the horosphere is Euclidean.

3. Even if hyperbolic geometry were proved to be true, it would still be impossible to solve a triangle given two sides and the enclosed angle without knowing the value of the linear constant.

4. All problems of area and volume in hyperbolic geometry can be solved by proper use of the calculus.

5. If hyperbolic geometry is true, then it is possible to square the circle .

## §34

*From the point D, a half-line $\overrightarrow{DM}$ parallel to $\overrightarrow{AN}$ can be drawn in the following way.*

Let $DB$ be the perpendicular from $D$ to $AN$. At any point $A$ of the line $AB$, draw $AE$, the perpendicular to $AN$ in the plane $ABD$. Further, let

$$DE \perp AE.$$

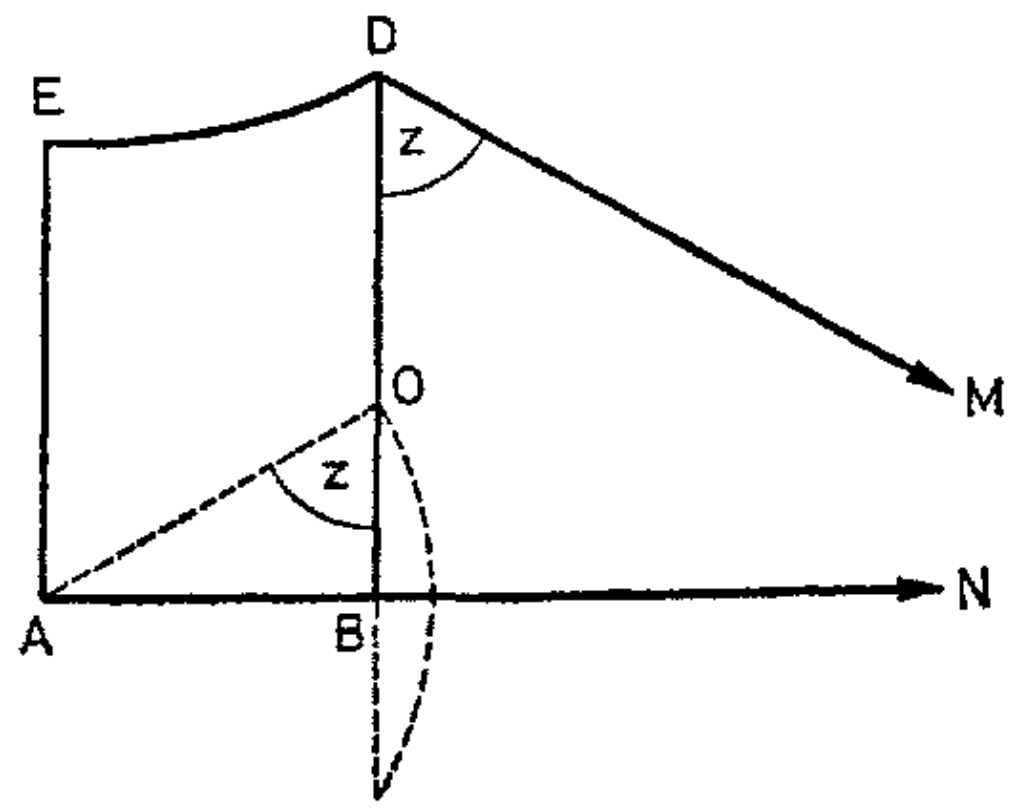

**Figure 23**

Now, if

$$\overrightarrow{DM} \parallel \overrightarrow{BN},$$

then by §27

$$\circ\overline{ED} : \circ\overline{AB} = 1 : \sin z.$$

But $\sin z \leq 1$, so that $\overline{AB} \leq \overline{DE}$. Therefore the quadrant of radius $\overline{DE}$ drawn from the centre $A$ in the angular domain $BAE$ has a point $B$ or $O$ in common with the half-line $\overrightarrow{BD}$. In the first case evidently

$$z = R.$$

In the second case, by §25,

$$(\circ\overline{AO} = \circ\overline{ED}) : \circ\overline{AB} = 1 : \sin AOB;$$

hence

$$z = \sphericalangle AOB.$$

Thus, if we take the angle $z$ to be equal to $\sphericalangle AOB$ then

$$\overrightarrow{DM} \parallel \overrightarrow{BN}.$$

***Note on 34***

The last part of the *Appendix* is devoted to constructions. The first of these is **a (neutral) compass and straightedge construction to draw the parallel to a given line.**

In the figure AN is the given straight line and D is the given point. Bolyai's construction of the parallel DM is as follows. First draw DB perpendicular to AN. Next take any point A on AN and draw AE perpendicular to AN (where E is yet to be determined). Then draw DE perpendicular to AE (which determines E). In the hyperbolic plane ABDE will be a quadrilateral with right angles at B, A and E and an acute angle at D. Consequently ED will be greater than AB. Next draw AO = ED. The segment AO, being greater than AB, will intersect BD at a point O, as shown. Then ∠AOB = *z* will be the angle of parallelism corresponding to the segment BD. All that remains is to draw ∠BDM = *z* and draw DM, which is the required parallel to AN.

The proof of this construction relies on the result obtained in Section 10 on the rectification of the equidistant. This implies that, if DM || BN,

$$^{\circ}ED : {}^{\circ}AB = 1 : \sin z.$$

On the other hand, in triangle AOB the neutral sine rule implies that

$$^{\circ}AO : {}^{\circ}AB = 1 : \sin \angle AOB.$$

But AO = ED by construction. Therefore a comparison of the two equations confirms that ∠AOB = *z*, the magnitude of the angle of parallelism corresponding to the segment BD.

It is important to note that this construction is valid both in Euclidean and hyperbolic geometry. If the parallel postulate is assumed to be true, ABDE will be a rectangle and ED will be equal to AB. So the point O will coincide with B and the angle of parallelism will be a right angle.

## §35

If S is valid in reality, *a straight line that is perpendicular to one arm of an acute angle and parallel to the other can be constructed as follows.*

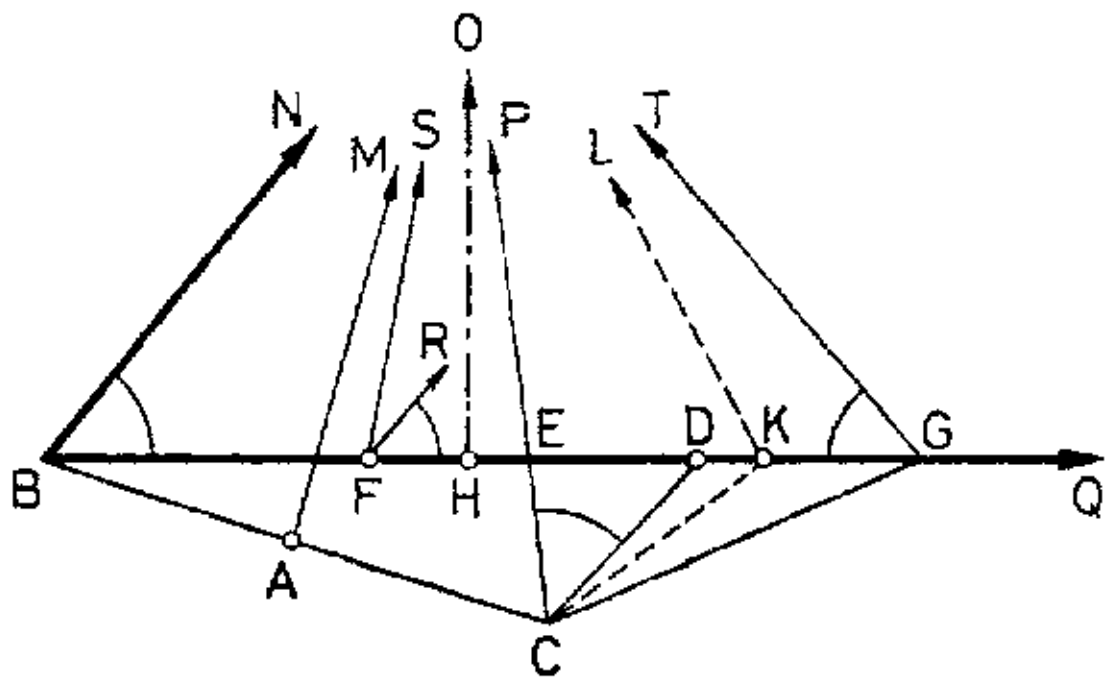

**Figure 24**

Let $AM \perp BC$. Choose the distance $\overline{AB} = \overline{AC}$ to be so small (see §19) that if constructing (by the aid of §34) the half-line $\overrightarrow{BN}$ parallel to $\overrightarrow{AM}$, the angle $ABN$ be greater than the angle originally given. Further, making use of §34 once more, draw the half-line $\overrightarrow{CP}$ that is parallel to $\overrightarrow{AM}$, and let each of $\sphericalangle NBQ$ and $\sphericalangle PCD$ be equal to the given angle. The half-lines $\overrightarrow{BQ}$ and $\overrightarrow{CD}$ intersect each other.

For let $\overrightarrow{BQ}$ (which lies in $\sphericalangle NBC$ by construction) intersect $\overrightarrow{CP}$ at the point $E$. Since $\overrightarrow{BN} \parallel \simeq \overrightarrow{CP}$, we have

$$\sphericalangle EBC < \sphericalangle ECB$$

and, consequently,

$$\overline{EC} < \overline{EB}.$$

Let

$$\overline{EF} = \overline{EC}, \quad \sphericalangle EFR = \sphericalangle ECD, \quad \text{and} \quad \overrightarrow{FS} \parallel \overrightarrow{EP}.$$

Then $\overrightarrow{FS}$ must fall into $\sphericalangle BFR$. Indeed, as

$$\overrightarrow{BN} \,]\, \overrightarrow{CP}$$

so that

$$\overrightarrow{BN} \parallel \overrightarrow{EP} \quad \text{and} \quad \overrightarrow{BN} \parallel \overrightarrow{FS},$$

by §14 we have

$$\sphericalangle FBN + \sphericalangle BFS < 2R = \sphericalangle FBN + \sphericalangle BFR$$

and, consequently,

$$\sphericalangle BFS < \sphericalangle BFR.$$

Thus $\overrightarrow{FR}$ intersects $\overrightarrow{EP}$, and therefore also $\overrightarrow{CD}$ intersects $\overrightarrow{EQ}$ at some point $D$.

Now let

$$\overline{DG} = \overline{DC}$$

and

$$\sphericalangle DGT = \sphericalangle DCP = \sphericalangle GBN.$$

Since $CD \simeq GD$, it follows that

$$BN \simeq GT \simeq CP.$$

If $K$ is the point on $\overrightarrow{BQ}$ (see §19) and $\overrightarrow{KL}$ is an axis of the L-line that corresponds to the axis $\overrightarrow{BN}$, then

$$BN \simeq KL.$$

Hence

$$\sphericalangle BKL = \sphericalangle BGT = \sphericalangle DCP$$

and, on the other hand,

$$KL \simeq CP.$$

So it is obvious that $K$ coincides with $G$ and

$$\overrightarrow{GT} \parallel \overrightarrow{BN}.$$

If, however, $\overrightarrow{HO}$ perpendicularly bisects $\overline{BG}$, then the half-line $\overrightarrow{HO}$ parallel to $\overrightarrow{BN}$ is constructed.

***Note on 35***

In this section Bolyai provides a straightedge and compass construction to draw the segment BH corresponding to a given angle of parallelism ∠NBQ and so draw the perpendicular OH parallel to BN.

## §36

Suppose there is given a half-line $\overrightarrow{CP}$ and a plane $MAB$. Let $CB$ be perpendicular to $MAB$, $BN$ be perpendicular to $BC$ in the plane $BCP$, and $\overrightarrow{CQ}$ be parallel to $\overrightarrow{BN}$ (see §34). If $\overrightarrow{CP}$ lies in $\sphericalangle BCQ$ then, in the plane $CBN$, the point of intersection of $\overrightarrow{CP}$ with $\overrightarrow{BN}$, and hence with $MAB$, can be found.

Next suppose we are given two planes $PCQ$ and $MAB$. Let

$$CB \perp MAB, \quad CR \perp PCQ,$$

and in the plane $BCR$

$$BN \perp BC, \quad CS \perp CR.$$

Then $BN$ and $CS$ lie in $MAB$ and $PCQ$, respectively. Taking the point of intersection of $BN$ and $CS$ (if it exists), the perpendicular at this point to $CS$ in $PCQ$ will apparently be the intersection of $MAB$ and $PCQ$.

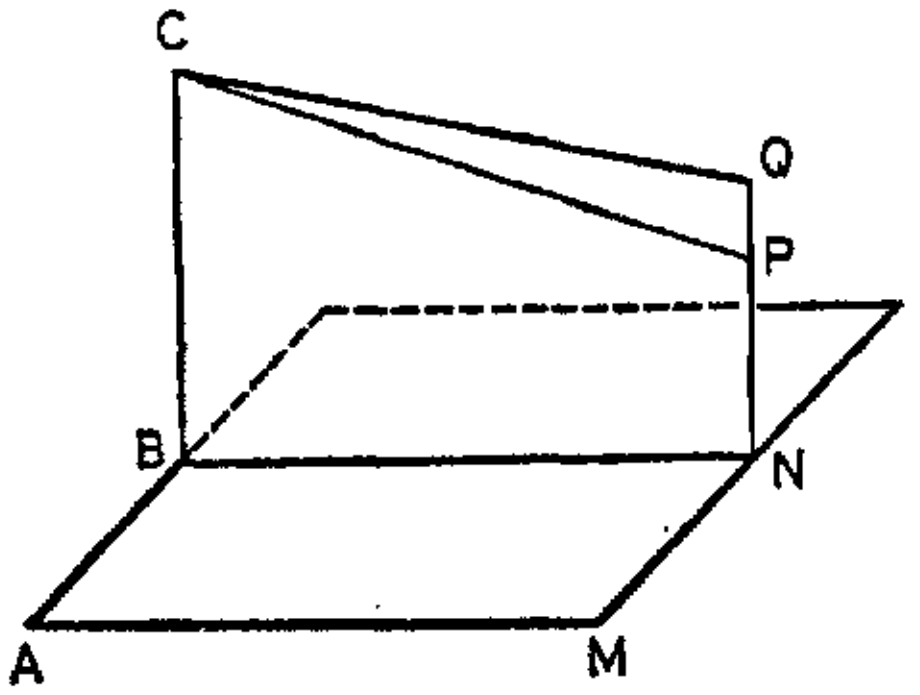

**Figure 25**

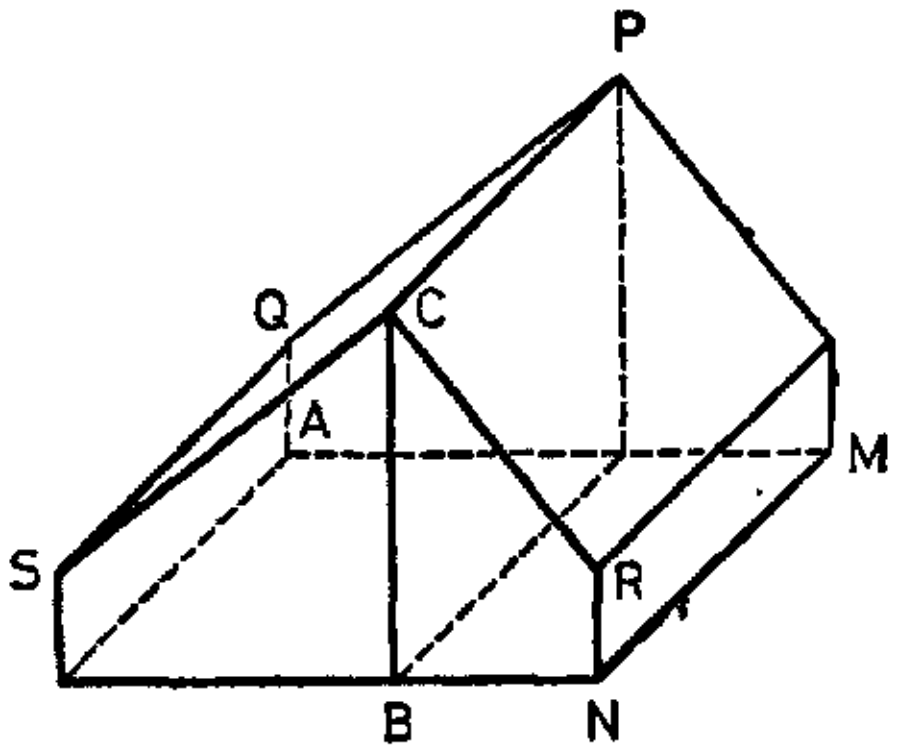

**Figure 26**

***Note on 36***

This section provides **a compass and straightedge construction to obtain the point where a given straight line CP intersects a given plane MAB and also to obtain the line of intersection between two planes MAB, PCQ.**

## §37

*On a half-line $\overrightarrow{AM}$, which is parallel to $\overrightarrow{BN}$, the point A can be chosen so that $AM \simeq BN$.*

For let us construct (see §34) a half-line $\overrightarrow{GT}$ that is parallel to $\overrightarrow{BN}$ and not in the plane *NBM*. Let

$$BG \perp GT, \quad \overline{GC} = \overline{GB}, \quad \text{and} \quad \overrightarrow{CP} \parallel \overrightarrow{GT}.$$

Choose the half-plane $|TG|D$ so as to form an angle with $|TG|B$ equal to the angle between $|PC|A$ and $|PC|B$. With the help of §36, determine the intersection *DQ* of the half-planes $|TG|D$ and $|NB|A$. Finally, let

$$BA \perp DQ.$$

Then the triangles formed by the **L**-lines which arise in the **F**-surface corresponding to $\overrightarrow{BN}$ are similar by §21. Therefore obviously

$$\overline{DB} = \overline{DA} \quad \text{and} \quad AM \simeq BN.$$

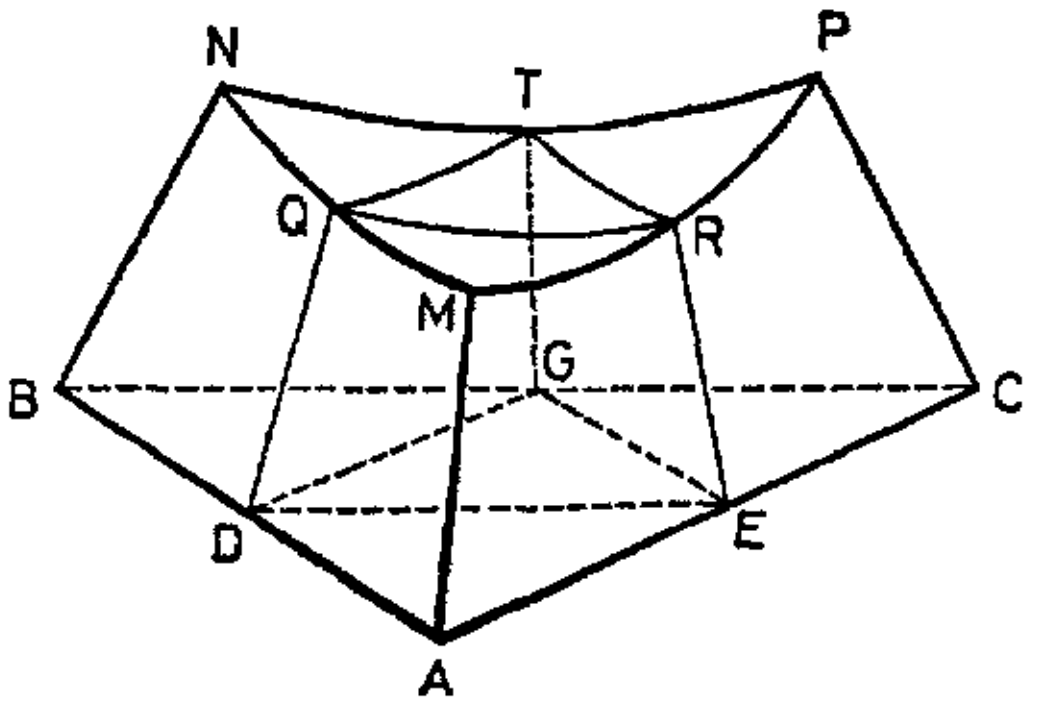

**Figure 27**

Hence it is easy to see that, *for* **L**-*lines given only by their endpoints, the endpoints of the fourth proportional and the geometric mean can be obtained, and all constructions which are, in system Σ, possible in the plane can so, without Axiom XI, be performed in the surface* **F**.

Thus, for instance, $4R$ can be divided into any number of equal parts by geometric construction if that division can be accomplished in system Σ.

### *Note on 37*

Given two parallel lines BN and AM, this section provides **a construction to obtain the point A on AM corresponding to B on BN.**

## §38

Let us construct, say, the angle $\sphericalangle NBQ = \frac{1}{3} R$ with the help of §37. Let $\overrightarrow{AM}$ be perpendicular to $\overrightarrow{BQ}$ and parallel to $\overrightarrow{BN}$ in system S (see §35). Determine $J$ by §37 so that $JM \simeq BN$. Then for the distance $\overline{JA} = x$ §28 yields

$$X = 1 : \sin \frac{1}{3} R = 2.$$

Thus $x$ is geometrically constructed.

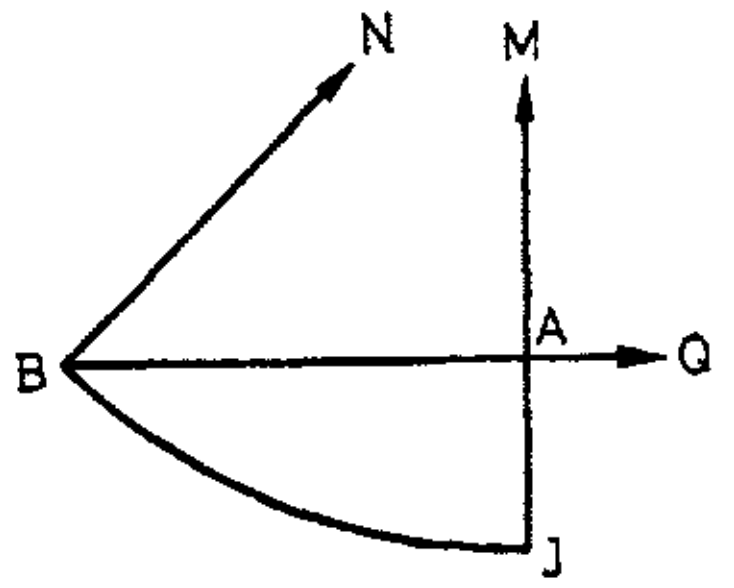

**Figure 28**

The angle $NBQ$ can be computed so that the difference between $\overline{JA}$ and $k$ be arbitrarily small. Indeed, it is sufficient to assure the validity of the relation

$$\sin NBQ = \frac{1}{e}.$$

***Note on 38***

In this construction Bolyai demonstrates **a method of constructing the radial distance *x* corresponding to a given ratio *X* between two concentric horocyclic arcs**:

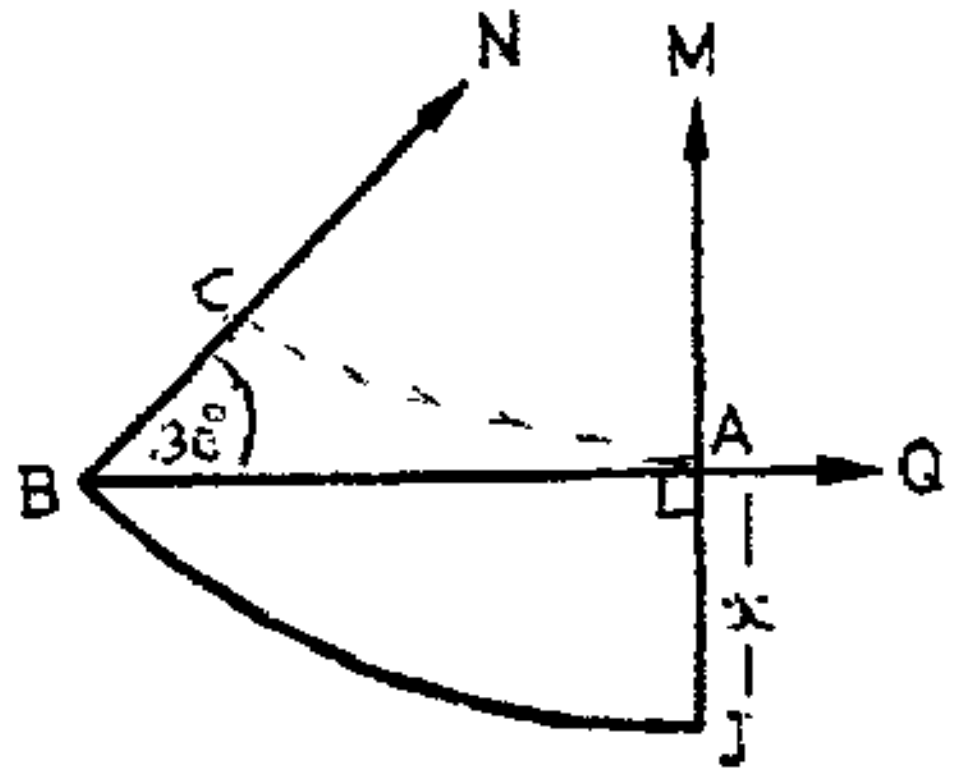

Suppose it is required to obtain the radial distance $x$ such that the ratio of two concentric arcs is, say, 2 : 1. On BN construct ∠NBQ = 30°. Using the construction given in Section 35, draw the segment BA on BQ corresponding to the angle of parallelism ∠NBA. Next draw AM perpendicular to BA. Thus AM || BN. Finally use the construction given in Section 37 to obtain the point J on AM corresponding to B. Then JB = $x$ is the required radial distance such that arc JB = 2 arc AC.

The proof follows from Section 28:

$$\text{arc JB} : \text{arc AC} = X$$

$$= 1 : \sin R/3$$

$$= 2 : 1.$$

However Bolyai points out that while it is possible to compute to any degree of accuracy the angle of parallelism ∠NBQ corresponding to a ratio $X = e$, which would occur when the radial distance $x = k$ (the linear constant), the construction would remain approximate.

## §39

If, in the plane, the curves $\widehat{PQ}$ and $\widehat{ST}$ are parallel to the straight line $MN$ (see §27), and the distances $\overline{AB}$ and $\overline{CD}$ are perpendicular to $MN$ and equal to each other, then evidently

$$\triangle DEC \equiv \triangle BEA.$$

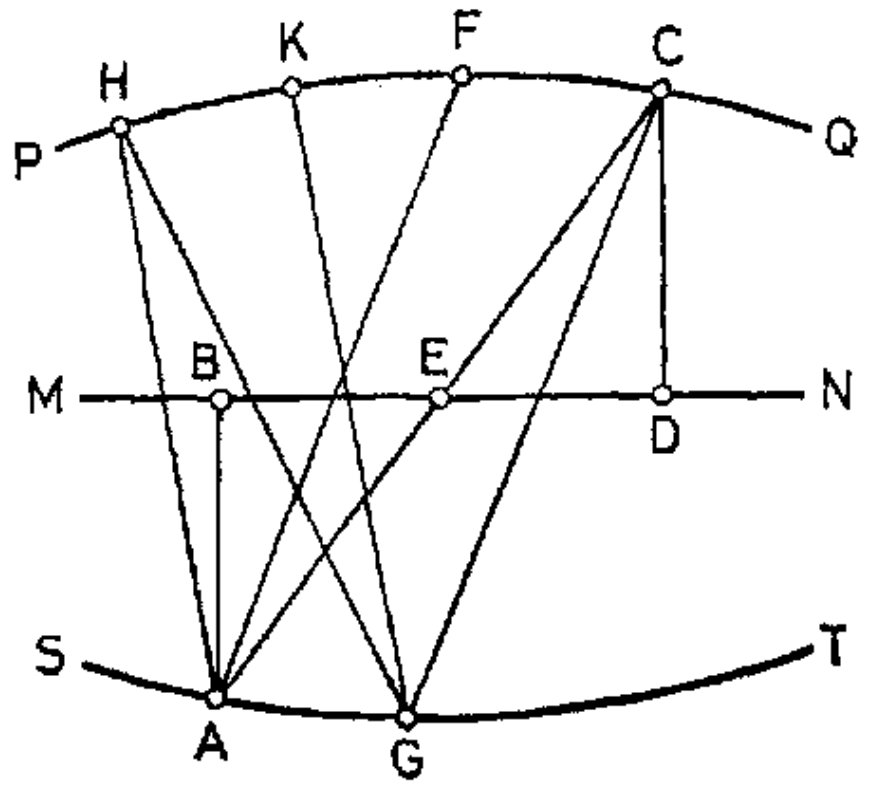

**Figure 29**

Therefore the angular domains $\sphericalangle ECP$ and $\sphericalangle EAT$ (which may be bounded by lines of mixed kind) are congruent, and

$$\overline{EC} = \overline{EA}.$$

Further, if

$$\widehat{CF} = \widehat{AG}$$

then

$$\triangle ACF \equiv \triangle CAG$$

and each of these triangles is half of the quadrilateral $FAGC$. If $FAGC$ and $HAGK$ are two quadrilaterals of this type based on $\widehat{AG}$ and lying between $\widehat{PQ}$ and $\widehat{ST}$, then (following Euclid) it can be seen that their areas are equal as are those of the triangles $AGC$ and $AGH$ erected on the same arc $\widehat{AG}$ and having one vertex on $\widehat{PQ}$.

Moreover,

$$\sphericalangle ACF = \sphericalangle CAG, \quad \sphericalangle GCQ = \sphericalangle CGA,$$

and by §32

$$\sphericalangle ACF + \sphericalangle ACG + \sphericalangle GCQ = 2R.$$

Thus

$$\sphericalangle CAG + \sphericalangle ACG + \sphericalangle CGA = 2R.$$

Consequently, in any such triangle $ACG$ the sum of the three angles is $2R$.

Whether the straight line $AG$ is in the curve $\widehat{AG}$ (the latter being parallel to $MN$) or not, it is now clear that the areas as well as the angle sums of the rectilinear triangles $AGC$, $AGH$ are equal.

***Note on 39***

The aim of this section is to a prove a theorem needed in the next sections. The theorem demonstrates that **if PQ and ST are equidistant curves situated symmetrically on either side of a straight line MN, and if the points A and G on ST are fixed while the point C on PQ can move along the curve, then the angle sum of triangle ACG is two right angles.**

## §40

*Two (from now on, rectilinear) triangles ABC and ABD of equal area and with one side equal have equal angle sums.*

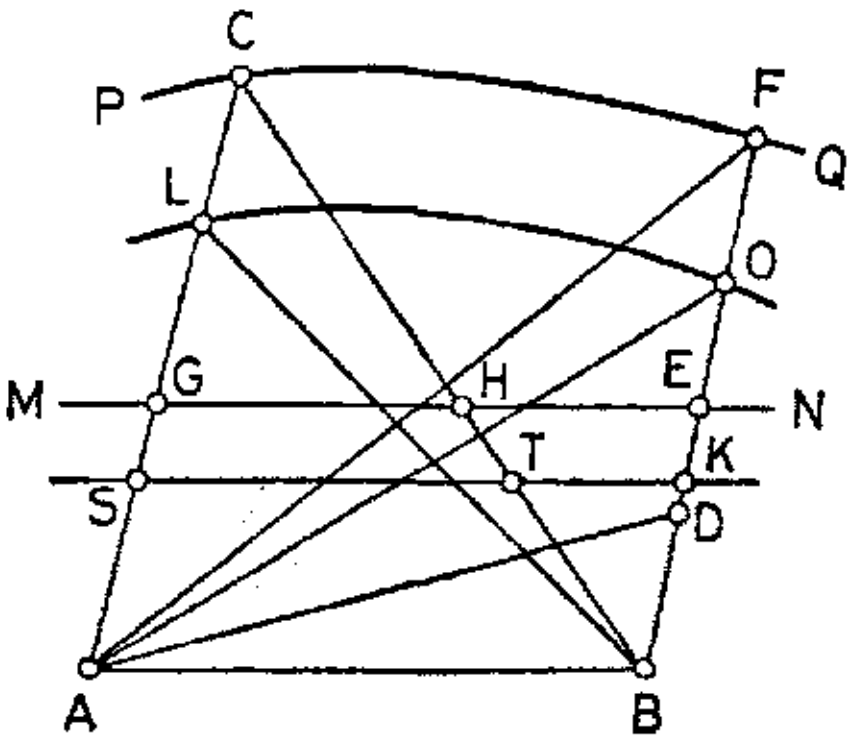

**Figure 30**

For let $MN$ bisect $\overline{AC}$ as well as $\overline{BC}$, and let the curve $\overset{\frown}{PQ}$ passing through the point $C$ be parallel to $MN$. Then $D$ lies on $\overset{\frown}{PQ}$.

In fact, if the half-line $\overrightarrow{BD}$ intersects the line $MN$ at the point $E$ and therefore, by §39, it intersects $\overset{\frown}{PQ}$ at distance $\overline{EF} = \overline{BE}$, then

$$\triangle ABC = \triangle ABF$$

and consequently

$$\triangle ABD = \triangle ABF;$$

thus $D$ coincides with $F$.

On the other hand, if $\overrightarrow{BD}$ does not intersect $MN$, let $C$ be the point where the perpendicular bisector of $\overline{AB}$ intersects $\overset{\frown}{PQ}$. Choose $\overline{GS} = \overline{HT}$ so that the line $ST$ intersects the produced line $BD$ at a point $K$ (the possibility of this choice follows as in §4). Moreover, let

$$\overline{SL} = \overline{SA}, \quad \overset{\frown}{LO} \parallel ST,$$

and let $O$ be the intersection of the straight line $BK$ and the curve $\overset{\frown}{LO}$. Then by §39

$$\triangle ABL = \triangle ABO$$

and therefore

$$\triangle ABC > \triangle ABD,$$

which contradicts the hypothesis.

***Note on 40***

This theorem proves that **if two triangles, ABC, ABD, having one side in common, are equal in area, then their angle sums are equal.**

## §41

*The angle sums of two triangles ABC and DEF of equal area are equal.*

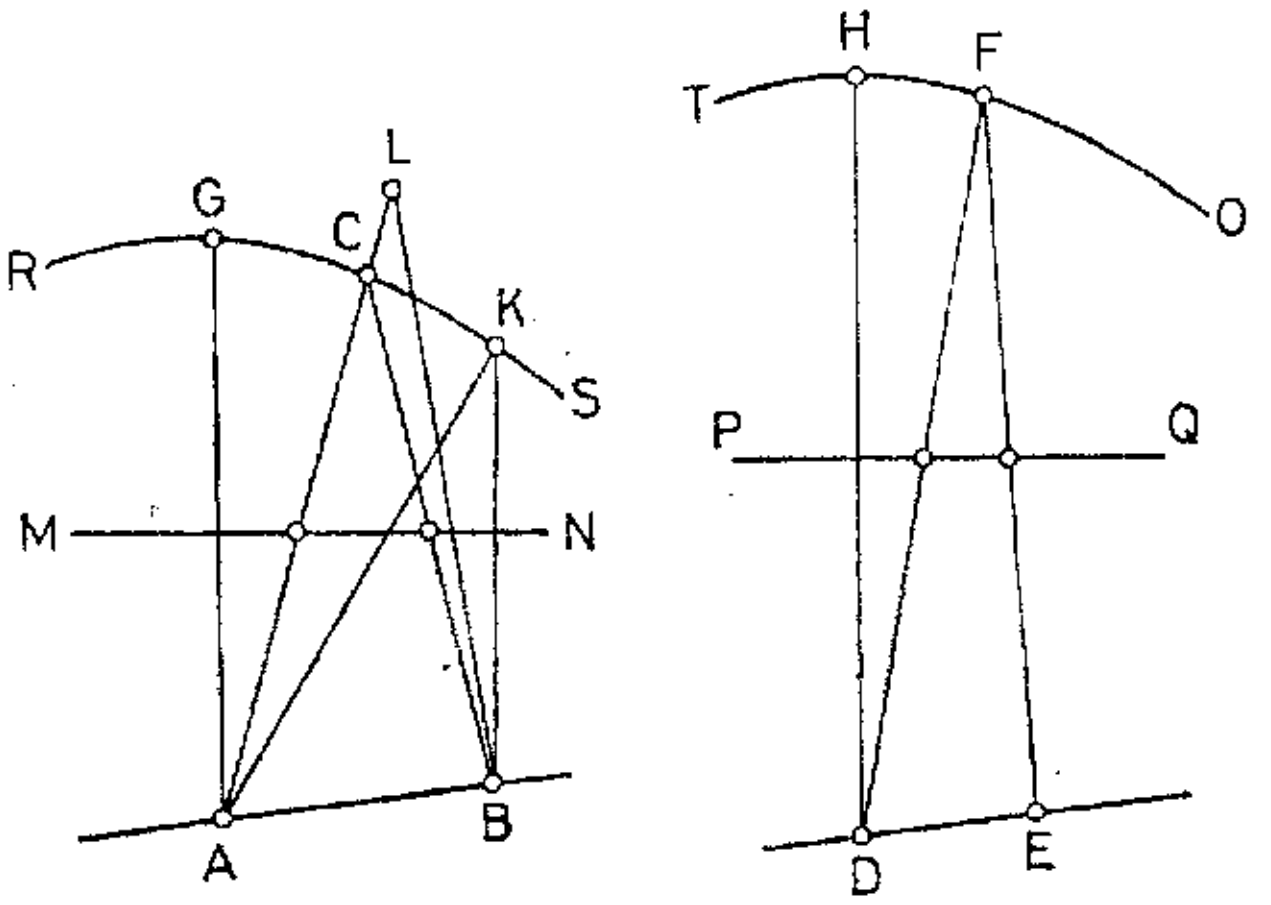

**Figure 31**

For let $MN$ bisect both $\overline{AC}$ and $\overline{BC}$, and $PQ$ bisect both $\overline{DF}$ and $\overline{EF}$. Further let

$$\widehat{RS} \parallel MN \quad \text{and} \quad \widehat{TO} \parallel PQ.$$

The distance $\overline{AG}$, perpendicular to $\widehat{RS}$, is either equal to the distance $\overline{DH}$, perpendicular to $\widehat{TO}$, or one of them, say $\overline{DH}$, is greater. In each case, the circle with centre $A$ and radius $\overline{DF}$ has a point $K$ in common with $\widehat{GS}$. Then by §39

$$\triangle ABK = \triangle ABC = \triangle DEF.$$

But the angle sum of $\triangle AKB$ is by §40 equal to that of $\triangle DFE$, and by §39 equal to that of $\triangle ABC$. Hence also $\triangle ABC$ and $\triangle DEF$ have equal angle sums.

In system S, this theorem may be reversed. For let the angle sums of $\triangle ABC$ and $\triangle DEF$ be equal and

$$\triangle BAL = \triangle DEF.$$

By the foregoing, the angle sum of one of the latter triangles is equal to that of the other. Consequently, also the angle sum of $\triangle ABC$ is equal to that of $\triangle ABL$. Hence obviously

$$\sphericalangle BCL + \sphericalangle BLC + \sphericalangle CBL = 2R.$$

But according to §31 the angle sum, in system S, of any triangle is less than $2R$. Thus $L$ coincides with $C$.

***Note on 41***

Section 41 generalises the result of the last theorem to prove that **if two triangles have equal areas, they have equal angle sums. In hyperbolic geometry the converse is also true, that if two triangles have equal areas, their angle sums are equal.**

**§42**

*If u and v are the supplements to 2R of the angle sums in △ABC and △DEF, respectively, then*

$$\triangle ABC : \triangle DEF = u : v.$$

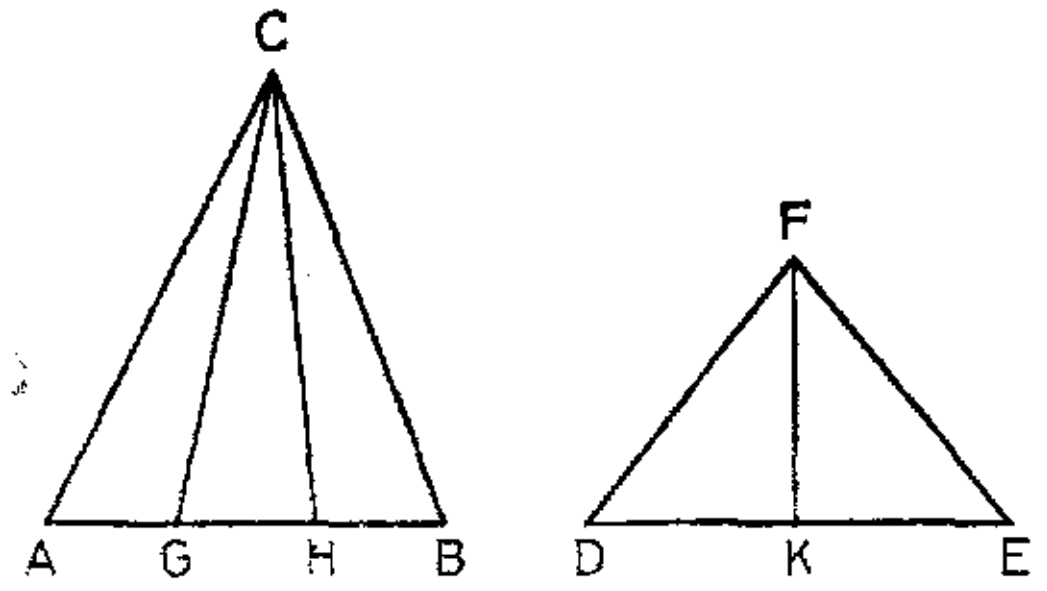

**Figure 32**

For let each of the triangles *ACG*, *GCH*, *HCB*, *DFK*, *KFE* have area *p* and

$$\triangle ABC = mp, \quad \triangle DEF = np.$$

Further let the angle sum of any triangle of area *p* be equal to *s*. Then manifestly

$$2R - u = ms - (m-1)2R = 2R - m(2R - s)$$

i.e.

$$u = m(2R - s),$$

and similarly

$$v = n(2R - s).$$

Thus

$$\triangle ABC : \triangle DEF = m : n = u : v.$$

It is easy to see that this extends to the case of triangles *ABC* and *DEF* whose areas are incommensurable.

It can be proved similarly that the areas of two spherical triangles are to each other as are the excesses over 2*R* of their angle sums. If two angles of a spherical triangle are right angles, then the third angle *z* is just the excess mentioned above. But the area of this spherical triangle is by §32, (VI) equal to

$$\frac{z}{2\pi} \frac{p^2}{2\pi},$$

where *p* denotes the circumference of the great circle. Consequently, the area of any spherical triangle with excess *z* is

$$\frac{zp^2}{4\pi^2}.$$

***Note on 42***

In hyperbolic geometry the sum of the angles of any triangle is strictly less than two right angles. The amount by which its angle sum falls short of two right angles is called its defect. Bolyai shows in this section and the next that **the area of any triangle is given by the product of its defect and the square of the linear constant *k*.**

His proof is in two parts. In Section 42 he proves first that the area of a triangle is proportional to its defect.

Suppose that the defect of triangle ABC = $u$ while the defect of triangle DEF = $v$. Suppose also that triangle ABC is divided into $m$ equal triangles, each of area $p$, while triangle DEF is divided into $n$ equal triangles, each of area $p$. Finally assume that the angle sum of any triangle of area $p$ is equal to $s$.

Considering triangle ABC, it is clear that its angle sum is equal to the sum of the angles in its $m$ constituent triangles minus the $(m-1)$ sets of adjacent angles along the base. Therefore

$$2R - u = ms - (m-1)\,2R$$

i.e. $u = m\,(2R - s)$.

Similarly

$$v = n\,(2R - s).$$

Consequently

$$\text{area triangle ABC} : \text{area triangle DEF} = m : n$$
$$= u : v.$$

This result is equivalent to the statement that the area of a triangle is equal to its defect times a constant of proportionality.

## §43

We now express the area, in system S, of a rectilinear triangle in terms of the angle sum.

If $\overline{AB}$ increases indefinitely*, then by **§42** the ratio

$$\triangle ABC:(R-u-v)$$

remains constant. On the other hand, by **§32**, (V)

$$\triangle ABC \to BACN$$

and by §1

$$R-u-v \to z.$$

Hence

$$BACN:z = \triangle ABC:(R-u-v) = BAC'N':z'$$

Furthermore, **§30** obviously yields

$$BDCN:BD'C'N' = r:r' = \operatorname{tg} z:\operatorname{tg} z'.$$

But for $y' \to 0$

$$\frac{BD'C'N'}{BAC'N'} \to 1$$

and also

$$\frac{\operatorname{tg} z'}{z'} \to 1.$$

It follows that

$$BDCN:BACN = \operatorname{tg} z:z.$$

By **§32**, however,

$$BDCN = rk = k^2 \operatorname{tg} z.$$

Therefore

$$BACN = zk^2.$$

* See Fig. 20.

### *Note on 43*

In the first part of Section 43 Bolyai completes the proof on the area of a triangle by showing that **the constant of proportionality relating the area of a triangle to its defect is the square of the linear constant *k*.** His argument can be illustrated using the figure given in the course of his rectification of the horocycle in Section 30:

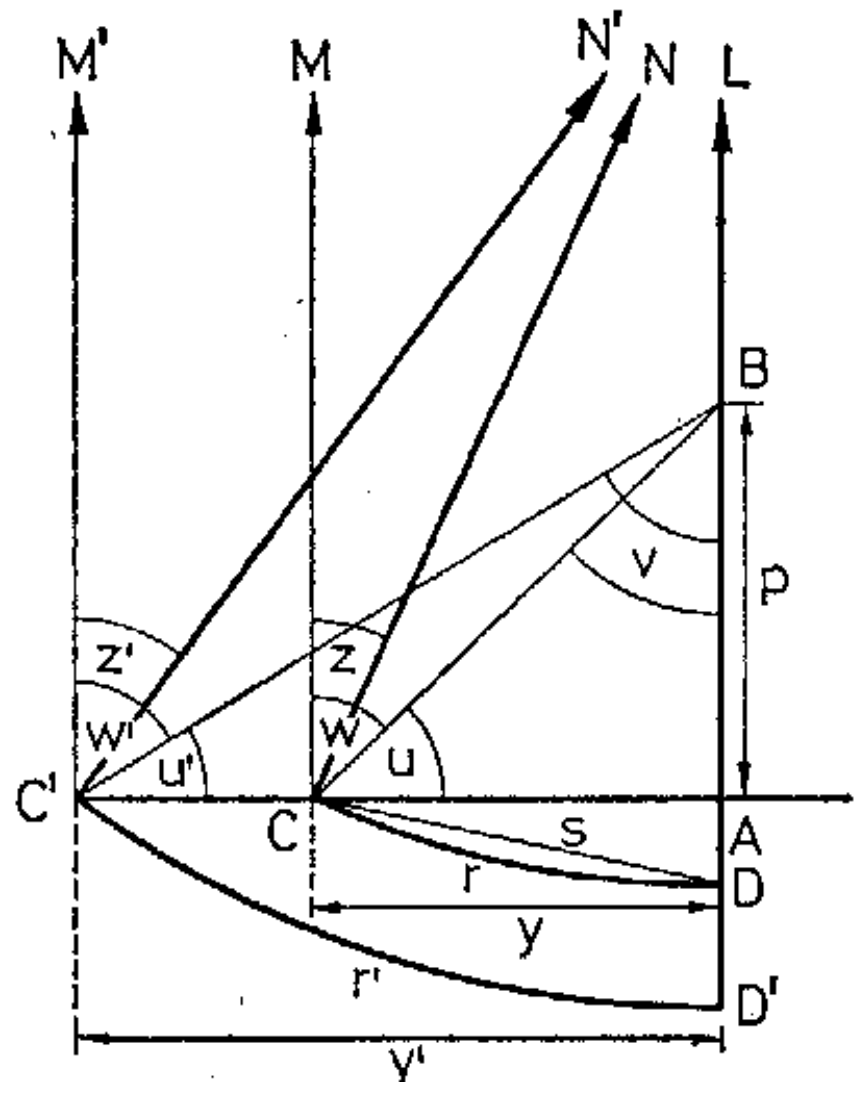

The strategy employed by Bolyai is first to express the ratio between the area of triangle ABC and its deficit in terms of the ratio between the area of the trilateral BACN and ∠MAN = $z$, the complement of the angle of parallelism corresponding to the segment AC. Then he expresses the area of BACN in terms of $z$ in order to evaluate this ratio.

As the point B recedes from A along AB to an infinitely distant point so that $p$ increases indefinitely,

$$\text{triangle ABC} \rightarrow \text{BACN}.$$

Moreover as CB approaches CN, $v \rightarrow 0$ and $u \rightarrow (R - z)$. Therefore

$$R - u - v \rightarrow z.$$

It follows that

$$\text{BACN} : z = \text{triangle ABC} : (R - u - v).$$

But given that ∠CAB is a right angle,

$$(R - u - v) = \text{deficit of triangle ABC},$$

so

$$\text{area triangle ABC} : \text{deficit of triangle ABC} = \text{BACN} : z.$$

Now

$$AC : AC' = (R - z) : (R - z') \quad \text{(segments corresponding to the angles of parallelism)}$$

$$= z : z',$$

therefore

$$BACN : BAC'N' = z : z'. \qquad \text{... (1)}$$

By the result given in Section 32, V,

$$\text{area } BDCN \to rk,$$

while

$$\text{area } BD'C'N' \to r'k,$$

so

$$BDCN : BDC'N' = r : r'$$

$$= \tan z : \tan z' \quad \text{(by the result given in Section 30).} \qquad \text{...(2)}$$

But as AC' = $y' \to 0$,

$$BD'C'N / BAC'N' \to 1 \qquad \text{... (3)}$$

and also

$$\tan z / \tan z' \to 1. \qquad \text{... (4)}$$

The four numbered equations imply that

$$BDCN : BACN = \tan z/z.$$

However, as shown above,

$$BDCN = rk$$

$$= k^2 \tan z.$$

Therefore

$$BACN = zk^2.$$

Finally

$$\text{area triangle } ABC : \text{deficit of triangle } ABC = BACN : z$$

$$= zk^2 : z$$

$$= k^2 : 1,$$

confirming that the constant of proportionality is $k^2$**. Thus the area of a triangle equals its deficit times the square of the linear constant *k*.**

## Section 43 (continued)

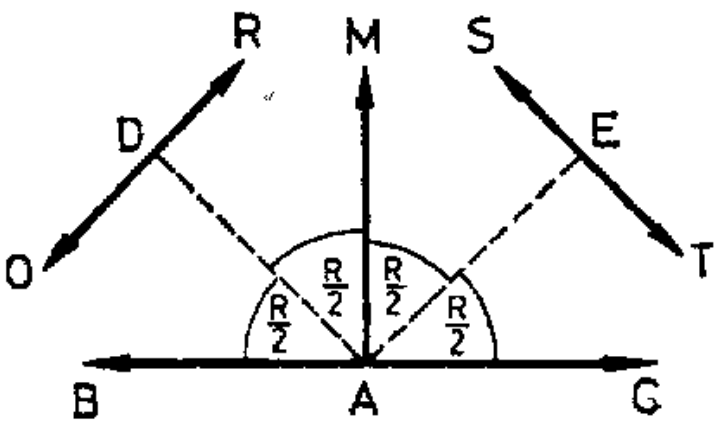

Figure 33

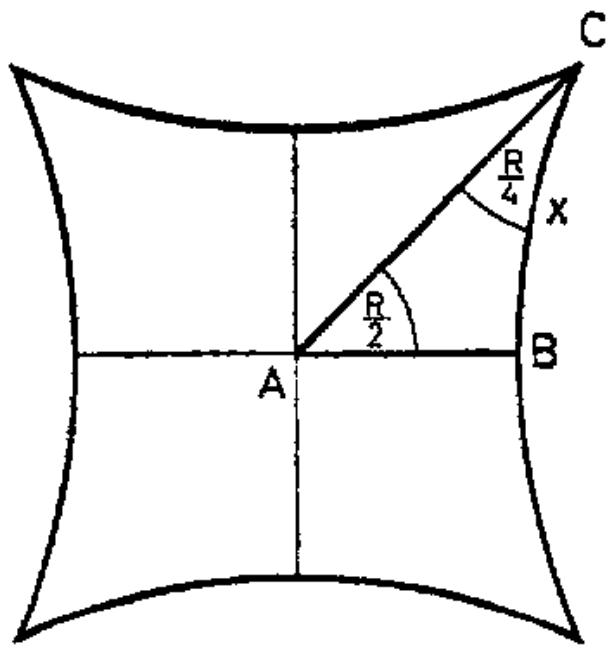

Figure 34

Thus, denoting from now onwards the area of any triangle the supplement to $2R$ of whose angle sum is $z$ simply by $\triangle$, we have

$$\triangle = zk^2.$$

Hence it easily follows that if*

$$\overrightarrow{OR} \parallel \overrightarrow{AM} \quad \text{and} \quad \overrightarrow{RO} \parallel \overrightarrow{AB}$$

then the area of the domain enclosed by the lines $OR$, $ST$ and $BC$, which is obviously the limit of the areas of indefinitely increasing rectilinear triangles i.e. the limit of $\triangle$ for $z \to 2R$, will be

$$\pi k^2 = \odot k$$

for a circular domain of the surface **F**. Denoting this limit by $\square$, from §30 and §21 we obtain

$$\pi r^2 = \mathrm{tg}^2 z \cdot \square = \odot r$$

for a circular domain of **F**. By §32, (VI) the latter value is equal to $\odot s$, where $s$ denotes the chord** $\overline{DC}$. Now if, perpendicularly bisecting the given radius $s$ of the circle in the plane (or the L-formed radius of the circle in the surface **F**) and leaning on §34, we construct a half-line $\overrightarrow{DB}$ with

$$\overrightarrow{DB} \parallel \simeq \overrightarrow{CN},$$

then by drawing the perpendicular $CA$ to $DB$ and the perpendicular $CM$ to $CA$ we obtain $z$. Hence, taking any L-formed radius for unit, tg $z$ can (by §37) be determined geometrically, with the help of two uniform arcs of the same curvature. If we know only their endpoints and construct their axes, such arcs can obviously be compared in measure just as straight line segments, and in this respect they may be considered equivalent to straight line segments.

Now a quadrilateral, say a square, of area $\square$ can be constructed as follows. Let***

$$\sphericalangle ABC = R, \quad \sphericalangle BAC = \frac{1}{2}R, \quad \sphericalangle ACB = \frac{1}{4}R,$$

and

$$\overline{BC} = x.$$

By §31, (II) $X$ can be expressed using merely square roots; by §37 it can be constructed as well. If we know $X$, we can determine $x$ with the help of §38 (or §29 and §35). Further it is clear that eight times the area of $\triangle ABC$ is equal to $\square$.

*Thus the geometrical quadrature of the plane circle of radius s is accomplished in terms of a rectilinear figure and uniform arcs of one and the same kind (the latter are equivalent, as concerns comparison, to segments); the complanation of a circular do-*

*main of* **F** *can be performed likewise. Consequently, either Axiom XI of* EUCLID *holds or the geometrical quadrature of the circle is possible;* though until now it has remained undecided which of these two cases takes place in reality.

If $\mathrm{tg}^2 z$ is an integer or a rational fraction whose denominator (after reducting the fraction to the simplest form) is either a prime number of the form $2^m+1$ (of which form is also $2=2^0+1$) or a product of any number of primes of this form where each prime, excepting 2 which alone may occur any number of times, appears as a factor only once, then by the theory of polygons due to the celebrated GAUSS (an outstanding discovery of our age and actually of all times) for these and only these values of $z$ even a rectilinear figure of area $\mathrm{tg}^2 z \cdot \square = \odot s$ can be constructed.

In fact, since the theorem of §42 can be easily extended to arbitrary polygons, division of the area $\square$ obviously requires dissection of $2R$. However, it is possible to show that this can be achieved in a geometrical way only under the condition stated. On the other hand, in all such cases the foregoing helps to reach the purpose easily. Moreover, if $n$ belongs to the Gauss class, then any rectilinear figure can be converted geometrically into an $n$-sided regular polygon of equal area.

To settle the matter in all respects, it has remained to prove the impossibility of deciding without any assumption whether $\Sigma$ or some (and which) **S** is valid. Nevertheless, we leave this to a more appropriate occasion.

***Note on 43 continued***

The final part of Section 43, with which the *Appendix* culminates, is devoted to Bolyai's astonishing proof that **if the parallel postulate is false, then it is possible to square the circle using just a compass and straightedge.**

In Euclidean geometry squaring the circle by these means is an impossible task. This is because it is not possible, using only compass and straightedge, to construct a segment of length $\sqrt{\pi}$, which is the side of a square equal in area to a unit circle. Only an algebraic number, that is to say, a number which is the root of a polynomial equation with rational coefficients, can be constructed with these tools. Unfortunately $\pi$ is not an algebraic number, though this fact was only proved in 1882 by Ferdinand von Lindemann.

But in hyperbolic geometry the task is made possible by the existence of two figures which have no equivalent in Euclidean geometry. One is the horosphere. Of particular importance is the fact that on the horosphere, as proved in Section 32, the area of a circular cap swept out by a horocyclic radius of length $k$ is $\pi k^2$.

The other is the triangle of maximum size, shown in Figure 33, defined as a triangle whose sides OR, ST, BC are mutually parallel and whose angle sum is zero, its vertices lying at infinitely distant points.

The existence of such a triangle follows from the fact that in hyperbolic geometry the defect of a triangle is proportional to its size. As its size increases, the angle sum decreases. But the angle sum clearly cannot cannot decrease below zero, at which point the maximum defect is two right angles. The existence of a maximum defect then implies the existence of a triangle of maximum size.

Now the area of a triangle is given by its defect multiplied by the square of the linear constant *k*. So in the case of the maximum triangle whose angle sum is zero and its defect $\pi$, its area is $\pi k^2$.

The fact that a circular figure, a cap on the horosphere, and a linear figure, the maximum triangle, have the same area, $\pi k^2$, gets round the need to construct a segment of length $\pi$, which is the stumbling block in Euclidean geometry. All that is required to obtain the quadrature of a given circle in hyperbolic geometry is

1) to determine the ratio of its area to a cap of radius *k* on the horosphere
2) to construct a square whose area bears the same ratio to the maximum triangle.

Bolyai outlines the method by which the first task can be completed. The figure shown below (which is not given in the original text but has been added in Halsted's translation) illustrates his construction to obtain the unique angle *z* associated with every plane circle of radius *s*:

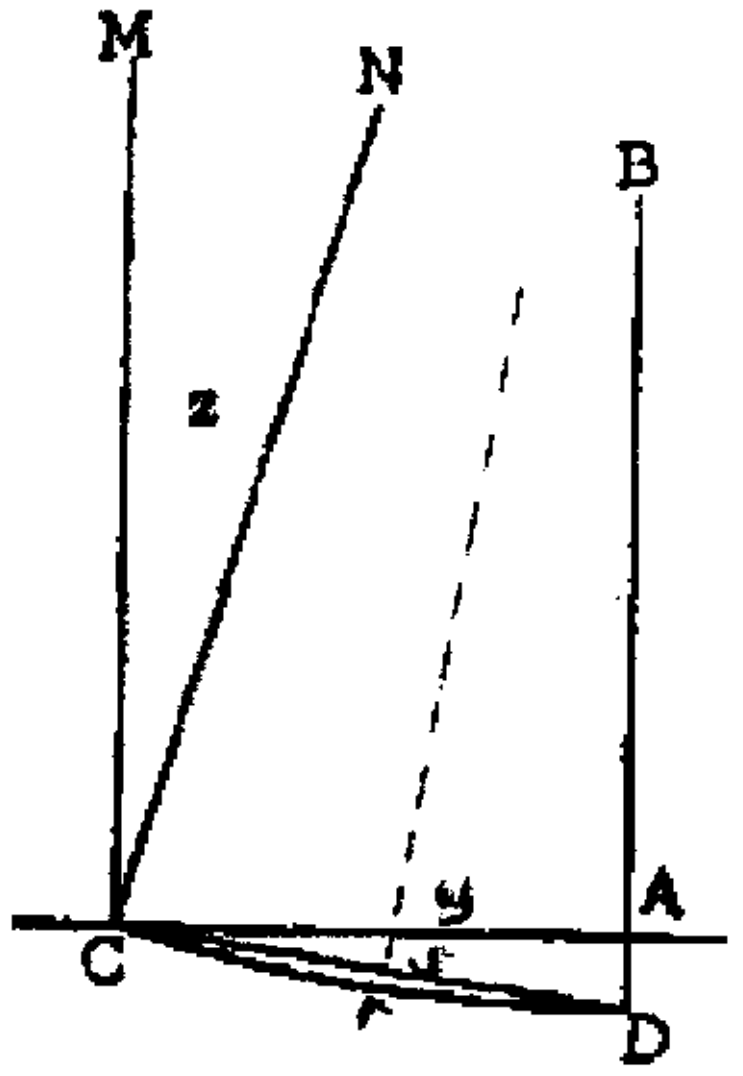

Draw the segment CD = *s* and draw the perpendicular bisector of CD (dotted). Using the construction given in Section 34, draw DB and CN parallel to this perpendicular. DB and CM will then be parallel to each another. Finally draw CA = *y* perpendicular to DB and CM perpendicular to CA. Then $\angle$MAN = *z* will be the required angle. It should be noted that the horocyclic arc CD = *r* is not part of the actual construction, but it can be supposed to pass through the corresponding points C and D.

The importance of the angle $z$ is that the area $A$ of a circle of radius $s$ is given by

$$A = \pi k^2 \tan^2 z.$$

Bolyai does not fully explain his derivation of this identity but it can be proved as follows. By the result given in Section 30,

$$r = k \sinh y/k,$$

while consideration of Figure 20 confirms that $s/2$ bears the same relation to $r/2$ as $y$ does to $r$, so

$$r/2 = k \sinh s/2k,$$

hence

$$\sinh s/2k = ½ \sinh y/k.$$

Considering the formula for the area of a circle of radius $s$ given in Section 32,

$$A = 4 \pi k^2 \sinh^2 s/2k$$

$$= 4 \pi k^2 ( ½ \sinh y/k)^2$$

$$= \pi k^2 \sinh^2 y/k$$

$$= \pi k^2 \tan^2 z \text{ (by the result given in Section 30).}$$

The term $\tan^2 z$ in this expression can be seen as a scaling factor which represents the ratio of the area of a plane circle of radius $s$ to the area $\pi k^2$ of a circular cap on the horosphere of radius $k$. In the case of the particular circle where the angle $z = \pi/4$ and $\tan^2 z = 1$, the scaling factor equals 1, so the area of this circle will be exactly equal to the area of the cap of radius $k$ and by extension to the area of the maximum triangle. All that remains to complete the quadrature of this circle is to construct a square equal in area to the maximum triangle.

Now a square in hyperbolic geometry has four equal sides and four equal acute angles. So if $w$ represents the magnitude of the angle of the square, the area of the square, in accordance with Bolyai's theorem, is given by its defect multiplied by $k^2$. The Euclidean formula for the angle sum of a square is $2\pi$, so its defect is $(2\pi - 4w)$. Therefore in the case of a square equal in area to the maximum triangle, Bolyai's formula from Section 42 requires that

$$k^2(2\pi - 4w) = \pi k^2,$$

which implies that $w = \pi/4 = 45°$. This square is shown in Figure 34 above.

The quadrature therefore requires the construction of the right angled triangle ABC with sides $a$, $b$, $c$ (where $a = x$), which represents 1/8 of the required square. The acute angles of the triangle, one equal to 45° and the other equal to 22½°, can be drawn by repeated bisection of a right angle using Euclid's neutral theorem (*Elements,* 1, 9) which is valid both in Euclidean and hyperbolic geometry. But to construct the triangle it is necessary to draw one of the sides, which is a more difficult task. Bolyai's instructions on how to do this are rather sketchy, but fortunately the noted Italian mathematician, Roberto Bonola, has provided a compass and straightedge construction.

His method depends on a preliminary construction, derived from Bolyai, to draw the straight line segment corresponding to a given angle of parallelism:

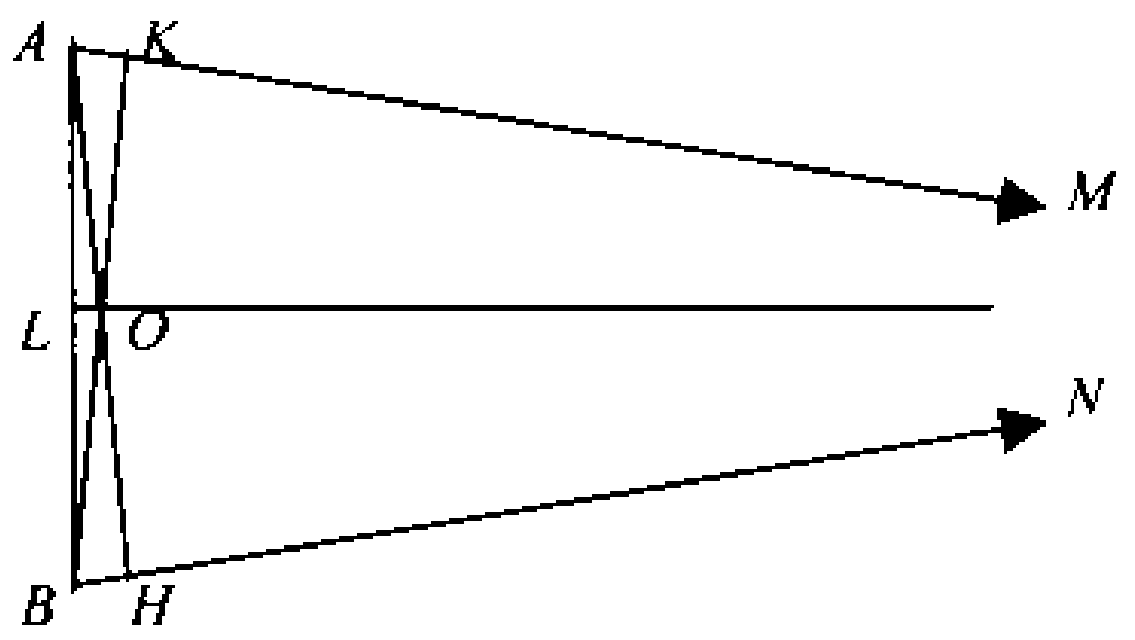

In order to draw the segment corresponding to the angle of parallelism ∠MAB, take a point B on AB such that the parallel BN to AM makes an acute angle ∠ABN with AB. Draw AH perpendicular to BN and BK perpendicular to AM. These two perpendiculars will meet at a point O. Draw OL perpendicular to AM. Then LO is parallel to AM and AL is the required segment.

The proof of this construction relies on the theorem that the perpendiculars drawn from each vertex of a triangle to the opposite side meet at a point. This theorem holds true even when one of the angles of the triangle, as here, lies at an infinitely distant point.

Now considering triangle ABC, the formula given in Section 31, I implies that

$$\cosh a/k = \cos \angle A / \sin \angle B$$
$$= \cos 45° / \sin 22 \tfrac{1}{2}°$$
$$= \sin 45° / \sin 22 \tfrac{1}{2}°.$$

At this point Bonola treats the angles 45° and 22 ½° as if they were angles of parallelism and uses the construction given above to draw the segments corresponding to these angles, denoting them by *b′* and *c′* respectively. In other words (using Lobachevsky's notation), 45° = $\Pi$ (*b′*) and 22 ½ ° = $\Pi$ (*c′*). Thus the previous equation can be written as

$$\cosh a/k = \sin \Pi (b') / \sin \Pi (c'). \qquad \ldots (1)$$

In the light of the identity (not given by Bolyai) that for any *y*

$$\cosh y/k = 1 / \sin \Pi (y),$$

equation (1) can be written as

$$\cosh a/k = \cosh c'/k / \cosh b'/k$$

or

$$\cosh a/k \,.\, \cosh b'/k = \cosh c'/k. \qquad \ldots (2)$$

Bonola now constructs an auxiliary right angled triangle A'B'C', where *c′* and b', obtained above, form, respectively, the hypotenuse and one of the sides.

The formula given by Bolyai in Section 31, III implies in respect of this triangle that

$$\cosh a'/k \,.\, \cosh b'/k = \cosh c'/k. \qquad \ldots (3)$$

A comparison of equations (2) and (3) shows that $a' = a$. In other words, the third side $a'$ of the auxiliary triangle A'B'C' is equal to the desired side $a$ of triangle ABC. This allows the construction of triangle ABC and therefore the construction by straightedge and compass of a square equal in area to the given circle.

The situation is more complicated when the scaling factor $\tan^2 z$ is a number other than 1. Bolyai observes that (expressed in modern notation) the quadrature of the circle is only feasible for those circles where $\tan^2 z$ is either an integer or a rational number whose denominator is of the form $2^m.\ p_1.\ p_2. \ldots$ , where $m$ is a non-negative integer and each $p_i$ is a distinct Fermat prime. A Fermat prime is a prime of the form $2^{2^{\wedge}n} + 1$ for $n$ = 0, 1, 2, 3, 4, … , though only these first five Fermat numbers are known to be prime.

The first case, when $\tan^2 z$ is an integer greater than 1, merits further examination. If, for example, $\tan^2 z = 2$, then the angle $w$ of a square equal in area to the given circle would be obtained from the equation

$$k^2(2\pi - 4w) = 2\pi k^2,$$

implying that $w = 0$. That is to say, the square in question would be the maximum square whose angle sum is zero and whose vertices lie at infinitely distant points, a figure which is obviously not constructible. Clearly if $\tan^2 z$ were an integer greater than 2, the angle of a square equal in area to the circle would be negative, which is impossible. Therefore a circle where $\tan^2 z$ is an integer equal to or greater than 2 could only be equated to a suitable polygon, as Bolyai implies in the text.

For example, a circle where $\tan^2 z = 3$ can be equated to a regular hexagon. The Euclidean formula for the sum of the angles of a hexagon is $4\pi$, so its defect is $(4\pi - 6v)$, where $v$ represents the magnitude of its internal angle. So in the case of a hexagon equal to a circle whose area is three times the area of the maximum triangle, Bolyai's formula requires that

$$k^2(4\pi - 6v) = 3\pi k^2,$$

which implies that $v = \pi/6 = 30^\circ$. This hexagon is shown below:

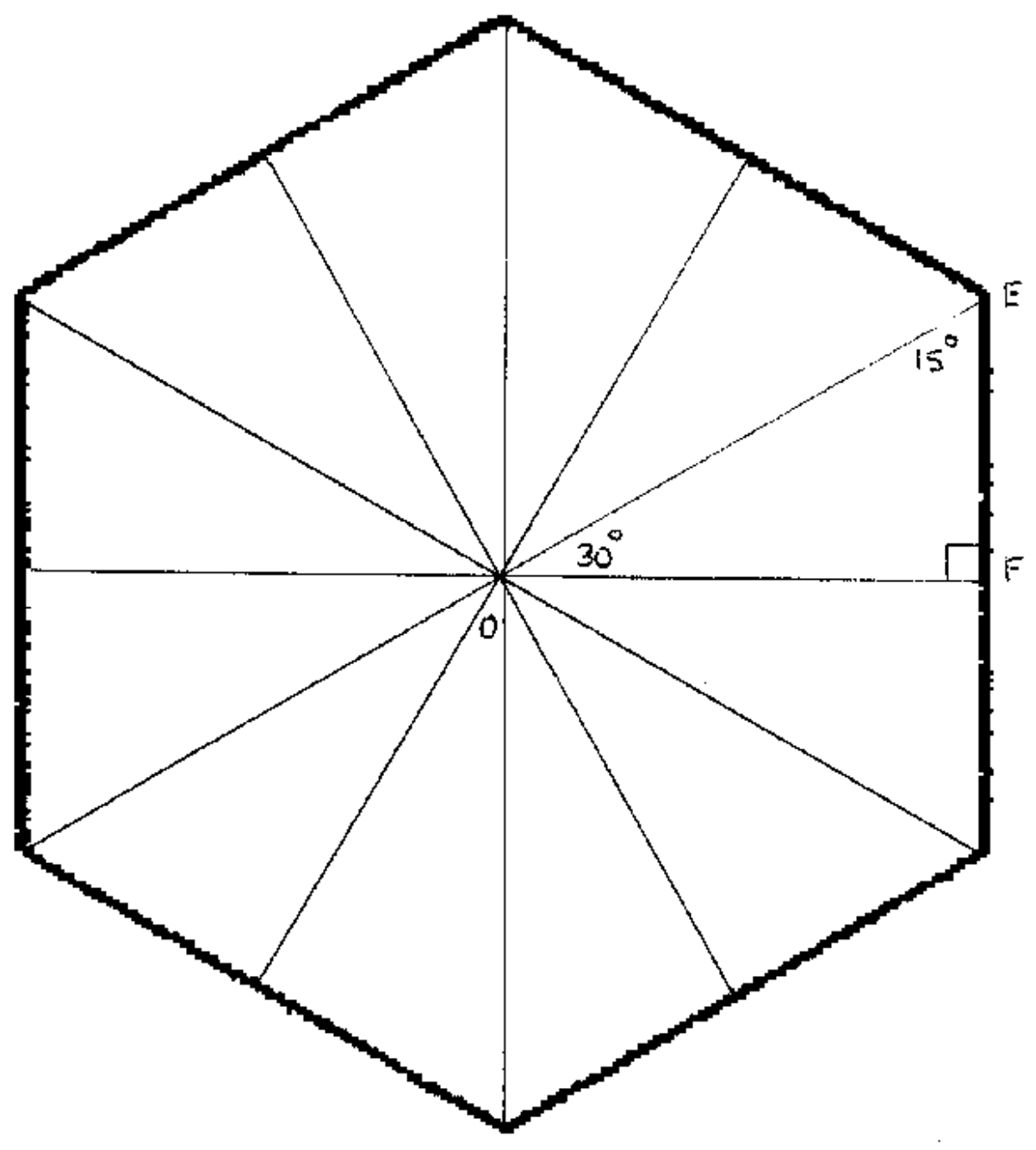

The quadrature requires the construction of the right angled triangle OEF, which represents 1/12 of the area of the hexagon. The acute angles of the triangle are 30° and 15°, which can be drawn without difficulty. The sides of the triangle can be constructed by compass and straightedge using the method outlined in the case of the square, after which the required hexagon can be drawn, completing the quadrature of the given circle.

The feasibility of quadrature in the second case, when $tan^2 z$ is a fraction, is determined by Gauss' theorem on constructible regular polygons. In the course of this theorem, which Bolyai praises as 'an outstanding discovery of the age, nay of all ages', Gauss shows that an angle $2\pi/k$ is only constructible by compass and straightedge if $k$ is a Fermat prime or a product of Fermat distinct Fermat primes times $2^m$. Subject to the condition that $tan^2 z$ is a rational number $j/k$ where $k$ is of the form stipulated by Gauss, the internal angle of the square or polygon can be constructed, after which the quadrature of the circle can be obtained.

The solution of this notorious problem using compass and straightedge, assuming the validity of hyperbolic geometry, remains an extraordinary intellectual feat, one of which Bolyai was obviously proud since he refers to it on the title page of the *Appendix.*

Bolyai ends the *Appendix* with a reference to **the impossibility of determining whether the true geometry is Euclidean or hyperbolic.** The first explicit proof that this question cannot be decided *a priori* was probably provided by Roberto Bonola in 1906. Henri Poincaré (1854-1912) had shown that the hyperbolic plane can be mapped onto an open Euclidean disc in which the straight lines of the plane are represented either by diameters or the arcs of circles orthogonal to the circumference of the disc:

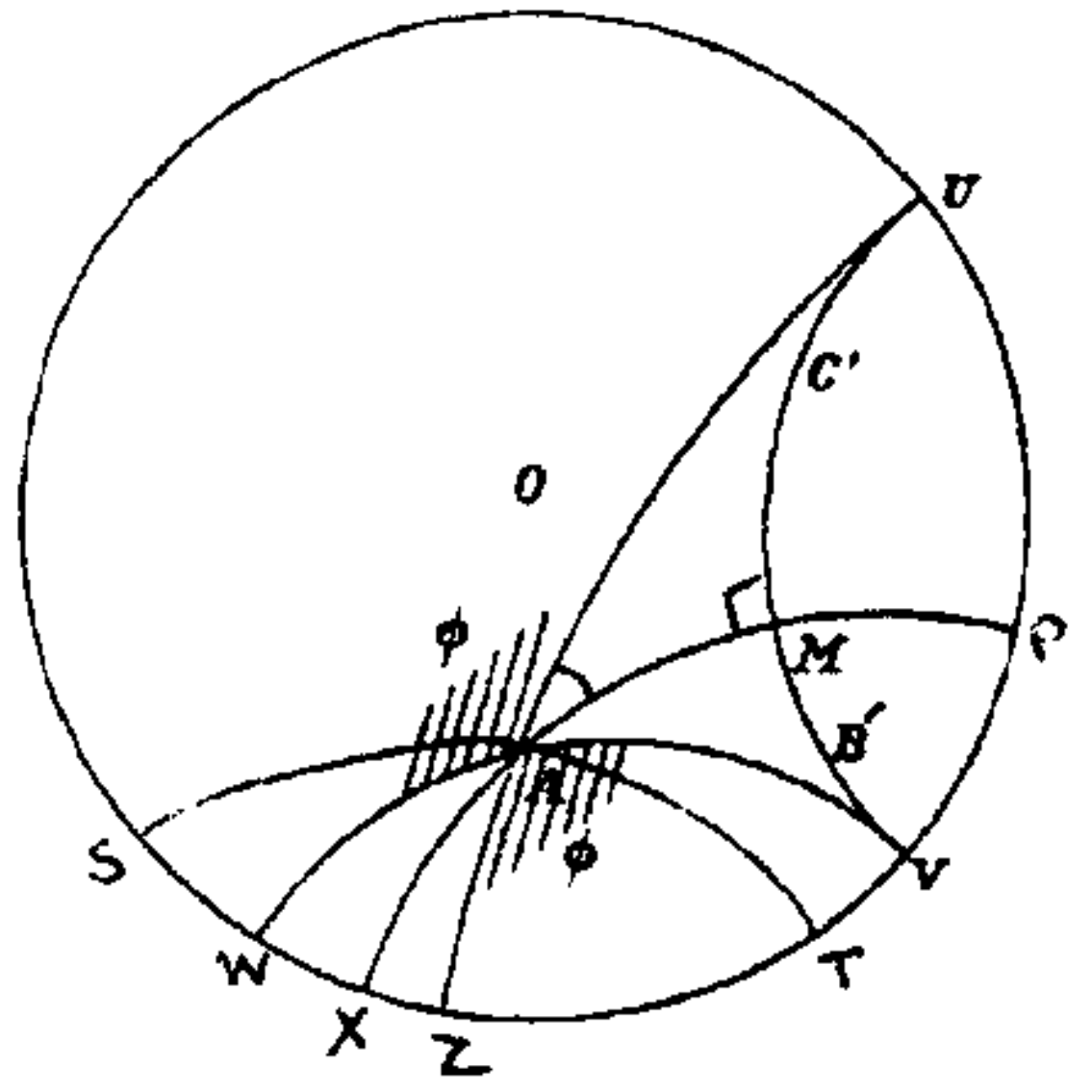

In order to map the whole of the hyperbolic plane onto the unit disc a suitable metric must be introduced to ensure that arcs such as XAMP are infinitely long. This can be achieved by defining thelength of a segment such as AM by the metric

AM = *log* (AP.XM /XA.MP),

where AP, XA, MP, XM are chords of the arc XAMP. This definition of length ensures that a segment AM becomes infinitely long as A approaches X or as M approaches P. In effect the circumference of the disc is deemed to lie at infinity. The definition also satisfies the usual properties of a distance function, among them the additive property, so that, given any point D (not shown) lying between A and M, AM = AD + DM. This follows from the law of logarithms, whereby

AM = *log* (AP.XM/XA.MP) = *log* {(AP.XD/XA.DP).(DP.XM/XD.MP)} = *log* (J.K) = *log* J + *log* K = AD + DM.

It can seen that all the postulates which are satisfied by straight lines in the hyperbolic plane apply equally to orthogonal circles on the Euclidean disc viz: 1) only one orthogonal arc can be drawn through any two points on the disc 2) an arc can be extended indefinitely in either direction 3) a (non-orthogonal) circle of any radius can be drawn on the disc 4) all right angles on the disc are equal and 5) the hyperbolic substitute for the parallel postulate applies to the orthogonal arcs. For given an arc UV and a point A not on this arc, the arcs WAV and ZAU drawn through A represent the left hand and right hand 'parallels' (so to speak) to UV, since each meets UV asymptotically at an infinitely distant point at the circumference. In contrast other arcs drawn through A are either intersecting with respect to UV, such as XAMP, or non-intersecting, such as SAT. Thus the 'parallels' WAV and ZAU separate the intersecting arcs from the non-intersecting arcs, just as parallels separate intersecting straight lines from non-intersecting straight lines in the hyperbolic plane.

It follows that all the properties of hyperbolic geometry will apply to figures on the disc. For example, the sum of the angles of any triangle on the disc formed by intersecting arcs, such as the doubly asymptotic triangle AUV, is less than two right angles. Furthermore in the asymptotic right angled triangle AMU, ∠UAM (marked) can be considered as the angle of parallelism corresponding to the segment AM. As a result it can be shown that the identity $\tan \frac{1}{2} \angle UAM = e^{-AM}$ and all the other trigonometric formulas of hyperbolic geometry hold true on the disc (where $k$ is taken as equal to 1). It is also possible to construct the equivalent of horocycles and equidistant lines on the disc.

Bonola pointed out that the possibility of representing the hyperbolic plane on a Euclidean disc implies that any contradiction in hyperbolic geometry would imply a contradiction in Euclidean geometry and vice versa. It is therefore impossible to establish the validity of one at the expense of the other. Thus the parallel postulate can neither be proved nor disproved. It is not a theorem. It is a postulate independent of the first four postulates, vindicating Euclid's decision to place it alongside the others.

17 July 2023